\documentclass[11pt, a4paper]{article}
\usepackage[cp1251]{inputenc}
\usepackage[russian]{babel}
\usepackage{amsmath} \usepackage{euscript}
\usepackage{mathrsfs}
\usepackage{longtable}
\usepackage{amssymb}
\usepackage[matrix, arrow, curve]{xy}
\begin{document}

\setcounter{MaxMatrixCols}{14}
\newcounter{bnomer} \newcounter{snomer} \newcounter{tnomer}
\newcounter{bsnomer}
\setcounter{bnomer}{0}
\setcounter{tnomer}{0}
\renewcommand{\thesnomer}{\thebnomer.\arabic{snomer}}
\renewcommand{\thebsnomer}{\thebnomer.\arabic{bsnomer}}
\renewcommand{\refname}{\begin{center}\large{\textbf{References}}\end{center}}

\newcommand{\sect}[1]{%
\setcounter{snomer}{0}\setcounter{bsnomer}{0}
\refstepcounter{bnomer}
\par\bigskip\begin{center}\large{\textbf{\arabic{bnomer}. {#1}}}\end{center}}
\newcommand{\sst}[1]{%
\refstepcounter{bsnomer}
\par\bigskip\textbf{\arabic{bnomer}.\arabic{bsnomer}. {#1}.}}
\newcommand{\defi}[1]{%
\refstepcounter{snomer}
\par\medskip\textbf{Definition \arabic{bnomer}.\arabic{snomer}. }{#1}\par\medskip}
\newcommand{\theo}[2]{%
\refstepcounter{snomer}
\par\textbf{Теорема \arabic{bnomer}.\arabic{snomer}. }{#2} {\emph{#1}}\hspace{\fill}$\square$\par}
\newcommand{\mtheop}[2]{%
\refstepcounter{snomer}
\par\textbf{Theorem \arabic{bnomer}.\arabic{snomer}. }{\emph{#1}}
\par\textsc{Proof}. {#2}\hspace{\fill}$\square$\par}
\newcommand{\mcorop}[2]{%
\refstepcounter{snomer}
\par\textbf{Corollary \arabic{bnomer}.\arabic{snomer}. }{\emph{#1}}
\par\textsc{Proof}. {#2}\hspace{\fill}$\square$\par}
\newcommand{\mtheo}[1]{%
\refstepcounter{snomer}
\par\medskip\textbf{Theorem \arabic{bnomer}.\arabic{snomer}. }{\emph{#1}}\par\medskip}
\newcommand{\mlemm}[1]{%
\refstepcounter{snomer}
\par\medskip\textbf{Lemma \arabic{bnomer}.\arabic{snomer}. }{\emph{#1}}\par\medskip}
\newcommand{\mprop}[1]{%
\refstepcounter{snomer}
\par\medskip\textbf{Proposition \arabic{bnomer}.\arabic{snomer}. }{\emph{#1}}\par\medskip}
\newcommand{\theobp}[2]{%
\refstepcounter{snomer}
\par\textbf{Теорема \arabic{bnomer}.\arabic{snomer}. }{#2} {\emph{#1}}\par}
\newcommand{\theop}[2]{%
\refstepcounter{snomer}
\par\textbf{Theorem \arabic{bnomer}.\arabic{snomer}. }{\emph{#1}}
\par\textsc{Proof}. #2\hspace{\fill}$\square$\par}
\newcommand{\exam}[1]{%
\refstepcounter{snomer}
\par\medskip\textbf{Example \arabic{bnomer}.\arabic{snomer}. }{#1}\par\medskip}
\newcommand{\deno}[1]{%
\refstepcounter{snomer}
\par\textbf{Definition \arabic{bnomer}.\arabic{snomer}. }{#1}\par}
\newcommand{\post}[1]{%
\refstepcounter{snomer}
\par\textbf{Предложение \arabic{bnomer}.\arabic{snomer}. }{\emph{#1}}\hspace{\fill}$\square$\par}
\newcommand{\postp}[2]{%
\refstepcounter{snomer}
\par\medskip\textbf{Proposition \arabic{bnomer}.\arabic{snomer}. }{\emph{#1}}%
\ifhmode\par\fi\textsc{Proof}. {#2}\hspace{\fill}$\square$\par\medskip}
\newcommand{\lemm}[1]{%
\refstepcounter{snomer}
\par\textbf{Lemma \arabic{bnomer}.\arabic{snomer}. }{\emph{#1}}\hspace{\fill}$\square$\par}
\newcommand{\lemmp}[2]{%
\refstepcounter{snomer}
\par\medskip\textbf{Lemma \arabic{bnomer}.\arabic{snomer}. }{\emph{#1}}
\par\textsc{Proof}. #2\hspace{\fill}$\square$\par\medskip}
\newcommand{\coro}[1]{%
\refstepcounter{snomer}
\par\textbf{Corollary \arabic{bnomer}.\arabic{snomer}. }{\emph{#1}}\hspace{\fill}$\square$\par}
\newcommand{\mcoro}[1]{%
\refstepcounter{snomer}
\par\textbf{Corollary \arabic{bnomer}.\arabic{snomer}. }{\emph{#1}}\par\medskip}
\newcommand{\corop}[2]{%
\refstepcounter{snomer}
\par\textbf{Corollary \arabic{bnomer}.\arabic{snomer}. }{\emph{#1}}
\par\textsc{Proof}. {#2}\hspace{\fill}$\square$\par}
\newcommand{\nota}[1]{%
\refstepcounter{snomer}
\par\medskip\textbf{Remark \arabic{bnomer}.\arabic{snomer}. }{#1}\par\medskip}
\newcommand{\propp}[2]{%
\refstepcounter{snomer}
\par\medskip\textbf{Proposition \arabic{bnomer}.\arabic{snomer}. }{\emph{#1}}
\par\textsc{Proof}. {#2}\hspace{\fill}$\square$\par\medskip}
\newcommand{\hypo}[1]{%
\refstepcounter{snomer}
\par\medskip\textbf{Conjecture \arabic{bnomer}.\arabic{snomer}. }{\emph{#1}}\par\medskip}
\newcommand{\prop}[1]{%
\refstepcounter{snomer}
\par\textbf{Proposition \arabic{bnomer}.\arabic{snomer}. }{\emph{#1}}\hspace{\fill}$\square$\par}
\newcommand{\tabl}[1]{%
\refstepcounter{tnomer}
\textbf{Table \arabic{tnomer}. {#1}}}

\newcommand{\Ind}[3]{%
\mathrm{Ind}_{#1}^{#2}{#3}}
\newcommand{\Res}[3]{%
\mathrm{Res}_{#1}^{#2}{#3}}
\newcommand{\epsi}{\epsilon}
\newcommand{\tri}{\triangleleft}
\newcommand{\Supp}[1]{%
\mathrm{Supp}(#1)}

\makeatletter
\def\iddots{\mathinner{\mkern1mu\raise\p@
\vbox{\kern7\p@\hbox{.}}\mkern2mu
\raise4\p@\hbox{.}\mkern2mu\raise7\p@\hbox{.}\mkern1mu}}
\makeatother

\newcommand{\reg}{\mathrm{reg}}
\newcommand{\empr}[2]{[-{#1},{#1}]\times[-{#2},{#2}]}
\newcommand{\sreg}{\mathrm{sreg}}
\newcommand{\codim}{\mathrm{codim}\,}
\newcommand{\chara}{\mathrm{char}\,}
\newcommand{\rk}{\mathrm{rk}\,}
\newcommand{\chr}{\mathrm{ch}\,}
\newcommand{\id}{\mathrm{id}}
\newcommand{\Ad}{\mathrm{Ad}}
\newcommand{\aad}{\mathrm{ad}}
\newcommand{\Ker}{\mathrm{Ker}\,}
\newcommand{\End}{\mathrm{End}}
\newcommand{\Ann}{\mathrm{Ann}\,}
\newcommand{\col}{\mathrm{col}}
\newcommand{\row}{\mathrm{row}}
\newcommand{\low}{\mathrm{low}}
\newcommand{\sur}{\twoheadrightarrow}
\newcommand{\inj}{\hookrightarrow}
\newcommand{\pho}{\hphantom{\quad}\vphantom{\mid}}
\newcommand{\fho}[1]{\vphantom{\mid}\setbox0\hbox{00}\hbox to \wd0{\hss\ensuremath{#1}\hss}}
\newcommand{\wt}{\widetilde}
\newcommand{\wh}{\widehat}
\newcommand{\ad}[1]{\mathrm{ad}_{#1}}
\newcommand{\tr}{\mathrm{tr}\,}
\newcommand{\GL}{\mathrm{GL}}
\newcommand{\Pf}{\mathrm{Pf}}
\newcommand{\Prim}{\mathrm{Prim}\,}
\newcommand{\Cent}{\mathrm{Cent}\,}
\newcommand{\SL}{\mathrm{SL}}
\newcommand{\SO}{\mathrm{SO}}
\newcommand{\UT}{\mathrm{UT}}
\newcommand{\Sp}{\mathrm{Sp}}
\newcommand{\Mat}{\mathrm{Mat}}
\newcommand{\Stab}{\mathrm{Stab}}
\newcommand{\ilm}{\varinjlim}
\newcommand{\lee}{\leq}
\newcommand{\gee}{\geq}

\newcommand{\vfi}{\varphi}
\newcommand{\vpi}{\varpi}
\newcommand{\teta}{\vartheta}
\newcommand{\Bfi}{\Phi}
\newcommand{\Fp}{\mathbb{F}}
\newcommand{\Rp}{\mathbb{R}}
\newcommand{\Zp}{\mathbb{Z}}
\newcommand{\Cp}{\mathbb{C}}
\newcommand{\Np}{\mathbb{N}}
\newcommand{\Xt}{\mathfrak{X}}
\newcommand{\ut}{\mathfrak{u}}
\newcommand{\at}{\mathfrak{a}}
\newcommand{\nt}{\mathfrak{n}}
\newcommand{\mt}{\mathfrak{m}}
\newcommand{\htt}{\mathfrak{h}}
\newcommand{\spt}{\mathfrak{sp}}
\newcommand{\slt}{\mathfrak{sl}}
\newcommand{\glt}{\mathfrak{gl}}
\newcommand{\ot}{\mathfrak{so}}
\newcommand{\rt}{\mathfrak{r}}
\newcommand{\hei}{\mathfrak{hei}}
\newcommand{\rad}{\mathfrak{rad}}
\newcommand{\bt}{\mathfrak{b}}
\newcommand{\gt}{\mathfrak{g}}
\newcommand{\vt}{\mathfrak{v}}
\newcommand{\pt}{\mathfrak{p}}
\newcommand{\kt}{\mathfrak{k}}
\newcommand{\Po}{\mathcal{P}}
\newcommand{\Uo}{\EuScript{U}}
\newcommand{\Fo}{\EuScript{F}}
\newcommand{\Do}{\EuScript{D}}
\newcommand{\Eo}{\EuScript{E}}
\newcommand{\Iu}{\mathcal{I}}
\newcommand{\Du}{\mathcal{D}}
\newcommand{\Uu}{\mathcal{U}}
\newcommand{\Mo}{\mathcal{M}}
\newcommand{\Nu}{\mathcal{N}}
\newcommand{\Ro}{\mathcal{R}}
\newcommand{\Co}{\mathcal{C}}
\newcommand{\Lo}{\mathcal{L}}
\newcommand{\Ou}{\mathcal{O}}
\newcommand{\Au}{\mathcal{A}}
\newcommand{\Vu}{\mathcal{V}}
\newcommand{\Bu}{\mathcal{B}}
\newcommand{\Sy}{\mathcal{Z}}
\newcommand{\Sb}{\mathcal{F}}
\newcommand{\Gr}{\mathcal{G}}
\newcommand{\rtc}[1]{C_{#1}^{\mathrm{red}}}

\author{Mikhail V. Ignatyev \and Aleksandr A. Shevchenko}

\date{}
\title{\Large{Centrally generated primitive ideals of $U(\nt)$ for exceptional types{$\vphantom{1}$}\footnotetext{The work was supported by the Foundation for the Advancement of Theoretical Physics and Mathematics ``BASIS'', grant no. 18--1--7--2--1.}}
}\maketitle

\begin{center}
\begin{tabular}{p{15cm}}
\small{\textsc{Abstract}. Let $\gt$ be a complex semisimple Lie algebra, let $\bt$ be a Borel subalgebra of $\gt$, let $\nt$ be the nilradical of $\bt$, and let $U(\nt)$ be the universal enveloping algebra of $\nt$. We study primitive ideals of~$U(\nt)$. Almost all primitive ideals are centrally generated, i.e., are generated by their intersections with the center $Z(\nt)$ of $U(\nt)$. We present an explicit characterization of the centrally generated primitive ideals of $U(\nt)$ in terms of the Dixmier map and the Kostant cascade in the case when $\gt$ is a simple algebra of exceptional type. (For classical simple Lie algebras, a similar characterization was obtained by Ivan Penkov and the first author.) As a corollary, we establish a classification of centrally generated primitive ideals of $U(\nt)$ for an arbitrary semisimple algebra $\gt$.}
\end{tabular}
\end{center}

\sect{Introduction}

The ground field is the field $\Cp$ of complex numbers. The theory of primitive ideals in enveloping algebras of Lie algebras can be considered as a part of the re\-pre\-sen\-tation theory of Lie algebras. In general, it is impossible to classify all irreducible representations of Lie algebras (except in a few very special cases), while a classification of annihilators of irreducible representations, i.e., of primitive ideals, can be achieved in much greater generality. In the case when $\nt$ is a finite-dimensional nilpotent Lie algebra, the primitive ideals in the universal enveloping algebra $U(\nt)$ can be described in terms of the Dixmier map assigning to any linear form $f\in\nt^*$ a primitive ideal $J(f)$ of $U(\nt)$.

If $\nt$ is abelian, $J(f)$ is simply the annihilator of $f$. For a general finite-dimensional nilpotent Lie algebra $\nt$, the theory of primitive ideals retains many properties from the abelian case: in particular, $J(f)$ is always a maximal ideal and every primitive ideal in $U(\nt)$ is of the form $J(f)$ for some $f\in\nt^*$. Moreover, $J(f)=J(f')$ if and only if $f$ and $f'$ belong to the same coadjoint orbit in $\nt^*$ \cite[Theo\-rem~6.2.4]{Dixmier3}. Note that, according to the Kostant--Kirillov orbit method, the set of all unitary irreducible representations of the Lie group $\exp\nt$ is also in one-to-one correspondence with the set of all coadjoint orbits in $\nt^*$ \cite{Kirillov1}. Note also that there exists a natural topology, called the Jacobson topology, on the set $\Prim U(\nt)$ of all primitive ideals of $U(\nt)$. It turns out that $\Prim U(\nt)$ is irreducible, and each primitive ideal from a certain (dense) open subset of $\Prim U(\nt)$ is generated by its intersection with the center $Z(\nt)$ of $U(\nt)$; such an ideal is called centrally generated (see Section~\ref{sect:cent_gen_ideals} for the details).

Suppose that $\nt$ is the nilradical of a Borel subalgebra $\bt$ of a complex finite-dimensional semisimple Lie algebra $\gt$. The description of $Z(\nt)$ goes back to Dixmier, Joseph and Kostant. It turns out that $Z(\nt)$ is a polynomial algebra whose generators are parametrized by the positive roots from the Kostant cascade $\Bu$, a certain strongly orthogonal subset of the set $\Phi^+$ of positive roots of $\gt$ with respect to $\bt$ \cite[Theorem 7]{Kostant2} (see Section~\ref{sect:center} for the precise definition). Let $\{e_{\alpha},~\alpha\in\Phi^+\}$ be a basis of $\nt$ consisting of root vectors. We say that a linear form $f\in\nt^*$ is a Kostant form if $f(e_{\alpha})=0$ for $\alpha\notin\Bu$ and $f(e_{\beta})\neq0$ for $\beta\in\Bu\setminus\Delta$, where $\Delta$ is the set of simple roots. Note that the coadjoint orbit of a Kostant form has maximal possible dimension.

It was proved in \cite[Theorem~3.1]{IgnatyevPenkov1} and \cite[Theorem 2.4]{Ignatyev2} that, when $\Phi$ is of classical type (i.e., $\Phi=A_{n-1}$, $B_n$, $C_n$ or $D_n$), $J\in\Prim U(\nt)$ is centrally generated if and only if $J=J(f)$ for a certain Kostant form~$f$. In this paper, we prove that this fact is also true when $\Phi$ is of exceptional type, i.e., $\Phi=E_6$, $E_7$, $E_8$, $F_4$ or $G_2$. Namely, let $\Delta_{\beta}$, $\beta\in\Bu$, be the set of canonical generators of $Z(\nt)$ (see Section~\ref{sect:center}). Let $J$ be a primitive ideal of $U(\nt)$. Since $J$ is the annihilator of a~simple $\nt$-module, given $\beta\in\Bu$, there exists unique $c_{\beta}\in\Cp$ such that $\Delta_{\beta}-c_{\beta}\in J$. Our main result, Theorem~\ref{theo:ideal_exceptional}, claims that the following conditions are equivalent:
\begin{equation*}
\begin{split}
&\text{\textup{i)} $J$ is centrally generated};\\
&\text{\textup{ii)} all scalars $c_{\beta}$, $\beta\in\Bu\setminus\Delta$, are nonzero};\\
&\text{\textup{iii)} $J=J(f)$ for a Kostant form $f$}.
\end{split}
\end{equation*}
If these conditions are satisfied, then $f$ can be explicitly reconstructed from $J$. As a corollary, we conclude that the same is true for arbitrary root system.

The paper is organized as follows. In Sec\-tion~\ref{sect:center}, we briefly recall Kostant's characterization of $Z(\nt)$ and present a (more or less) explicit description of the canonical generators of $Z(\nt)$ based on Panov's work \cite{Panov2}. Using this description, in Section~\ref{sect:cent_gen_ideals} we prove that certain centrally generated ideals are primitive (in fact, it is the key ingredient in the proof of the main result, see~Proposition~\ref{prop:J_c_primitive_exceptional}). Section~\ref{sect:distinct_orbits} is devoted to some particular classes of coadjoint orbits. Namely, we prove that certain orbits are disjoint, see Proposition~\ref{prop:distinc_orbits_precise}. Finally, in Section~\ref{sect:proof_main}, combining our results from two previous sections, we prove the main result, Theorem~\ref{theo:ideal_exceptional}. As an immediate corollary, we obtain that the similar result is true for an arbitrary semisimple Lie algebra, see Theorem~\ref{theo:ideal_all}.

\textsc{Acknowledgments}. We express our gratitude to A. Panov and I. Penkov for useful discussions. The work was supported by the Foundation for the Advancement of Theoretical Physics and Mathema\-tics ``BASIS'', grant no. 18--1--7--2--1.

\sect{The center of $U(\nt)$}

Let $G$ be\label{sect:center} a complex semisimple algebraic group, let $H$ be a Cartan subgroup of $G$, let $B$ be a Borel subgroup of $G$ containing $H$, and let $N$ be the unipotent radical of $B$. We denote by $\Phi$ the root system of $G$ with respect to $B$, and by $\Phi^+$ the set of positive roots with respect to $B$. Let $\gt$ (respectively, $\htt$, $\bt$ and $\nt$) be the Lie algebra of $G$ (respectively, of $H$, $B$ and $N$), so that $\bt=\htt\oplus\nt$ as vector spaces. The Lie algebra $\nt$ has a basis consisting of root vectors $e_{\alpha}$, $\alpha\in\Phi^+$. We denote the dual basis of the dual space $\nt^*$ by $\{e_{\alpha}^*,~\alpha\in\Phi^+\}$.

Let $\Rp^n$ be the $n$-dimensional Euclidean space with the standard inner product $(\cdot,\cdot)$, and $\{\epsi_i\}_{i=1}^n$ be the standard basis of $\Rp^n$. If $\Phi$ is irreducible then we identify $\Phi^+$ with the following subset of $\Rp^n$ \cite{Bourbaki1}:
\begin{equation*}\predisplaypenalty=0
\begin{split}
A_{n-1}^+&=\{\epsi_i-\epsi_j,~1\leq i<j\leq n\},\\
B_n^+&=\{\epsi_i-\epsi_j,~1\leq i<j\leq n\}\cup\{\epsi_i+\epsi_j,~1\leq i<j\leq n\}\cup\{\epsi_i,~1\leq i\leq n\},\\
C_n^+&=\{\epsi_i-\epsi_j,~1\leq i<j\leq n\}\cup\{\epsi_i+\epsi_j,~1\leq i<j\leq n\}\cup\{2\epsi_i,~1\leq i\leq n\},\\
D_n^+&=\{\epsi_i-\epsi_j,~1\leq i<j\leq n\}\cup\{\epsi_i+\epsi_j,~1\leq i<j\leq n\},\\
E_6^+&=\{\pm\epsi_i+\epsi_j,~1\leq i<j\leq5\}\cup\left\{\dfrac{1}{2}\left(\epsi_8-\epsi_7-\epsi_6+\sum\limits_{i=1}^5(-1)^{\nu(i)}\epsi_i\right),~\sum\limits_{i=1}^5\nu(i)\text{ is even}\right\},\\
E_7^+&=\{\pm\epsi_i+\epsi_j,~1\leq i<j\leq6\}\cup\{\epsi_8-\epsi_7\}\cup\left\{\dfrac{1}{2}\left(\epsi_8-\epsi_7+\sum\limits_{i=1}^{6}(-1)^{\nu(i)}\epsi_i \right),~\sum\limits_{i=1}^{6}\nu(i)\text{ is even}\right\},\\
E_8^+&=\{\pm\epsi_i+\epsi_j,1\leq i<j\leq8\}\cup\left\{\dfrac{1}{2}\left(\epsi_8+\sum\limits_{i=1}^{7}(-1)^{\nu(i)}\epsi_i\right),~\sum\limits_{i=1}^{7}\nu(i)\text{ is even}\right\},\\
F_4^+&=\{\epsi_i\pm\epsi_j,~1\leq i<j\leq4\}\cup\{(\epsi_1\pm\epsi_2\pm\epsi_3\pm\epsi_4)/2\}\cup\{\epsi_i,1\leq i\leq 4\},\\
G_2^+&=\{\epsi_1-\epsi_2,~-2\epsi_1+\epsi_2+\epsi_3,~-\epsi_1+\epsi_3,-\epsi_2+\epsi_3,~\epsi_1-2\epsi_2+\epsi_3,~-\epsi_1-\epsi_2+2\epsi_3\}.\\
\end{split}
\end{equation*}

Under this identification, the set $\Delta\subset\Phi^+$ of the simple roots has the following form:
\begin{equation}\label{formula:simple_roots}
\Delta=\begin{cases}\bigcup\nolimits_{i=1}^{n-1}\{\alpha_i=\epsi_i-\epsi_{i+1}\}&\text{for }A_{n-1},\\
\bigcup\nolimits_{i=1}^{n-1}\{\alpha_i=\epsi_i-\epsi_{i+1}\}\cup\{\alpha_n=\epsi_n\}&\text{for }B_n,\\
\bigcup\nolimits_{i=1}^{n-1}\{\alpha_i=\epsi_i-\epsi_{i+1}\}\cup\{\alpha_n=2\epsi_n\}&\text{for }C_n,\\
\bigcup\nolimits_{i=1}^{n-1}\{\alpha_i=\epsi_i-\epsi_{i+1}\}\cup\{\alpha_n=\epsi_{n-1}+\epsi_n\}&\text{for }D_n,\\
\{\alpha_1=(\epsi_1+\epsi_8-\sum\nolimits_{k=2}^7\epsi_k)/2,\\
\alpha_2=\epsi_1+\epsi_2\}\cup\bigcup\nolimits_{i=1}^4\{\alpha_{i+2}=\epsi_{i+1}-\epsi_i\}&\text{for }E_6,\\
\{\alpha_1=(\epsi_1+\epsi_8-\sum\nolimits_{k=2}^7\epsi_k)/2,\\
\alpha_2=\epsi_1+\epsi_2\}\cup\bigcup\nolimits_{i=1}^5\{\alpha_{i+1}=\epsi_{i+1}-\epsi_i\}&\text{for }E_7,\\
\{\alpha_1=(\epsi_1+\epsi_8-\sum\nolimits_{k=2}^7\epsi_k)/2,\\
\alpha_2=\epsi_1+\epsi_2\}\cup\bigcup\nolimits_{i=1}^6\{\alpha_{i+2}=\epsi_{i+1}-\epsi_i\}&\text{for }E_8,\\
\{\alpha_1=\epsi_2-\epsi_3,~\alpha_2=\epsi_3-\epsi_4,~\alpha_3=\epsi_4,\\
\alpha_4=(\epsi_1-\epsi_2-\epsi_3-\epsi_4)/2\}&\text{for }F_4,\\
\{\alpha_1=\epsi_1-\epsi_2,~\alpha_2=-2\epsi_1+\epsi_2+\epsi_3\}&\text{for }G_2.\\
\end{cases}
\end{equation}

Recall that there exists a natural partial order on $\Phi$: by definition, $\alpha<\beta$ if $\beta-\alpha$ can be represented as a sum of positive roots. Denote by $\Bu$ the subset of $\Phi^+$ constructed by the following inductive procedure. Let $\Bu_1$ be the set consisting of the maximal roots of all irreducible components of $\Phi$. For $n\geq2$, we denote $\Phi_n=\{\alpha\in\Phi\mid\alpha\perp\beta\text{ for all }\beta\in\Bu_1\cup\ldots\cup\Bu_{n-1}\}$, and set $\Bu_n$ to be the set of the maximal roots of all irreducible components of $\Phi_n$. Finally, we denote by~$\Bu$ the union of all $\Bu_n$'s. Note that $\Bu$ is a maximal strongly orthogonal subset of $\Phi^+$, i.e., $\Bu$ is maximal with the property that if $\alpha,~\beta\in\Bu$ then neither $\alpha-\beta$ nor $\alpha+\beta$ belongs to $\Phi^+$.

\defi{We call $\Bu$ the \emph{Kostant cascade} of orthogonal roots in~$\Phi^+$.}

Let $[n/2]$ be the largest integer not exceeding $n/2$. For irreducible $\Phi$, $\Bu$ has the following form:
\begin{equation*}
\Bu=\begin{cases}\bigcup\nolimits_{i=1}^{[n/2]}\{\beta_i=\epsi_i-\epsi_{n-i+1}\}&\text{for }A_{n-1},\\
\bigcup\nolimits_{i=1}^{n/2}\{\beta_{2i-1}=\epsi_{2i-1}+\epsi_{2i},~\beta_{2i}=\epsi_{2i-1}-\epsi_{2i}\}&\text{for }B_n,~n\text{ even},\\
\bigcup\nolimits_{i=1}^{[n/2]}\{\beta_{2i-1}=\epsi_{2i-1}+\epsi_{2i},~\beta_{2i}=\epsi_{2i-1}-\epsi_{2i}\}\cup\{\beta_n=\epsi_n\}&\text{for }B_n,~n\text{ odd},\\
\bigcup\nolimits_{i=1}^n\{\beta_i=2\epsi_i\}&\text{for }C_n,\\
\bigcup\nolimits_{i=1}^{[n/2]}\{\beta_{2i-1}=\epsi_{2i-1}+\epsi_{2i},~\beta_{2i}=\epsi_{2i-1}-\epsi_{2i}\}&\text{for }D_n,\\
\{\beta_1=(\epsi_1+\epsi_2+\epsi_3+\epsi_4+\epsi_5-\epsi_6-\epsi_7+\epsi_8)/2,\\\beta_2=(-\epsi_1-\epsi_2-\epsi_3-\epsi_4+\epsi_5-\epsi_6-\epsi_7+\epsi_8)/2,\\\beta_3=
-\epsi_1+\epsi_4,~\beta_4=-\epsi_2+\epsi_3\}&\text{for }E_6,\\
\{\beta_1=-\epsi_7+\epsi_8,~\beta_2=\epsi_5+\epsi_6,~\beta_3=\epsi_3+\epsi_4,\beta_4=-\epsi_5+\epsi_6,\\\beta_5=\epsi_1+\epsi_2,\beta_6=-\epsi_1+\epsi_2,~\beta_7=-\epsi_3+\epsi_4\}&\text{for }E_7,\\
\{\beta_1=\epsi_7+\epsi_8,~\beta_2=-\epsi_7+\epsi_8,~\beta_3=\epsi_5+\epsi_6,~\beta_4=\epsi_3+\epsi_4,\\~\beta_5=-\epsi_5+\epsi_6,~\beta_6=\epsi_1+\epsi_2,~\beta_7=-\epsi_1+\epsi_2,~\beta_8=-\epsi_3+\epsi_4\}&\text{for }E_8,\\
\{\beta_1=\epsi_1+\epsi_2,~\beta_2=\epsi_1-\epsi_2,~\beta_3=\epsi_3+\epsi_4,~\beta_4=\epsi_3-\epsi_4\}&\text{for }F_4,\\
\{\beta_1=-\epsi_1-\epsi_2+2\epsi_3,~\beta_2=\epsi_1-\epsi_2\}&\text{for }G_2.\\
\end{cases}
\end{equation*}

Denote by $U(\nt)$ the enveloping algebra of $\nt$, and by $S(\nt)$ the symmetric algebra of $\nt$. Then $\nt$ and $S(\nt)$ are $B$-modules as $B$ normalizes $N$. Denote by $Z(\nt)$ the center of $U(\nt)$. It is well-known that the restriction of the symmetrization map
\begin{equation}
\sigma\colon S(\nt)\to U(\nt),~x^k\mapsto x^k,~x\in\nt,~k\in\Zp_{\geq0},\label{formula:symm_map}
\end{equation}
to the algebra $Y(\nt)$ of $N$-invariants in $S(\nt)$ is an algebra isomorphism between $Y(\nt)$ and $Z(\nt)$.\newpage

We next present a canonical set of generators of $Z(\nt)$ (or, equivalently, of $Y(\nt)$), whose description goes back to Dixmier, Joseph and Kostant \cite{Dixmier1}, \cite{Joseph1}, \cite{Kostant1}, \cite{Kostant2}. Let $\Xt$ be the group of rational multiplicative characters of $H$, i.e., of algebraic group morphisms from $H$ to $\Cp^{\times}=\Cp\setminus\{0\}$. We can consider $\Zp\Phi$, the $\Zp$-linear span of $\Phi$, as a subgroup of $\Xt$. Recall that a vector $\lambda\in\Rp^n$ is called a \emph{weight} if $c(\alpha,\lambda):=2(\alpha,\lambda)/(\alpha,\alpha)$ is an integer for any $\alpha\in\Phi^+$. A weight $\lambda$ is called \emph{dominant} if $c(\alpha,\lambda)\geq0$ for all $\alpha\in\Phi^+$. An element $a$ of an $H$-module is called a \emph{weight vector} if there exists $\nu\in\Xt$ such that $h\cdot a=\nu(h)a$ for all $h\in H$. By \cite[Theorems 6, 7]{Kostant2}, every weight occurs in $Y(\nt)$ with multiplicity at most 1. Furthermore, there exist unique (up to scalars) prime polynomials $\xi_{\beta}\in Y(\nt)$, $\beta\in\Bu$, such that each $\xi_{\beta}$ is a weight polynomial of a dominant weight $\mu_{\beta}$ belonging to the $\Zp$-linear span $\Zp\Bu$ of $\Bu$. A remarkable fact is that
\begin{equation}
\xi_{\beta},~\beta\in\Bu,\text{ are algebraically independent generators of }Y(\nt),\label{formula:xi_i}
\end{equation}
so $Y(\nt)$ and $Z(\nt)$ are polynomial rings \cite[Theorem 7]{Kostant2}. It turns out that the weights $\mu_{\beta}$'s have the following form \cite[Theorem 2.12]{Panov2}.

\begin{center}
\tabl{Weights $\mu_{\beta}$ for $\beta\in\Bu$}\label{table:weights_of_xi_i}
\begin{tabular}{|l|l|}
\hline
$\Phi=A_{n-1}$&$\mu_{\beta_i}=\epsi_1+\ldots+\epsi_i-\epsi_{n-i+1}-\ldots-\epsi_n$,~$1\leq i\leq[n/2]$\\
\hline
$\Phi=B_n$&$\mu_{\beta_i}=\begin{cases}2\epsi_1+\ldots+2\epsi_{i-1}&\text{for even $i$},\\
\epsi_1+\ldots+\epsi_i+\epsi_{i+1}&\text{for odd $i<n$,}\\
\epsi_1+\ldots+\epsi_i&\text{for odd $i=n$,}
\end{cases}$\\
\hline
$\Phi=C_n$&$\mu_{\beta_i}=2\epsi_1+\ldots+2\epsi_i$,~$1\leq i\leq n$\\
\hline
$\Phi=D_n$&$\mu_{\beta_i}=\begin{cases}2\epsi_1+\ldots+2\epsi_{i-1}&\text{for even $i<n$},\\
\epsi_1+\ldots+\epsi_{n-1}-\epsi_n&\text{for even $i=n$},\\
\epsi_1+\ldots+\epsi_{i+1}&\text{for odd $i$}
\end{cases}$\\
\hline
$\Phi=E_6$&
\begin{tabular}{l}
$\mu_{\beta_1}=(\epsilon_{1} + \epsilon_{2} + \epsilon_{3} + \epsilon_{4} + \epsilon_{5} - \epsilon_{6} - \epsilon_{7} + \epsilon_{8})/2,$\\
$\mu_{\beta_2}=\epsilon_{5} - \epsilon_{6} - \epsilon_{7} + \epsilon_{8},$\\
$\mu_{\beta_3}=(-\epsilon_{1} + \epsilon_{2} + \epsilon_{3} + 3\epsilon_{4} + 3\epsilon_{5} - 3\epsilon_{6} - 3\epsilon_{7} + 3\epsilon_{8})/2,$\\
$\mu_{\beta_4}=2\epsilon_{3} + 2\epsilon_{4} + 2\epsilon_{5} - 2\epsilon_{6} - 2\epsilon_{7} + 2\epsilon_{8}$
\end{tabular}
\\
\hline
$\Phi=E_7$&
\begin{tabular}{l}
$\mu_{\beta_1}=-\epsilon_{7} + \epsilon_{8}$\\
$\mu_{\beta_2}=\epsilon_{5} + \epsilon_{6} - \epsilon_{7} + \epsilon_{8},$\\
$\mu_{\beta_3}=\epsilon_{3} + \epsilon_{4} + \epsilon_{5} + \epsilon_{6} - 2\epsilon_{7} + 2\epsilon_{8},$\\
$\mu_{\beta_4}=2\epsilon_{6} - \epsilon_{7} + \epsilon_{8},$\\
$\mu_{\beta_4}=\epsilon_{1} + \epsilon_{2} + \epsilon_{3} + \epsilon_{5} + \epsilon_{5} + \epsilon_{6} - 2\epsilon_{7} + 2\epsilon_{8},$\\
$\mu_{\beta_5}=-\epsilon_{1} + \epsilon_{2} + \epsilon_{3} + \epsilon_{6} + \epsilon_{5} + \epsilon_{6} - 3\epsilon_{7} + 3\epsilon_{8},$\\
$\mu_{\beta_7}=2\epsilon_{4} + 2\epsilon_{5} + 2\epsilon_{6} - 3\epsilon_{7} + 3\epsilon_{8}$\\
\end{tabular}\\
\hline
$\Phi=E_8$&
\begin{tabular}{l}
$\mu_{\beta_1}=\epsilon_{7} + \epsilon_{8},~\mu_{\beta_2}=2\epsilon_{8},$\\
$\mu_{\beta_3}=\epsilon_{5} + \epsilon_{6} + \epsilon_{7} + 3\epsilon_{8}$,\\
$\mu_{\beta_4}=\epsilon_{3} + \epsilon_{4} + \epsilon_{5} + \epsilon_{6} + \epsilon_{7} + 5\epsilon_{8},$\\
$\mu_{\beta_5}=2\epsilon_{6} + 2\epsilon_{7} + 4\epsilon_{8},$\\
$\mu_{\beta_6}=\epsilon_{1} + \epsilon_{2} + \epsilon_{3} + \epsilon_{4} + \epsilon_{5} + \epsilon_{6} + \epsilon_{7} + 5\epsilon_{8},$\\
$\mu_{\beta_7}=-\epsilon_{1} + \epsilon_{2} + \epsilon_{3} + \epsilon_{4} + \epsilon_{5} + \epsilon_{6} + \epsilon_{7} + 7\epsilon_{8},$\\
$\mu_{\beta_8}=2\epsilon_{4} + 2\epsilon_{5} + 2\epsilon_{6} + 2\epsilon_{7} + 8\epsilon_{8}$\\
\end{tabular}\\
\hline
$\Phi=F_4$&
\begin{tabular}{l}
$\mu_{\beta_1}=\epsilon_{1} + \epsilon_{2},~\mu_{\beta_2}=2\epsilon_{1},$\\
$\mu_{\beta_3}=3\epsilon_{1} + \epsilon_{2} + \epsilon_{3} + \epsilon_{4},$\\
$\mu_{\beta_4}=4\epsilon_{1} + 2\epsilon_{2} + 2\epsilon_{3}$
\end{tabular}\\
\hline
$\Phi=G_2$&$\mu_{\beta_1}=-\epsilon_{1}-\epsilon_{2}+2\epsi_3,~\mu_{\beta_2}=-2\epsi_2+2\epsi_3$\\
\hline
\end{tabular}
\end{center}

\newpage

For the sequel, we need to express the weights $\mu_{\beta}$'s as linear combination of simple roots. Such expressions are presented in Table \ref{table:weights_of_xi_i_basis}.

\begin{center}
\tabl{Weights $\mu_{\beta}$, $\beta\in\Bu$, as linear combinations of simple roots}\label{table:weights_of_xi_i_basis}

\vspace{0.5cm}
\begin{tabular}{|l|l|}
\hline
$\Phi=A_{n-1}$&$\mu_{\beta_i}=\sum\nolimits_{k=1}^{i-1}k\alpha_k+\sum\nolimits_{k=i}^{n-i}i\alpha_k+\sum\nolimits_{k=n-i+1}^{n-1}(n-k)\alpha_k$,~$1\leq i\leq[n/2]\vphantom{1_1}$\\&\\
\hline
$\Phi=B_n$&$\mu_{\beta_i}=\begin{cases}2\sum\nolimits_{k=1}^{i-1}k\alpha_k+2(i-1)\sum\nolimits_{k=i}^n\alpha_k&\text{for even $i$},\\
\sum\nolimits_{k=1}^ik\alpha_k+(i+1)\sum\nolimits_{k=i+1}^n\alpha_k&\text{for odd $i$}
\end{cases}$\\&\\
\hline
$\Phi=C_n$&$\mu_{\beta_i}=\begin{cases}2\sum_{k=1}^ik\alpha_k+2i\sum_{k=i+1}^{n-1}\alpha_k+i\alpha_n&\text{for $i<n$},\\
2\sum_{k=1}^{n-1}k\alpha_k+n\alpha_n&\text{for $i=n$}
\end{cases}$\\&\\
\hline
$\Phi=D_n$&$\mu_{\beta_i}=\begin{cases}2\sum_{k=1}^{i-1}k\alpha_k+2(i-1)\sum_{k=i}^{n-2}\alpha_k+i(\alpha_{n-1}+\alpha_n)&\text{for even $i<n$},\\
\sum_{k=1}^{n-2}k\alpha_k+n\alpha_{n-1}/2+(n-2)\alpha_n/2&\text{for even $i=n$},\\
\sum_{k=1}^ik\alpha_k+(i+1)\sum_{k=i+1}^{n-2}\alpha_k+(i+1)(\alpha_{n-1}+\alpha_n)/2&\text{for odd $i<n-1$,}\\
\sum_{k=1}^{n-2}k\alpha_k+(n-2)\alpha_{n-1}/2+n\alpha_n/2&\text{for odd $i=n-1$}\\
\end{cases}$\\&\\
\hline
$\Phi=E_6$&
\begin{tabular}{l}
$\mu_{\beta_1}=\alpha_1+2\alpha_2+2\alpha_3+3\alpha_4+2\alpha_5+\alpha_6,$\\
$\mu_{\beta_2}=2\alpha_1+2\alpha_2+3\alpha_3+4\alpha_4+3\alpha_5+2\alpha_6,$\\
$\mu_{\beta_3}=3\alpha_1+4\alpha_2+6\alpha_3+8\alpha_4+6\alpha_5+3\alpha_6,$\\
$\mu_{\beta_4}=4\alpha_1+6\alpha_2+8\alpha_3+12\alpha_4+8\alpha_5+4\alpha_6$\\\\
\end{tabular}\\
\hline
$\Phi=E_7$&
\begin{tabular}{l}
$\mu_{\beta_1}=2\alpha_1+2\alpha_2+3\alpha_3+4\alpha_4+3\alpha_5+2\alpha_6+\alpha_7,$\\
$\mu_{\beta_2}=2\alpha_1+3\alpha_2+4\alpha_3+6\alpha_4+5\alpha_5+4\alpha_6+2\alpha_7,$\\
$\mu_{\beta_3}=4\alpha_1+6\alpha_2+8\alpha_3+12\alpha_4+9\alpha_5+6\alpha_6+3\alpha_7,$\\
$\mu_{\beta_4}=2\alpha_1+3\alpha_2+4\alpha_3+6\alpha_4+5\alpha_5+4\alpha_6+3\alpha_7,$\\
$\mu_{\beta_5}=4\alpha_1+7\alpha_2+8\alpha_3+12\alpha_4+9\alpha_5+6\alpha_6+3\alpha_7,$\\
$\mu_{\beta_6}=6\alpha_1+8\alpha_2+12\alpha_3+16\alpha_4+12\alpha_5+8\alpha_6+4\alpha_7,$\\
$\mu_{\beta_7}=6\alpha_1+9\alpha_2+12\alpha_3+18\alpha_4+15\alpha_5+10\alpha_6+5\alpha_7$\\\\
\end{tabular}\\
\hline
$\Phi=E_8$&
\begin{tabular}{l}
$\mu_{\beta_1}=2\alpha_1+3\alpha_2+4\alpha_3+6\alpha_4+5\alpha_5+4\alpha_6+3\alpha_7+2\alpha_8,$\\
$\mu_{\beta_2}=4\alpha_1+5\alpha_2+7\alpha_3+10\alpha_4+8\alpha_5+6\alpha_6+4\alpha_7+2\alpha_8,$\\
$\mu_{\beta_3}=6\alpha_1+9\alpha_2+12\alpha_3+18\alpha_4+15\alpha_5+12\alpha_6+8\alpha_7+4\alpha_8,$\\
$\mu_{\beta_4}=10\alpha_1+15\alpha_2+20\alpha_3+30\alpha_4+24\alpha_5+18\alpha_6+12\alpha_7+6\alpha_8,$\\
$\mu_{\beta_5}=8\alpha_1+12\alpha_2+16\alpha_3+24\alpha_4+20\alpha_5+16\alpha_6+12\alpha_7+6\alpha_8,$\\
$\mu_{\beta_6}=10\alpha_1+16\alpha_2+20\alpha_3+30\alpha_4+24\alpha_5+18\alpha_6+12\alpha_7+6\alpha_8,$\\
$\mu_{\beta_7}=14\alpha_1+20\alpha_2+28\alpha_3+40\alpha_4+32\alpha_5+24\alpha_6+16\alpha_7+8\alpha_8,$\\
$\mu_{\beta_8}=16\alpha_1+24\alpha_2+32\alpha_3+48\alpha_4+40\alpha_5+30\alpha_6+20\alpha_7+10\alpha_8$\\\\
\end{tabular}\\
\hline
$\Phi=F_4$&
\begin{tabular}{l}
$\mu_{\beta_1}=2\alpha_1+3\alpha_2+4\alpha_3+2\alpha_4,~\mu_{\beta_2}=2\alpha_1+4\alpha_2+6\alpha_3+4\alpha_4,$\\
$\mu_{\beta_3}=4\alpha_1+8\alpha_2+12\alpha_3+6\alpha_4,$\\
$\mu_{\beta_4}=6\alpha_1+12\alpha_2+16\alpha_3+8\alpha_4$\\\\
\end{tabular}\\
\hline
$\Phi=G_2$&$\mu_{\beta_1}=3\alpha_1+2\alpha_2,~\mu_{\beta_2}=4\alpha_1+2\alpha_2$\\&\\
\hline
\end{tabular}
\end{center}

\newpage Further, we also need to express the weights $\mu_{\beta_i}$'s as linear combinations of roots from $\Bu$ (it is possible since $\mu_{\beta_i}\in\Zp\Bu$ for all $i$). Such expressions are presented in Table \ref{table:weights_of_xi_i_cascade}.

\begin{center}
\tabl{Weights $\mu_{\beta}$, $\beta\in\Bu$, as linear combinations of roots from $\Bu$}\label{table:weights_of_xi_i_cascade}
\begin{tabular}{|l|l|}
\hline
$\Phi=A_{n-1}$&$\mu_{\beta_i}=\beta_1+\ldots+\beta_i$,~$1\leq i\leq[n/2]$\\
\hline
$\Phi=B_n$&$\mu_{\beta_i}=\begin{cases}2\beta_1+2\beta_3+\ldots+2\beta_{i-2}+\beta_{i-1}+\beta_i&\text{for even $i$},\\
\beta_1+\beta_3\ldots+\beta_i&\text{for odd $i<n$}\\
\end{cases}$\\
\hline
$\Phi=C_n$&$\mu_{\beta_i}=\beta_1+\ldots+\beta_i$,~$1\leq i\leq n$\\
\hline
$\Phi=D_n$&$\mu_{\beta_i}=\begin{cases}2\beta_1+2\beta_3+\ldots+2\beta_{i-2}+\beta_{i-1}+\beta_i&\text{for even $i<n$},\\
\beta_1+\beta_3+\ldots+\beta_{n-3}+\beta_n&\text{for even $i=n$},\\
\beta_1+\beta_3\ldots+\beta_i&\text{for odd $i$}
\end{cases}$\\
\hline
$\Phi=E_6$&$\mu_{\beta_1}=\beta_1,~\mu_{\beta_2}=\beta_1+\beta_2,~\mu_{\beta_3}=2\beta_1+\beta_2+\beta_3,~\mu_{\beta_4}=3\beta_1+\beta_2+\beta_3+\beta_4$\\
\hline
$\Phi=E_7$&
\begin{tabular}{l}
$\mu_{\beta_1}=\beta_1,~\mu_{\beta_2}=\beta_1+\beta_2,$\\
$\mu_{\beta_3}=2\beta_1+\beta_2+\beta_3,~\mu_{\beta_4}=\beta_1+\beta_2+\beta_4$\\
$\mu_{\beta_5}=2\beta_1+\beta_2+\beta_3+\beta_5,\mu_{\beta_6}=3\beta_1+\beta_2+\beta_3+\beta_6,$\\
$\mu_{\beta_7}=3\beta_1+2\beta_2+\beta_3+\beta_7$\\
\end{tabular}\\
\hline
$\Phi=E_8$&
\begin{tabular}{l}
$\mu_{\beta_1}=\beta_1,~\mu_{\beta_2}=\beta_1+\beta_2,$\\
$\mu_{\beta_3}=2\beta_1+\beta_2+\beta_3,~\mu_{\beta_4}=3\beta_1+2\beta_2+\beta_3+\beta_4,$\\
$\mu_{\beta_5}=3\beta_1+\beta_2+\beta_3+\beta_5,~\mu_{\beta_6}=3\beta_1+2\beta_2+\beta_3+\beta_4+\beta_6,$\\
$\mu_{\beta_7}=4\beta_1+3\beta_2+\beta_3+\beta_4+\beta_7,~\mu_{\beta_8}=5\beta_1+3\beta_2+2\beta_3+\beta_4+\beta_8$\\
\end{tabular}\\
\hline
$\Phi=F_4$&
\begin{tabular}{l}
$\mu_{\beta_1}=\beta_1,~\mu_{\beta_2}=\beta_1+\beta_2,$\\
$\mu_{\beta_3}=2\beta_1+\beta_2+\beta_3,~\mu_{\beta_4}=3\beta_1+\beta_2+\beta_3+\beta_4$\\
\end{tabular}\\
\hline
$\Phi=G_2$&$\mu_{\beta_1}=\beta_1,~\mu_{\beta_2}=\beta_1+\beta_2$\\
\hline
\end{tabular}
\end{center}

\nota{i) In fact, we\label{nota:values_of_xi_beta} will use Tables \ref{table:weights_of_xi_i}, \ref{table:weights_of_xi_i_basis}, \ref{table:weights_of_xi_i_cascade} only for exceptional root systems, but for the reader's convenience we describe the weights $\mu_{\beta_i}$ for all irreducible root systems.

ii) Note that the correspondence between $\beta_i$ and $\mu_{\beta_i}$ is uniquely determined by the fact that $$\mu_{\beta_i}\in\langle\beta_1,~\ldots,~\beta_i\rangle_{\Rp}\setminus\langle\beta_1,~\ldots,~\beta_{i-1}\rangle_{\Rp},$$ where $\langle\cdot\rangle_{\Rp}$, as usual, denotes the linear span over $\Rp$. Furthermore, each $\beta\in\Bu$ occurs in $\mu_{\beta}$ with coefficient $1$. Note also that if $\beta\in\Bu\cap\Delta$ then $\mu_{\beta}$ is the unique weight in Table~\ref{table:weights_of_xi_i_cascade} in which expression $\beta$ occurs. Our enumeration of the weights $\mu_{\beta_i}$ slightly differs from \cite{Panov2}, \cite{Ignatyev2} and \cite{IgnatyevPenkov1}.

iii) Recall that $\{e_{\alpha}^*,~\alpha\in\Phi^+\}$ is the basis of $\nt^*$ dual to the basis $\{e_{\alpha},~\alpha\in\Phi^+\}$ of $\nt$. Put $$R=\left\{t=\sum\nolimits_{\beta\in\Bu}t_{\beta}e_{\beta}^*,~t_{\beta}\in\Cp^{\times}\right\}$$ and denote by $X$ the union of all $N$-orbits in $\nt^*$ of elements of $R$. In fact, $X$ is a single $B$-orbit in $\nt^*$, and the $N$-orbits of two distinct point of $R$ are disjoint \cite[Theorem 2.5]{Kostant1}. Kostant \cite[Theorems 2.8, 3.11]{Kostant1} proved that $X$ is a Zariski dense subset of $\nt^*$, and for $t\in R$, up to scalar, for each $\beta_i\in\Bu$,
\begin{equation}
\xi_{\beta_i}(t)=\prod\nolimits_{\beta\in\Bu}t_{\beta}^{r_{\beta_i}(\beta)},~1\leq i\leq|\Bu|,\text{ where }\label{formula:xi_of_t}
r_{\beta_i}(\beta)=(\mu_{\beta_i},\beta)/(\beta,\beta).
\end{equation}
(Here we identify $S(\nt)$ with the algebra of polynomial functions on $\nt^*$.) Clearly, $r_{\beta_i}(\beta)$ is nothing but the coefficient at $\beta$ in the expression of $\mu_{\beta_i}$ in Table~\ref{table:weights_of_xi_i_cascade}.}

We fix the generators $\xi_{\beta}$, $\beta\in\Bu$, so that the formulas (\ref{formula:xi_of_t}) are satisfied (without any additional scalars). Recall the definition of the map $\sigma$ from (\ref{formula:symm_map}). For $\beta\in\Bu$, we denote $\Delta_{\beta}=\sigma(\xi_{\beta})\in Z(\nt)$. Explicit formulas for $\xi_{\beta}$ and $\Delta_{\beta}$ for classical root systems can be found in \cite[Subsection 2.1]{IgnatyevPenkov1}.

\defi{We call $\Delta_{\beta}$ (respectively, $\xi_{\beta}$), $\beta\in\Bu$, the \emph{canonical generators} of~the algebra $Z(\nt)$ (respectively, of the algebra $Y(\nt)$).}
\newpage

\sect{Centrally generated ideals}

Let\label{sect:cent_gen_ideals} $\gt$, $\nt$, $\Phi$, $\Delta$, etc., be as in Sec\-tion~\ref{sect:center}. A (two-sided) ideal $J\subset U(\nt)$ is called \emph{primitive} if $J$ is the annihilator of a simple $\nt$-module. An ideal $J$ is called \emph{centrally generated} if $J$ is generated (as an ideal) by its intersection $J\cap Z(\nt)$ with the center $Z(\nt)$ of $U(\nt)$.

In the 1960s Kirillov, Kostant and Souriau discovered that the orbits of the coadjoint action play a crucial role in the representation theory of $B$ and $N$ (see, e.g., \cite{Kirillov1}, \cite{Kirillov2}). Works of Dixmier, Duflo, Vergne, Mathieu, Conze and Rentschler led to the result that the orbit method provides a nice description of primitive ideals of the universal enveloping algebra of a nilpotent Lie algebra (in particular, of $\nt$). Below we briefly recall this description.

To any linear form $\lambda\in\nt^*$ one can assign a bilinear form $\beta_{\lambda}$ on $\nt$ by putting $\beta_{\lambda}(x,y)=\lambda([x,y])$ for $x,y\in\nt$. A subalgebra $\pt\subseteq\nt$ is a~\emph{polarization of $\nt$ at} $\lambda$ if it is a maximal $\beta_{\lambda}$-isotropic subspace. By \cite{Vergne}, such a subalgebra always exists. Let $\pt$ be a polarization of $\nt$ at $\lambda$, and let $W$ be the one-dimensional representation of $\pt$ defined by $x\mapsto\lambda(x)$. Then the annihilator $J(\lambda)=\Ann_{U(\nt)}{V}$ of the induced representation $V = U(\nt) \otimes_{U(\pt)} W$ is a primitive two-sided ideal of $U(\nt)$. It turns out that $J(\lambda)$ depends only on $\lambda$ and not on the choice of polarization. Further, $J(\lambda)=J(\mu)$ if and only if the coadjoint $N$-orbits of $\lambda$ and $\mu$ coincide. Finally, the \emph{Dixmier map} $$\Du\colon\nt^*\to\Prim U(\nt),~\lambda\mapsto J(\lambda),$$ induces a homeomorphism between $\nt^*/N$ and $\Prim U(\nt)$, where the latter set is endowed with the Jacobson topology. (See \cite{Dixmier2}, \cite{Dixmier3}, \cite{BorhoGabrielRentschler} for the details.)

In addition, it is well known that the following conditions on an ideal $J\subset U(\nt)$ are equivalent\break \cite[Proposition 4.7.4, Theorem 4.7.9]{Dixmier3}:
\begin{equation}
\begin{split}
&\text{i) $J$ is \label{formula:conditions_Prim_fd}primitive;}\\
&\text{ii) $J$ is maximal;}\\
&\text{iii) the center of $U(\nt)/J$ is trivial;}\\
&\text{iv) $U(\nt)/J$ is isomorphic to a Weyl algebra of finitely many variables.}
\end{split}
\end{equation}
Recall that the Weyl algebra $\Au_s$ of $2s$ variables is the unital associative algebra with generators $p_i$,~$q_i$ for $1\leq i\leq s$, and relations $[p_i,q_i]=1$, $[p_i,q_j]=0$ for $i\neq j$, $[p_i,p_j]=[q_i,q_j]=0$ for all $i,~j$. Furthermore, in (\ref{formula:conditions_Prim_fd}) we have $U(\nt)/J\cong\Au_s$ where $s$ equals one half of the dimension of the coadjoint $N$-orbit of $\lambda$, given that $J=J(\lambda)$.

\defi{To a map $\xi\colon\Bu\to\Cp$ we assign the linear form $f_{\xi}=\sum_{\beta\in\Bu}\xi(\beta)e_{\beta}^*\in\nt^*$. We call a form $f_{\xi}$ a \emph{Kostant form} if $\xi(\beta)\neq0$ for any $\beta\in\Bu\setminus\Delta$.}

Let $V$ be a simple (and hence at most countably dimensional) $\nt$-module and let $J=\Ann_{U(\nt)}{V}$ be the corresponding primitive ideal of $U(\nt)$. By a version of Schur's Lemma \cite{Dixmier4}, each central element of $U(\nt)$ acts on $V$ as a scalar operator. Given a tuple $c=(c_{\beta})_{\beta\in\Bu}$ of complex numbers, we denote by $J_c$ the ideal of $U(\nt)$ generated by all $\Delta_{\beta}-c_\beta$, $\beta\in\Bu$. If, for each $\beta\in\Bu$, $\Delta_{\beta}$ acts on $V$ by the scalar $c_{\beta}$, then, clearly, $J_c\subseteq J$. Further, since $Z(\nt)$ is a polynomial ring and the center of $U(\nt)/J$ is trivial, $J$~is centrally generated if and only if $J=J_c$.

The following result was proved in \cite[Theorem 3.1]{IgnatyevPenkov1} and \cite[Theorem 2.4]{Ignatyev2}.

\mtheo{Suppose $\Phi$ is an irreducible root system of classical type\textup, i.e.\textup, $\Phi=A_{n-1}$\textup, $B_n$\textup, $C_n$ or $D_n$. The following conditions on a primitive ideal $J\subset U(\nt)$ are equivalent\textup{:}
\begin{equation*}
\begin{split}
&\text{\textup{i)} $J$ is centrally generated \textup{(}or\textup{,} equivalently\textup{,} $J=J_c$\textup{)}};\\
&\text{\textup{ii)} the scalars $c_{\beta}$\textup, $\beta\in\Bu\setminus\Delta$\textup, are nonzero};\\
&\text{\textup{iii)} $J=J(f_{\xi})$ for a Kostant form $f_{\xi}\in\nt^*$}.\\
\end{split}
\end{equation*}
If these\label{theo:ideal_classical} conditions are satisfied\textup{,} then the map $\xi$ can be reconstructed from $J$.}\newpage

The main result of the paper is to prove that this is also true for exceptional root systems, see Theorem~\ref{theo:ideal_exceptional} in Section~\ref{sect:proof_main}. One of the key ingredients in the proof of Theorem~\ref{theo:ideal_classical} was to check that if condition (ii) is satisfied then $J_c$ is primitive. To do this for $A_{n-1}$ and $C_n$, in \cite{IgnatyevPenkov1} an explicit set of generators of the quotient algebra $U(\nt)/J_c$ was constructed. It turns out that these generators satisfy (up to scalars) the defining relations of the Weyl algebra $\Au_s$ for $s=|\Phi^+\setminus\Bu|/2$. Since $\Au_s$ is simple and, as one can check, $J_c\neq U(\nt)$, we conclude that $U(\nt)/J_c\cong\Au_s$, and, consequently, $J_c$ is primitive. On the other hand, for $B_n$ and $D_n$, in \cite{Ignatyev2} an explicit set of generators for $U(\nt)/J_c$ was constructed a posteriori (see \cite[Theorem 2.9]{Ignatyev2}), while primitivity of $J_c$ was established by another argument. In this section we modify the idea from \cite[Proposition 2.5]{Ignatyev2} to check that $J_c$ is primitive if $c_{\beta}\neq0$ for $\beta\in\Bu\setminus\Delta$ for exceptional types.

To do this, we need some additional notation. From now to the end of this section we assume that $\Phi$ is an irreducible root system of exceptional type, i.e., $\Phi=E_6$, $E_7$, $E_8$, $F_4$ or $G_2$. Recall that $\beta_1$ is the maximal root with respect to the natural order on $\Phi$. It is obvious that $(\alpha,\beta_1)\geq0$ for all $\alpha\in\Phi^+$. We put $\wt\Phi=\{\alpha\in\Phi\mid(\alpha,\beta_1)=0\}$ and $\wt\Phi^+=\wt\Phi\cap\Phi^+$, $\wt\Delta=\wt\Phi\cap\Delta$. Then $\wt\Phi$ is of respective type $D_5$, $E_6$, $E_7$, $C_3$ or $A_1$. Denote
\begin{equation*}
\begin{split}
\wt\nt&=\langle e_{\alpha},~\alpha\in\Phi^+,~(\alpha,\beta_1)=0\rangle_{\Cp}=\langle e_{\alpha},~\alpha\in\wt\Phi^+\rangle_{\Cp},\\
\kt&=\langle e_{\alpha},~\alpha\in\Phi^+,~(\alpha,\beta_1)>0\rangle_{\Cp},\\
\end{split}
\end{equation*}
where $\langle\cdot\rangle_{\Cp}$ is the $\Cp$-linear span. Then $\wt\nt$ is a Lie subalgebra of $\nt$ isomorphic to the nilradical of the Borel subalgebra $\wt\bt=\wt\gt\cap\bt$ of the simple Lie algebra $\wt\gt$ with the root system $\wt\Phi$, where $\wt\gt$ is the subalgebra of $\gt$ generated by the root vectors $e_{\alpha}$, $\alpha\in\wt\Phi$.

On the other hand, $\kt$ is an ideal of $\nt$ isomorphic to the Heisenberg Lie algebra $\hei_s$, where ${s=(|\Phi^+\setminus\wt\Phi^+|-1)/2}$, with the center $\Cp e_{\beta_1}$. (This follows from the fact that if $\alpha\in\Phi^+$ and $(\alpha,\beta_1)>0$ then $\beta_1-\alpha$ is again a positive root, see \cite[Corollary 2.3]{Joseph2} for the details.) Recall that $\hei_s$ is the $(2s+1)$-dimensional Lie algebra with basis $\{z,~x_i,~y_i,~1\leq i\leq s\}$ and relations $[x_i,~y_i]=z$ for all $i$, $[x_i,z]=[y_j,z]=[x_i,y_j]=0$ for all $i\neq j$.

Given $c_1\in\Cp^{\times}$, denote by $J_1$ the ideal of $U(\kt)$ generated by $e_{\beta_1}-c_1$, so, clearly, $U(\kt)/J_1\cong\Au_s$. Since $\kt$ is an ideal of the Lie algebra $\nt$, given $x\in\wt\nt$, one can consider $\ad{x}$ as a derivation of $\kt$. Since $e_{\beta_1}-c_1$ is a central element of $U(\nt)$, one has $\ad{x}(J_1)\subseteq J_1$, so $\ad{x}$ can be considered as a derivation of~$\Au_s$. It is well known (see, e.g., \cite[10.1.1--10.1.4]{Dixmier3}) that there exists unique $\theta(x)\in\Au_s$ such that $\ad{x}(y)=[\theta(x),y]$ for all $y\in\Au_s$, and $\theta\colon\wt\nt\to\Au_s$ is a morphism of Lie algebras. Furthermore, there exists the unique epimorphism of associative algebras $r\colon U(\nt)\sur U(\wt\nt)\otimes\Au_s$ such that $r(y)=1\otimes \overline y$ for $y\in\kt$ and $r(x)=x\otimes1+1\otimes\theta(x)$  for $x\in\wt\nt$. (Here $\overline a$ is the image of an element $a\in U(\kt)$ under the canonical projection $U(\kt)\sur U(\kt)/J_1\cong\Au_s$.) It turns out that the kernel of the epimorphism $r$ coincides with the ideal $J_0$ of $U(\nt)$ generated by $e_{\beta_1}-c_1$ \cite[Lemma 10.1.5]{Dixmier3}.

\propp{Let\label{prop:J_c_primitive_exceptional} $J_c$ be the ideal of $U(\nt)$ generated by $\Delta_{\beta}-c_{\beta}$ for all $\beta\in\Bu$ with $c_{\beta}\neq0$ for $\beta\in\Bu\setminus\Delta$. Then $J_c$ is primitive.}{Put $c_1=c_{\beta_1}$. Since $r$ is surjective, $r(J_c)$ is an ideal of $U(\wt\nt)\otimes\Au_s$ generated by $r(\Delta_{\beta})-c_{\beta}$, $\beta\in\wt\Bu=\Bu\setminus\{\beta_1\}$. Note that $\wt\Bu$ is the Kostant cascade of $\wt\Phi$. Denote by $\wt\Delta_{\beta}$, $\beta\in\wt\Phi$, the set of canonical generators of $Z(\wt\nt)$. We will show that, up to nonzero scalar, $r(\Delta_{\beta})$ coincides with $\wt\Delta_{\beta}\otimes1$ for all $\beta\in\wt\Bu$.

To check this fact, we will use Tables \ref{table:weights_of_xi_i}, \ref{table:weights_of_xi_i_basis}, \ref{table:weights_of_xi_i_cascade}. Pick a root $\beta\in\wt\Bu$. Since $r$ is surjective, $r(\Delta_{\beta})$ is central in $U(\wt\nt)\otimes\Au_s$. The center of this algebra has the form $Z(\wt\nt)\otimes\Cp$, so in fact $r(\Delta_{\beta})\in Z(\wt\nt)\otimes\Cp$. Denote $\wt\htt=\wt\gt\cap\htt$, so $\wt\htt$ is a Cartan subalgebra of $\wt\gt$ and $\wt\bt=\wt\htt\oplus\wt\nt$ as vector spaces. By \cite[Theorem~6]{Kostant2} (see also \cite[Lemma 4.4]{Joseph1}), $Z(\nt)$ (respectively, $Z(\wt\nt)$) is a direct sum of one-dimensional weight spaces of~$\htt$ (respectively, of $\wt\htt$) with respect to the adjoint action of the corresponding Cartan subalgebras. Since $[h,e_{\beta_1}]=0$ for all $h\in\wt\htt$, the algebra $\wt\htt$ naturally acts on $\Au_s$, and so on $U(\wt\nt)\otimes\Au_s$. We define the result of this action by $h.x$, $h\in\wt\htt$, $x\in U(\wt\nt)\otimes\Au_s$. Hence it is enough to check that, given $\beta\in\wt\Bu$, $r(\Delta_{\beta})$ is a nonzero $\wt\htt$-weight element of weight $\wt\mu_{\beta}=\mu_{\beta}-r_{\beta}(\beta_1)\beta_1$, because $r_{\beta}(\beta_1)$ is just the coefficient of $\beta_1$ in the expression of $\mu_{\beta}$ as a linear combination of roots from $\Bu$.

To prove that $r(\Delta_{\beta})$ is an $\wt\htt$-weight element of weight $\wt\mu_{\beta}$, denote the result of the natural (adjoint) action of $\wt\htt$ on $U(\nt)$ by $h\cdot x$, $h\in\wt\htt$, $x\in U(\nt)$. As above, since $\wt\htt\cdot e_{\beta_1}=0$, the algebra $\wt\htt$ naturally acts on $U(\nt)/J_0$ by the formula $h\cdot r(x)=r(h\cdot x)$. We claim that this action coincides with the action of $\wt\htt$ on $U(\wt\nt)\otimes\Au_s$ defined above, i.e., $h.x=h\cdot x$ for all $h\in\wt\htt$, $x\in U(\wt\nt)\otimes\Au_s$.

Indeed, if $y\in\kt$, then $$h\cdot r(y)=r([h,y])=1\otimes\overline{[h,y]}=h.(1\otimes\overline{y})=h.r(y),$$ as required. On the other hand, if $e_{\alpha}\in\wt\nt$ for some root $\alpha\in\Phi^+$, then, by \cite[Subsection~4.8]{Joseph1}, $\theta(e_{\alpha})$ is a linear combination of elements of the form $\overline{e}_{\alpha+\gamma}\overline{e}_{\beta_1-\gamma}$, $\gamma\in\Phi^+\setminus\wt\Phi^+$ (i.e., $(\gamma,\beta_1)>0$). We conclude that $$h(1\otimes\theta(e_{\alpha}))=(\beta_1+\alpha)(h)1\otimes\theta(\alpha)=\alpha(h)1\otimes\theta(e_{\alpha}),$$ because $\beta_1(\wt\htt)=0$. Thus, we obtain
$$h.r(e_{\alpha})=[h,e_{\alpha}]\otimes1+1\otimes\alpha(h)\theta(e_{\alpha})=\alpha(h)r(e_{\alpha})=h\cdot r(e_{\alpha}).$$ It remains to note that $\Delta_{\beta}$ is an $\htt$-weight element of $U(\nt)$ of weight $r_{\beta}(\beta_1)\beta_1+\wt\mu_{\beta}$, but $\beta_1(\wt\htt)=0$.

To show that $r(\Delta_{\beta})\neq0$, recall that the kernel of $r$ is $J_0$. If $\Delta_{\beta}\in J_0$ (i.e., if $\Delta_{\beta}=(e_{\beta_1}-c_{\beta_1})a$ for some $a\in U(\nt)$), then, clearly, $a\in Z(\nt)$. But this contradicts the fact that $\Delta_{\beta}$ and $\Delta_{\beta_1}$ are algebraically independent, because, as one can deduce from \cite{Panov2}, $\Delta_{\beta_1}=e_{\beta_1}$ for all irreducible root systems.

So, given $\beta\in\wt\Bu$, there exists unique $a_{\beta}\in\Cp^{\times}$ such that $r(\Delta_{\beta})=a_{\beta}\wt\Delta_{\beta}\otimes1$. Consequently, $r(J_c)$ is generated by $\wt\Delta_{\beta}\otimes1-\wt c_{\beta}$, $\beta\in\wt\Bu$, where $\wt c_{\beta}=a_{\beta}^{-1}c_{\beta}$. In particular, $\wt c_{\beta}\neq0$ if $\beta\in\wt\Bu$ is not a simple root of $\wt\Phi^+$. Now we will use the induction on $\rk\Phi$ to prove that $J_c$ is primitive. The base (i.e., the case of classical $\Phi$ of low rank) immediately follows from \cite[Theorem 3.1]{IgnatyevPenkov1} and \cite[Theorem 2.4]{Ignatyev2}. Denote by $\wt J_c$ the ideal of $U(\wt\nt)$ generated by $\wt\Delta_{\beta}-\wt c_{\beta}$, $\beta\in\wt\Bu$. By the inductive assumption, $\wt J_c$ is a primitive ideal of $U(\wt\nt)$, so $U(\wt\nt)/\wt J_c\cong\Au_t$ for certain $t$. We conclude that
\begin{equation*}
U(\nt)/J_c\cong(U(\nt)/J_0)/r(J_c)\cong(U(\wt\nt)\otimes\Au_s)/r(J_c)=(U(\wt\nt)/\wt J_c)\otimes\Au_s\cong\Au_t\otimes\Au_s\cong\Au_{t+s}.
\end{equation*}
Thus, $J_c$ is primitive. The proof is complete.}

\sect{Distinct coadjoint orbits}\label{sect:distinct_orbits}

Recall that, given a primitive ideal $J$ in $U(\nt)$, there exist unique scalars $c_{\beta}\in\Cp$ such that $\Delta_{\beta}-c_{\beta}\in J$ for all $\beta\in\Bu$. To prove our main result, Theorem~\ref{theo:ideal_exceptional}, we need to check that if $J$ is centrally generated then $c_{\beta}\neq0$ for $\beta\in\Bu\setminus\Delta$. To do this, we will prove that certain coadjoint $N$-orbits on $\nt^*$ are distinct.

Namely, let $D$ be a subset of $\Phi^+$. To each map $\xi\colon D\to\Cp^{\times}$ one can assign the linear form $$f_{D,\xi}=\sum_{\beta\in D}\xi(\beta)e_{\beta}^*\in\nt^*.$$ Denote by $\Omega_{D,\xi}$ the coadjoint $N$-orbit of $f_{D,\xi}$. We say that $f_{D,\xi}$ and $\Omega_{D,\xi}$ are associated with the subset~$D$. For example, $f_{D,\xi}$ is a Kostant form if and only if $\Bu\setminus\Delta\subset D\subset\Bu$.

It was proved in \cite[Corollary 1.4]{Panov1} that if $\Phi=A_{n-1}$, $D$ is an orthogonal subset (i.e., $(\alpha,\beta)=0$ for all $\alpha,~\beta\in D$, $\alpha\neq\beta$) and $\xi_1$, $\xi_2$ are two distinct maps from $D$ to~$\Cp^{\times}$ then $\Omega_{D,\xi_1}\neq\Omega_{D,\xi_2}$. It is not hard to deduce from this result that the same is true for all classical root systems, see the proofs\break of \cite[Theorem 3.1]{IgnatyevPenkov1} and \cite[Theorem 2.4]{Ignatyev2}. But for exceptional types this is not an immediate consequence of the result for $A_{n-1}$. In this section, we prove that if $\xi_1\neq\xi_2$ then $\Omega_{D,\xi_1}$ and $\Omega_{D,\xi_2}$ are distinct for some particular orthogonal subsets $D$ and some particular maps $\xi_1$, $\xi_2$, which will be used in the next section in the proof of our main result.\newpage

To do this, we need to introduce the notion of singular roots.

\defi{Let $\beta$, $\alpha$ be positive roots. We say that $\alpha$ is $\beta$-\emph{singular} (or \emph{singular} for $\beta$) if there exists $\gamma\in\Phi^+$ such that $\beta=\alpha+\gamma$. The set of all $\beta$-singular roots is denoted by $S(\beta)$.}

Note that if $\Phi$ is irreducible and simple-laced (i.e., if all roots in $\Phi$ have the same length) then, given $\beta>\alpha$, $\alpha$ is $\beta$-singular if and only if $(\alpha,\beta)>0$. It turns out that if $D$ is an orthogonal subset of $\Phi^+$, $\xi$ is a map from $D$ to $\Cp^{\times}$, and $\beta$, $\beta'\in D$ are such that $\beta'\in S(\beta)$ then $\Omega_{D,\xi}=\Omega_{D',\xi'}$, where $D'=D\setminus\beta'$ and $\xi'$ is the restriction of $\xi$ to $D'$ \cite[Lemma 1.3]{Ignatyev3}.

\propp{Let $\Phi$ be\label{prop:distinct_orbits} an irreducible root system\textup, and  let $D$ be a subset of $\Phi^+$ such that if $\beta_1$, $\beta_2\in D$ then $\beta_1\notin S(\beta_2)$. Let $\beta_0$ be a root in $D$\textup, and let $\xi_1$\textup, $\xi_2$ be maps from $D$ to $\Cp^{\times}$ for which $\xi_1(\beta_0)\neq\xi_2(\beta_0)$. Fix a total order $\leq_t$ on $\{e_{\alpha},~\alpha\in\Phi^+\}$ such that if $\alpha>\beta$ then $e_{\alpha}<_t e_{\beta}$. Assume that there exists a simple root $\alpha_0\in\Delta$ satisfying $(\alpha_0,\beta_0)\neq0$ and $(\alpha_0,\beta)=0$ for all $\beta\in D$ such that $e_{\beta}<_t e_{\beta_0}$. Then $\Omega_{D,\xi_1}\neq\Omega_{D,\xi_2}$.}{As usual, given a vector space $V$, we denote by $\glt(V)$ the Lie algebra of all linear operators on $V$. Recall that $\Phi$ is the root system of the Lie algebra $\gt$. Denote by $\aad\colon\gt\to\glt(\gt)$ the adjoint representation of the Lie algebra $\gt$, i.e., $\aad(x)=\ad{x}$. It is well known that the adjoint representation is faithful, so $\aad(\gt)$ and $\gt$ are isomorphic as Lie algebras. Let $\rk\Phi=n$. To each simple root $\alpha_i$, $1\leq i\leq n$, one can assign unique $h_{\alpha_i}\in\htt$ such that $\alpha_i(h_{\alpha_j})=2(\alpha_i,\alpha_j)/(\alpha_j,\alpha_j)$ for all $i,~j$. Recall that $\{e_{\alpha},~\alpha\in\Phi^+\}$ is a basis of $\nt$. It can be uniquely extended to the Chevalley basis $\{e_{\alpha},~\alpha\in\Phi^+\}\cup\{h_{\alpha_i},~1\leq i\leq n\}\cup\{e_{-\alpha},~\alpha\in\Phi^+\}$ of $\gt$. We extend the total order $\leq_t$ on the set $\{e_{\alpha}.~\alpha\in\Phi^++\}$ to a total order $\leq_t$ on this basis such that $e_{\alpha}<_th_{\alpha_i}<_te_{-\beta}$ for all $\alpha,\beta\in\Phi^+$, $1\leq i\leq n$, and $e_{\alpha}<_te_{\beta}$ if $\alpha,\beta\in\Phi$ and $\alpha>\beta$. This identifies $\glt(\gt)$ with the Lie algebra $\glt_{\dim\gt}(\Cp)$, and $\aad(\nt)$ with a subalgebra of the Lie algebra $\ut$ of all upper-triangular matrices from $\glt_{\dim\gt}(\Cp)$ with zeroes on the diagonal.

Let $\GL(V)$ be the group of all invertible linear operators on a vector space $V$. Since we fixed a basis in $\gt$, the group $\GL(\gt)$ is identified with the group $\GL_{\dim\gt}(\Cp)$, and $\exp\aad(\nt)\cong N$ is identified with a subgroup of the group $U$ of all upper-triangular matrices from $\GL_{\dim\gt}(\Cp)$ with $1$'s on the diagonal. Furthermore, using the Killing form on $\gt$ and the trace form on $\glt(\gt)$, one can identify $\nt^*$ with the space $\nt_-=\langle e_{-\alpha},~\alpha\in\Phi^+\rangle_{\Cp}$ and $\ut^*$ with the space $\ut_-=\ut^T$, where the superscript $T$ denote the transposed matrix. Under these identifications, it is enough to check that the coadjoint $U$-orbits of the linear forms $\wt f_{D,\xi_1}$ and $\wt f_{D,\xi_2}$ are distinct. Here, given $\xi\colon D\to\Cp^{\times}$, we denote by $\wt f_{D,\xi}$ the matrix $$\wt f_{D,\xi}=\left(\aad\left(\sum\nolimits_{\beta\in D}\xi(\beta)e_{\beta}\right)\right)^T\in\ut_-\cong\ut^*.$$

To do this, we will study the matrix $f=\wt f_{D,\xi}$ in more detail. The rows and the columns of matrices from $\glt(\gt)$ are now indexed by the elements of the Chevalley basis fixed above. Given a matrix $x$ from $\glt(\gt)$ and basis elements $a,~b$, we will denote by $x_{a,b}$ the entry of $x$ lying in the $a$th row and the $b$th column. Since $$\ad{e_{\beta_0}}(h_{\alpha_0})=[e_{\beta_0,h_{\alpha_0}}]=-\dfrac{2(\alpha_0,\beta_0)}{(\alpha_0,\alpha_0)}e_{\beta_0},$$ we obtain $f_{h_{\alpha_0},e_{\beta_0}}=-\xi(\beta_0)\dfrac{2(\alpha_0,\beta_0)}{(\alpha_0,\alpha_0)}\neq0$. One may assume without loss of generality that $h_{\alpha_0}>_th_{\alpha_i}$ for all $\alpha_i\neq\alpha_0$. We claim that
\begin{equation}\label{formula:condition_Andre}
f_{h_{\alpha_0},e_{\alpha}}=f_{e_{-\gamma},e_{\beta_0}}=0\text{ for all }e_{\alpha}<_te_{\beta_0}\text{ and all }e_{-\gamma},~\alpha,\gamma\in\Phi^+.
\end{equation}
Indeed, if $\alpha\notin D$ then, evidently, $f_{h_{\alpha_0},e_{\alpha}}=0$. If $\alpha=\beta\in D$ and $e_{\beta}<_te_{\beta_0}$ then $(\alpha_0,~\beta)=0$, hence $$f_{h_{\alpha_0},e_{\beta}}=-\xi(\beta)\dfrac{2(\alpha_0,\beta)}{(\alpha_0,\alpha_0)}=0,$$ because $(\alpha_0,\beta)=0$. On the other hand, if $f_{e_{-\gamma},e_{\beta_0}}\neq0$ for some $\gamma\in\Phi^+$ then $\beta_0=\beta-\gamma$. This contradicts the condition $\beta_0\notin S(\beta)$.\newpage

Thus, $(\wt f_{D,\xi_1})_{h_{\alpha_0},e_{\alpha}}$ and $(\wt f_{D,\xi_2})_{h_{\alpha_0},e_{\alpha}}$ are different nonzero scalars, and (\ref{formula:condition_Andre}) is satisfied both for $f=\wt f_{D,\xi_1}$ and for $f=\wt f_{D,\xi_2}$. Now it follows immediately from the proof of \cite[Proposition 3]{Andre1} that the coadjoint $U$-orbits of these matrices are distinct, and, consequently, $\Omega_{D,\xi_1}\neq\Omega_{D,\xi_2}$, as required.}

Now, for exceptional $\Phi$, we will prove an existence of certain subsets $D\subset\Phi^+$. To each such subset $D$ from this list we assign its subset $D'\subset D$. Using Proposition~\ref{prop:distinct_orbits}, we will show that if $\xi_1$ and $\xi_2$ are two maps from $D$ to~$\Cp^{\times}$ such that $\xi_1(\beta_0)\neq\xi_2(\beta_0)$ for some root $\beta_0\in D'$ then $\Omega_{D,\xi_1}\neq\Omega_{D,\xi_2}$, see Proposition~\ref{prop:distinc_orbits_precise} below. We will consider all exceptional root systems subsequently. For brevity, we use the following notation. If $\beta=\sum_{i=1}^nm_i\alpha_i\in D$ then we write $m_1\ldots m_n$ instead of $\beta$ (our enumeration of simple roots is as in (\ref{formula:simple_roots})). Note that in all cases, except 11, 12 and 14 for $F_4$ (see Table~\ref{table:D_D'_F_4} in the proof of Proposition~\ref{prop:subsets D_conditions} below), $D$ is an orthogonal subset of $\Phi$, while in cases 11, 12, 14 for $F_4$ all inner products of distinct roots from~$D$ are non-positive. Note also that all these subsets are linearly independent.

\propp{Let $\Phi$ be\label{prop:subsets D_conditions} an irreducible root system of exceptional type. Denote by $c$ an $m$-tuple $c=(c_{\beta_1},~\ldots,~c_{\beta_m})$\textup, $m=|\Bu|$\textup, and assume that $c_{\beta_i}=0$ for some $\beta_i\in\Bu\setminus\Delta$. There exist subsets $D,~D'\subset\Phi^+$ such that $D$ is linearly independent\textup, $D'\subset D$\textup, and the following conditions are satisfied\textup:
\begin{equation*}
\begin{split}
&\text{\textup{i)} $D$ is linear independent \textup(in fact\textup, orthogonal\textup, except cases $11$\textup, $12$ and $14$ for $F_4$ in Table \textup{\ref{table:D_D'_F_4}}\textup)};\\
&\text{\textup{ii)} if $c_{\beta_i}=0$ for some $\beta_i\in\Bu$ then $\mu_{\beta_i}$ does not belong to the $\Zp_{\geq0}$-linear span of $D$};\\
&\text{\textup{iii)} if $c_{\beta_i}\neq0$ for some $\beta_i\in\Bu$ then $\mu_{\beta_i}$ belongs to the $\Zp_{\geq_0}$-linear span of $D$};\\
&\text{\textup{iv)} there exist $\alpha_0\in\Delta$ and $\beta_0\in D'$ such that $(\alpha_0,\beta_0)\neq0$ and $(\alpha_0,\beta)=0$ for all $\beta\in D\setminus D'$}.\\
\end{split}
\end{equation*}
Furthermore\textup, let $\mu_{\beta_i}=\sum\nolimits_{\gamma\in D}a_{\gamma,\beta_i}\gamma,~1\leq i\leq m,~c_{\beta_i}\neq0$. Then\textup, for $\beta_0\in D$\textup, the solution space of the system of linear equations $\sum\nolimits_{\gamma\in D}a_{\gamma,\beta_i}y_{\gamma}=b_i,~1\leq i\leq m,~c_{\beta_i}\neq0,$ is not orthogonal to the axis $y_{\beta_0}$ if and only if $\beta_0\in D'$.}{The proof is case-by-case. Below we present the list of the subsets $D$, $D'$ and the cor\-res\-pon\-dence between these subsets and $m$-tuples $c$ for $E_6$, $F_4$ and $G_2$.}

\begin{center}
\tabl{List of the subsets $D$ and $D'$ for $E_6$}\label{table:D_D'_E_6}
\end{center}
\vspace{-0.7cm}\begin{longtable}{||l|l|l||l|l|l||l|l|l||}
\hline
&$D$&$D'$&&$D$&$D'$&&$D$&$D'$\\
\hline\hline
1
&\begin{tabular}[t]{l}$010000$\end{tabular}&
\begin{tabular}[t]{l}$010000$\end{tabular}&
2
&\begin{tabular}[t]{l}$011211$,\\
$111221$,\\
$112210$
\end{tabular}&
\begin{tabular}[t]{l}\\
$111221$,\\
$112210$
\end{tabular}&
3
&\begin{tabular}[t]{l}$011111$,\\
$111110$,\\
$112321$
\end{tabular}&
\begin{tabular}[t]{l}\\
$111110$,\\
$112321$
\end{tabular}\\
\hline
4
&\begin{tabular}[t]{l}$011211$,\\
$111210$,\\
$112221$
\end{tabular}&
\begin{tabular}[t]{l}$011211$,\\
$111210$
\end{tabular}&
5
&\begin{tabular}[t]{l}$111111$,\\
$112321$
\end{tabular}&
\begin{tabular}[t]{l}$111111$,\\
$112321$
\end{tabular}&
6
&\begin{tabular}[t]{l}$011210$,\\
$111221$,\\
$112211$
\end{tabular}&
\begin{tabular}[t]{l}\\
$111221$,\\
$112211$
\end{tabular}\\
\hline
7
&\begin{tabular}[t]{l}$011110$,\\
$111111$,\\
$112321$
\end{tabular}&
\begin{tabular}[t]{l}$011110$,\\
$111111$,\\
$112321$
\end{tabular}&
8
&\begin{tabular}[t]{l}$010000$,\\
$011210$,\\
$111211$,\\
$112221$
\end{tabular}&
\begin{tabular}[t]{l}$010000$
\end{tabular}&
9
&\begin{tabular}[t]{l}$001000$,\\
$122321$
\end{tabular}&
\begin{tabular}[t]{l}$001000$
\end{tabular}\\
\hline
10
&\begin{tabular}[t]{l}$001100$,\\
$000111$,\\
$122321$,\\
$101110$
\end{tabular}&
\begin{tabular}[t]{l}$001100$,\\
$000111$
\end{tabular}&
11
&\begin{tabular}[t]{l}$001111$,\\
$122321$,\\
$101110$
\end{tabular}&
\begin{tabular}[t]{l}$001111$,\\
\\
$101110$
\end{tabular}&
12
&\begin{tabular}[t]{l}$001111$,\\
$000100$,\\
$122321$,\\
$101110$
\end{tabular}&
\begin{tabular}[t]{l}$001111$,\\
\\
\\
$101110$
\end{tabular}\\
\hline
13
&\begin{tabular}[t]{l}$001000$,\\
$101111$,\\
$122321$
\end{tabular}&
\begin{tabular}[t]{l}$001000$
\end{tabular}&
14
&\begin{tabular}[t]{l}$001100$,\\
$000110$,\\
$122321$,\\
$101111$
\end{tabular}&
\begin{tabular}[t]{l}$001100$,\\
$000110$
\end{tabular}&&&\\
\hline
\end{longtable}

\newpage $\hphantom{1}$
\vspace{-1.6cm}
\begin{center}
\tabl{List of the subsets $D$ and $D'$ for $F_4$}\label{table:D_D'_F_4}
\end{center}
\vspace{-0.7cm}\begin{longtable}{||l|l|l||l|l|l||l|l|l||}
	\hline
	&$D$&$D'$&&$D$&$D'$&&$D$&$D'$\\
	\hline\hline
1&\begin{tabular}[t]{l}
	$1000$\end{tabular}&\begin{tabular}[t]{l}
	$1000$\end{tabular}&2&\begin{tabular}[t]{l}
	$1231$,\\$1222$\end{tabular}&\begin{tabular}[t]{l}
	$1231$,\\$1222$\end{tabular}&3&\begin{tabular}[t]{l}
	$1121$,\\$1342$\end{tabular}&\begin{tabular}[t]{l}
	$1121$,\\$1342$\end{tabular}\\\hline4&\begin{tabular}[t]{l}
	$1000$,\\$1221$,\\$1242$\end{tabular}&\begin{tabular}[t]{l}
	$1000$\\ \\ \end{tabular}&5&\begin{tabular}[t]{l}
	$1000$,\\$1232$\end{tabular}&\begin{tabular}[t]{l}
	$1000$\\ \end{tabular}&6&\begin{tabular}[t]{l}
	$1000$,\\$1220$,\\$1232$\end{tabular}&\begin{tabular}[t]{l}
	$1000$\\ \\ \end{tabular}\\\hline7&\begin{tabular}[t]{l}
	$1120$,\\$1122$,\\$1342$\end{tabular}&\begin{tabular}[t]{l}
	$1120$,\\$1122$,\\$1342$\end{tabular}&8&\begin{tabular}[t]{l}
	$1000$,\\$1220$,\\$1222$,\\$1242$\end{tabular}&\begin{tabular}[t]{l}
	$1000$\\ \\ \\ \end{tabular}&9&\begin{tabular}[t]{l}
	$0100$,\\$2342$\end{tabular}&\begin{tabular}[t]{l}
	$0100$\\ \end{tabular}\\\hline10&\begin{tabular}[t]{l}
	$0111$,\\$0120$,\\$2342$\end{tabular}&\begin{tabular}[t]{l}
	$0111$,\\$0120$\\ \end{tabular}&11&\begin{tabular}[t]{l}
	$1121$,\\$1342$,\\ $2342$\end{tabular}&\begin{tabular}[t]{l}
	$1121$,\\ $1342$ \\ \end{tabular}&12&\begin{tabular}[t]{l}
	$1100$,\\$1121$,\\$1342$,\\ $2342$\end{tabular}&\begin{tabular}[t]{l}
	$1100$,\\$1121$,\\$1342$\\ \end{tabular}\\\hline13&\begin{tabular}[t]{l}
	$0100$,\\$0122$,\\$2342$\end{tabular}&\begin{tabular}[t]{l}
	$0100$\\ \\ \end{tabular}&14&\begin{tabular}[t]{l}
	$1110$,\\$0122$,\\$1342$,\\$2342$\end{tabular}&\begin{tabular}[t]{l}
	$1110$,\\\\$1342$ \\ \end{tabular}&&\begin{tabular}[t]{l}
	\\ \\ \\\end{tabular}&\begin{tabular}[t]{l}
	\\ \\ \\ \end{tabular}\\\hline

\end{longtable}

\vspace{-0.7cm}\begin{center}
\tabl{Correspondence between $c$ and $D$ for $E_6$ and $F_4$}\label{table:c_D_E_6_F_4}\vspace{-0.2cm}

\vspace{0.0cm}\begin{longtable}{||l|l||l|l||l|l||}
\hline
$D$&Type of $c$&$D$&Type of $c$&$D$&Type of $c$\\
\hline\hline
1&$(0,0,0,0)$&2&$(0,0,0,c_{\beta_4})$&3&$(0,0,c_{\beta_3},0)$\\
\hline
4&$(0,0,c_{\beta_3},c_{\beta_4})$&5&$(0,c_{\beta_2},0,0)$&6&$(0,c_{\beta_2},0,c_{\beta_4})$\\
\hline
7&$(0,c_{\beta_2},c_{\beta_3},0)$&8&$(0,c_{\beta_2},c_{\beta_3},c_{\beta_4})$&9&$(c_{\beta_1},0,0,0)$\\
\hline
10&$(c_{\beta_1},0,0,c_{\beta_4})$&11&$(c_{\beta_1},0,c_{\beta_3},0)$&12&$(c_{\beta_1},0,c_{\beta_3},c_{\beta_4})$\\
\hline
13&$(c_{\beta_1},c_{\beta_2},0,0)$&14&$(c_{\beta_1},c_{\beta_2},0,c_{\beta_4})$&&\\
\hline
\end{longtable}
\end{center}

\begin{center}
\tabl{List of the subsets $D$ and $D'$ and correspondence between $c$ and $D$ for $G_2$}\label{table:c_D_D'_G_2}
\end{center}
\vspace{-0.7cm}\begin{longtable}{||l|l|l|l||l|l|l|l||}
	\hline
	&$D$&$D'$&Type of $c$&&$D$&$D'$&Type of $c$\\
	\hline\hline		1&\begin{tabular}[t]{l}
		$31$\end{tabular}&\begin{tabular}[t]{l}
		$31$\end{tabular}&$(0,~0)$&2&\begin{tabular}[t]{l}
		$31,$\\$11$\end{tabular}&\begin{tabular}[t]{l}
		$31,$\\$11$\end{tabular}&$(0,~c_{\beta_2})$\\
	\hline
\end{longtable}

The root systems $E_7$, $E_8$ can be considered similarly. We list the subsets $D$, $D'$ and prove the correspondence between $c$ and $D$ for $E_7$ and $E_8$ in Appendix A. (In fact, we just checked all orthogonal subsets of $\Phi^+$ to find subsets satisfying the required conditions using Python Programming Lan\-guage~\cite{Python}. The listing of the code is presented in Appendix B.)

\propp{Let $\Phi$ be\label{prop:distinc_orbits_precise} an irreducible root system of exceptional type\textup, and let $D\subset\Phi^+$ be a subset from Proposition~$\ref{prop:subsets D_conditions}$. Let $\xi_1$, $\xi_2$ be maps from $D$ to $\Cp^{\times}$. There exists a root $\beta_0\in D'$ such that if $\xi_1(\beta_0)\neq\xi_2(\beta_0)$ then $\Omega_{D,\xi_1}\neq\Omega_{D,\xi_2}$.}{Note that if $\beta_1,~\beta_2\in D$ then $\beta_1\notin S(\beta_2)$. Indeed, for $E_6$, $E_7$ and $E_8$ this follows from the fact that $(\beta_1,~\beta_2)\leq0$ for all distinct roots $\beta_1$, $\beta_2\in D$ and the comment before Proposition~\ref{prop:distinct_orbits}. For $F_4$ and $G_2$, this can be checked directly. According to Proposition~\ref{prop:distinct_orbits}, it is enough to check that there exist a simple root $\alpha_0$ and a root $\beta_0\in D'$ such that $(\alpha_0,\beta_0)\neq0$ and $(\alpha_0,\beta)=0$ for all $\beta\in D\setminus D'$. (Indeed, one can pick $\beta_0$ such that $e_{\beta_0}$ is minimal with respect to the total order $\leq_t$ from Proposition~\ref{prop:distinct_orbits} among all $e_{\beta}$, $\beta\in D'$, for which $\beta$ is not orthogonal to $\alpha_0$.) The existence of such roots follows from Proposition~\ref{prop:subsets D_conditions} (iv). For example, consider case 10 for $\Phi=E_6$. Here
\begin{equation*}
\begin{split}
D&=\{\alpha_3+\alpha_4,~\alpha_4+\alpha_5+\alpha_6,~\alpha_1+2\alpha_2+2\alpha_3+3\alpha_4+2\alpha_5+\alpha_6,~\alpha_1+\alpha_3+\alpha_4+\alpha_5\},\\
D'&=\{\alpha_3+\alpha_4,~\alpha_4+\alpha_5+\alpha_6\}.
\end{split}
\end{equation*}
Then $\beta_0=\alpha_3+\alpha_4$ and $\alpha_0=\alpha_4$ satisfy all the conditions of Proposition~\ref{prop:distinct_orbits}, hence $\Omega_{D,\xi_1}\neq\Omega_{D,\xi_2}$.}\newpage

\sect{Proof of the main result}

We are now ready to\label{sect:proof_main} formulate and, using the previous sections, prove our main result, Theorem~\ref{theo:ideal_exceptional} (cf. Theorem~\ref{theo:ideal_classical}). Note that each element of $S(\nt)$ can be considered as a polynomial function on $\nt^*$ via the natural isomorphism $(\nt^*)^*\cong\nt$.

\theop{Suppose $\Phi$ is an irreducible root system of exceptional type\textup, i.e.\textup, $\Phi=E_6$\textup, $E_7$\textup, $E_8$\textup, $F_4$ or $G_2$. The following conditions on a primitive ideal $J\subset U(\nt)$ are equivalent\textup{:}
\begin{equation*}
\begin{split}
&\text{\textup{i)} $J$ is centrally generated};\\
&\text{\textup{ii)} the scalars $c_{\beta}$\textup, $\beta\in\Bu\setminus\Delta$\textup, are nonzero};\\
&\text{\textup{iii)} $J=J(f_{\xi})$ for a Kostant form $f_{\xi}\in\nt^*$}.\\
\end{split}
\end{equation*}
If these\label{theo:ideal_exceptional} conditions are satisfied\textup{,} then the map $\xi$ can be reconstructed from $J$.}{$\mathrm{(ii)}\Longrightarrow\mathrm{(iii)}$. Recall that two $N$-orbits of distinct linear form from $$R=\left\{t=\sum\nolimits_{\beta\in\Bu}t_{\beta}e_{\beta}^*,~t_{\beta}\in\Cp^{\times}\right\}$$ are disjoint, and the union  $X$ of all such orbits is an open dense subset of $\nt^*$; in fact, $X$ is a single $B$-orbit. Given $\beta_i\in\Bu$ and $t\in R$, $$\xi_{\beta_i}(t)=\prod_{\beta\in\Bu}t_{\beta}^{r_{\beta_i}(\beta)},\text{ where }\mu_{\beta_i}=\sum_{\beta\in\Bu}r_{\beta_i}(\beta)\beta,$$ see Remark~\ref{nota:values_of_xi_beta} (iii). Furthermore, from Remark~\ref{nota:values_of_xi_beta} (i), (ii) we see that actually $$\mu_{\beta_i}=\beta_i+\sum_{j<i}r_{\beta_i}(\beta_j)\beta_j,\text{ so }\xi_{\beta_i}(t)=t_{\beta_i}\prod_{j<i}t_{\beta_j}^{r_{\beta_i}(\beta_j)}.$$

We claim that there exists a unique map $\xi\colon\Bu\to\Cp^{\times}$ such that $\xi(\beta)\neq0$ for $\beta\in\Bu\setminus\Delta$ and the Kostant form $f_{\xi}$ satisfies $\xi_{\beta_i}(f_{\xi})=c_{\beta_i}$ for all $\beta_i\in\Bu$. Indeed, since $$f_{\xi}=\sum_{j=1}^m\xi(\beta_j)e_{\beta_j}^*,~m=|\Bu|,$$ belongs to $\overline{R}$, the Zariski closure of $R$ in $\nt^*$, we obtain that $$\xi_{\beta_i}(f_{\xi})=\xi(\beta_i)\prod_{j<i}\xi(\beta_j)^{r_{\beta_i}(\beta_j)}$$ for all $i$ from $1$ to $m$. Since $\xi_{\beta_1}=e_{\beta_1}$, we must set $\xi(\beta_1)=c_{\beta_1}$. Now, assume that $i>1$ and that $\xi_{\beta_j}$ is already defined for all $j<i$ so that $\xi_{\beta_j}(f_{\xi})=c_{\beta_j}$. Then one can put $\xi(\beta_i)$ to be equal to $c_{\beta_i}\prod_{j<i}\xi(\beta_j)^{-r_{\beta_i}(\beta_j)}$, so that $\xi_{\beta_i}(f_{\xi})=c_{\beta_i}$, as required. Note that, according to Remark~\ref{nota:values_of_xi_beta}~(i), if $\xi(\beta_j)=0$ for some $j<i$ then $\beta_j\in\Bu\cap\Delta$, and, consequently, $r_{\beta_i}(\beta_j)=0$. Thus, $\beta_j$ does not actually occurs in the expression of $\mu_{\beta_i}$ and we do not divide by zero in the definition of $\xi(\beta_i)$.

Now, let $\xi\colon\Bu\to\Cp^{\times}$ be such that $\xi_{\beta_i}(f_{\xi})=c_{\beta_i}$ for all $i$. Then, by \cite[6.6.9 (c)]{Dixmier3}, $J(f_{\xi})$ contains $\Delta_{\beta_i}-c_{\beta_i}$ for all $i$, hence $J(f_{\xi})$ contains the centrally generated ideal $J_c$, which is generated by definition by all $\Delta_{\beta_i}-c_{\beta_i}$, $1\leq i\leq m$. But, thanks to Proposition~\ref{prop:J_c_primitive_exceptional}, $J_c$ is primitive. Thus, both $J(f_{\xi})$ and $J_c$ are primitive (and so maximal), hence $J(f_{\xi})=J_c$.

$\mathrm{(iii)}\Longrightarrow\mathrm{(i)}$. Again by \cite[6.6.9 (c)]{Dixmier3},
\begin{equation}
c_{\beta_i}\label{formula:reconstruct_xi_by_c}=\xi_{\beta_i}(f_{\xi})=\xi(\beta_i)\prod_{j<i}\xi(\beta_j)^{r_{\beta_i}(\beta_j)}
\end{equation}
 for all $\beta_i\in\Bu$, hence $c_{\beta_i}\neq0$ for $\beta_i\in\Bu\setminus\Delta$. The ideal $J=J(f_{\xi})$ contains the centrally generated ideal $J_c$, and condition (ii) is satisfied, so $J_c$ is primitive and $J=J(f_{\xi})=J_c$, as required.

\newpage$\mathrm{(i)}\Longrightarrow\mathrm{(ii)}$. This is the most interesting and complicated part of the proof. Let $J$ be a centrally generated primitive ideal of $U(\nt)$. Denote by $c$ the $m$-tuple $c=(c_{\beta_1},~\ldots,~c_{\beta_m})$ and assume that $c_{\beta_i}=0$ for some $\beta_i\in\Bu\setminus\Delta$. To such an $m$-tuple we assign subsets $D'\subset D\subset\Phi^+$ from Proposition~\ref{prop:subsets D_conditions}.

Our next claim is that there exists a map $\eta$ from~$D$ to the field of rational functions $\Cp(x)$ satisfying the following conditions (here, given $z\in\Cp^{\times}$, we denote by $\eta(z)$ the map from~$D$ to $\Cp^{\times}$ defined by $D\ni\gamma\mapsto\eta(\gamma)(z)$):
\begin{equation*}
\begin{split}
\text{i) }&\eta(\gamma)\text{ is non-constant for all }\gamma\in D',\\
\text{ii) }&\eta(\gamma)(z)\text{ is well-defined and nonzero for all }z\neq0,~\gamma\in D,\\
\text{iii) }&\xi_{\beta_i}(f_{D,\eta(z)})=c_{\beta_i}\text{ for all }\beta_i\in\Bu\text{ and }z\neq0.
\end{split}
\end{equation*}

To construct such a map $\eta$, we first note that if $c_{\beta_i}\neq0$ then $\mu_{\beta_i}$ can be uniquely expressed as a $\Zp_{\geq0}$-linear combination of the roots from $D$, i.e., there exist unique $a_{\gamma,\beta_i}\in\Zp_{\geq0}$, $\gamma\in D$, satisfying
\begin{equation*}
\mu_{\beta_i}=\sum\nolimits_{\gamma\in D}a_{\gamma,\beta_i}\gamma,~1\leq i\leq m,~c_{\beta_i}\neq0.
\end{equation*}
Indeed, the uniqueness follows from the linear independence of $D$ (Proposition~\ref{prop:subsets D_conditions} (i)), while the existence is implied by Proposition~\ref{prop:subsets D_conditions} (iii). For instance, if the root system $\Phi$ is of type $E_6$ and
\begin{equation*}
\begin{split}
D=\{\gamma_1&=\alpha_2,~\gamma_2=\alpha_2+\alpha_3+2\alpha_4+\alpha_5,~\gamma_3=\alpha_1+\alpha_2+\alpha_3+2\alpha_4+\alpha_5+\alpha_6,\\
\gamma_4&=\alpha_1+\alpha_2+2\alpha_3+2\alpha_4+2\alpha_5+\alpha_6\}
\end{split}
\end{equation*}
is the $8$th subset from Table~\ref{table:D_D'_E_6} then $$\mu_{\beta_2}=\gamma_3+\gamma_4,~\mu_{\beta_3}=\gamma_2+\gamma_3+2\gamma_4,~\mu_{\beta_4}=2\gamma_2+2\gamma_3+2\gamma_4.$$ This means that, given a map $\xi\colon D\to\Cp^{\times}$, the value of at most one monomial of $\xi_{\beta_i}$ at the linear form $f_{D,\xi}$ is non-zero (precisely, the monomial $e_D=\prod_{\gamma\in D}e_{\gamma}^{a_{\gamma,\beta_i}}$). Indeed, the value of $f_{D,\xi}$ on a monomial is nonzero if and only if all variables involved in this monomial has the form $e_{\gamma}$ for $\gamma\in D$, but all monomials of $\xi_{\beta_i}$ have weight~$\mu_{\beta_i}$. (We see that if $c_{\beta_j}=0$ then $\xi_{\beta_j}(f_{D,\xi})=0$ by Proposition~\ref{prop:subsets D_conditions}~(ii).) Note, however, that a priori we do not know that the monomial $e_D$ in fact occurs in $\xi_{\beta_i}$.

To check that $e_D$ occurs in $\xi_{\beta_i}$, we consider the subset
\begin{equation*}
D_1=\begin{cases}D,&\text{if }\beta_1\in D,\\
D\cup\{\beta_1\},&\text{if }\beta_1\notin D.
\end{cases}
\end{equation*}
(Actually, $\beta_1\in D$ if and only if $c_{\beta_1}\neq0$.) Given a map $\xi\colon D\to\Cp^{\times}$, we consider a map $\xi_1\colon D_1\to\Cp^{\times}$ such that $\xi_1(\gamma)=\xi(\gamma)$ for $\gamma\in D$, and the linear form $f_1$, where
\begin{equation*}
f_1=f_{D_1,\xi_1}=\begin{cases}
f_{D,\xi},&\text{if $\beta_1\in D$},\\
f_{D,\xi}+\xi_1(\beta_1)e_{\beta_1}^*,&\text{if $\beta_1\notin D$.}
\end{cases}
\end{equation*}
Clearly, $e_D(f_{D,\xi})=e_D(f_1)$. On the other hand, one can easily check that the condition $f_1(e_{\beta_1})\neq0$ implies that there exists a linear form $\lambda$ in the coadjoint $N$-orbit of $f_1$ such that $\lambda(e_{\beta})\neq0$ for all $\beta\in\Bu$. By Remark~\ref{nota:values_of_xi_beta} (iii), $e_{\Bu}=\prod_{\beta\in\Bu}e_{\beta}^{r_{\beta_i}(\beta)}$ enters $\xi_{\beta_i}$ with coefficient~1, and $e_{\Bu}(\lambda)\neq0$. But $\lambda=g.f_1$ for a certain $g\in N$, where $g.f_1$ denotes the result of the coadjoint action of $N$ on $\nt^*$. The adjoint action of $N$ on the algebra $S(\nt)$ has the form $$(z.s)(\mu)=s(z^{-1}.\mu),~z\in N,~s\in S(\nt),~\mu\in\nt^*.$$ We see that $(g^{-1}.e_{\Bu})(f_1)=e_{\Bu}(\lambda)\neq0$. But $\xi_{\beta_i}$ in $N$-invariant, so the monomial $g^{-1}.e_{\Bu}$ occurs in $\xi_{\beta}$ (with coefficient 1). Thus, the latter monomial coincides with $c_D^ie_D$ for certain $c_D^i\in\Cp^{\times}$.\newpage

Next, we note that the (affine) solution space for the system of linear equations
\begin{equation}
\sum\nolimits_{\gamma\in D}a_{\gamma,\beta_i}y_{\gamma}=b_i,~1\leq i\leq m,~c_{\beta_i}\neq0.\label{formula:system_for_eta}
\end{equation}
is at least one-dimensional for all possible $b_i\in\Cp$. Indeed, $k=|D|$ is in fact not less than the number of non-zero scalars in $c$ plus one, so this system on $k$ variables contains at most $(k-1)$ equations. It follows from Panov's description of the weights $\mu_{\beta_i}$ \cite[p. 8]{Panov2} that the equations are linearly independent, and the rank of the system is at most $(k-1)$. Hence the solution space is at least one-dimensional, as required. Recall that $D'$ is exactly the set of roots $\gamma\in D$ for which the solution space is not orthogonal to the axis $y_{\gamma}$.

We are ready to construct a map $\eta$ satisfying the above conditions. Assume that such a map is already constructed. Then, given $z\in\Cp^{\times}$,
\begin{equation*}
\begin{split}
\xi_{\beta_i}(f_{D,\eta(z)})&=c_D^ie_D(f_{D,\eta(z)})=c_D^i\left(\prod_{\gamma\in D}e_{\gamma}^{a_{\gamma,\beta_i}}\right)(f_{D,\eta(z)})\\
&=c_D^i\prod_{\gamma\in D}(\eta(\gamma)(z))^{a_{\gamma,\beta_i}}=c_{\beta_i},~1\leq i\leq m.
\end{split}
\end{equation*}
Let $b_i=c_{\beta_i}-c_D^i$, and $y=(y_{\gamma})_{\gamma\in D}$ be a solution of system (\ref{formula:system_for_eta}). Then, clearly, $\prod_{\gamma\in D}x_{\gamma}^{a_{\gamma,\beta_i}}=c_{\beta_i}\left(c_D^i\right)^{-1}$ for all $i$, where $x_{\gamma}=\exp(y_{\gamma})$. Let $\{y_{\gamma}(t)\}_{\gamma\in D}$ be a one-dimensional affine subspace of the space of solutions of (\ref{formula:system_for_eta}) not orthogonal to the axes $y_{\gamma}$ for $\gamma\in D'$. Note that each $y_{\gamma}(t)$, $\gamma\in D$, is an affine-linear function in $t$. We set $x=\exp(t)$ and $\eta(\gamma)(x)=\exp(y_{\gamma}(t))$ for $\gamma\in D$, then $\eta(\gamma)$ is a rational function in $x$. Condition~(iii) is satisfied by the definition of $\eta$, condition~(ii) is obvious. Finally, condition~(i) is satisfied because the subspace $\{y_{\gamma}(t)\}_{\gamma\in D}$ is not orthogonal to the axis $y_{\gamma}$ for $\gamma\in D'$.

For example, if $\Phi=E_6$ and $D=\{\gamma_1,~\gamma_2,~\gamma_3,~\gamma_4\}$ is, as above, the 8th subset from Table~\ref{table:D_D'_E_6}, then one can put
\begin{equation*}
\begin{split}
\eta(\gamma_1)&=x,~\eta(\gamma_2)=c_D^2c_{\beta_2}^{-1}\left(c_D^4\right)^{-1/2}c_{\beta_4}^{1/2},\\
\eta(\gamma_3)&=
c_{\beta_2}\left(c_D^2\right)^{-1}c_D^3c_{\beta_3}^{-1}c_{\beta_4}^{1/2}\left(c_D^4\right)^{-1/2},~
\eta(\gamma_4)=c_{\beta_3}\left(c_D^3\right)^{-1}c_{\beta_4}^{-1/2}\left(c_D^4\right)^{1/2}.
\end{split}
\end{equation*}
(Here, given $a\in\Cp$, we denote by $a^{1/2}$ a complex number such that $(a^{1/2})^2=a$.) Another example: if $\Phi=E_6$ and $D=\{\gamma_1=\alpha_3+\alpha_4+\alpha_5+\alpha_6,~\gamma_2=\alpha_1+2\alpha_2+2\alpha_3+3\alpha_4+2\alpha_5+\alpha_6,~\gamma_3=\alpha_1+\alpha_3+\alpha_4+\alpha_5\}$ is the 11th subset from Table~\ref{table:D_D'_E_6}, then $D'=\{\gamma_1,~\gamma_3\}$ and
\begin{equation*}
\begin{split}
&\mu_{\beta_1}=\gamma_2,~\mu_{\beta_3}=\gamma_1+2\gamma_2+\gamma_3, ~\eta(\gamma_1)=x,~\eta(\gamma_2)=c_{\beta_1}\left(c_D^1\right)^{-1},~\eta(\gamma_3)=c_{\beta_3}\left(c_D^3\right)^{-1}\left(c_D^1\right)^2c_{\beta_1}^{-2}x^{-1}.
\end{split}
\end{equation*}

In general, let $\eta\colon D\to\Cp(x)$ be a map satisfying conditions (i), (ii), (iii). Since $\xi_{\beta_i}(f_{D,\eta(z)})=c_{\beta_i}$ for all $i$ and all $z\in\Cp^{\times}$, one has $J\subset J(f_{D,\eta(z)})$ for all $z\in\Cp^{\times}$. But both $J$ and $J(f_{D,\eta(z)})$ are maximal, hence $J=J(f_{D,\eta(z)})$ for all $z\in\Cp^{\times}$. Pick a root $\beta_0\in D'$ satisfying Proposition~\ref{prop:distinc_orbits_precise}. Since $\eta(\beta_0)$ is non-constant, there exist $z_1,~z_2\in\Cp^{\times}$ such that $\eta(\beta_0)(z_1)\neq\eta(\beta_0)(z_2)$. By Proposition~\ref{prop:distinc_orbits_precise}, the coadjoint orbits of $f_{D,\eta(z_1)}$ and $f_{D,\eta(z_2)}$ do not coincide. Hence $J(f_{D,\eta(z_1)})\neq J(f_{D,\eta(z_2)})$, a contradiction.

Finally, if conditions (i)--(iii) are satisfied then formula (\ref{formula:reconstruct_xi_by_c}) and the equality $\Delta_{\beta_1}=e_{\beta_1}$ together imply that the map $\xi$ can be reconstructed from $J$. The proof is complete.}

As an immediate corollary, we obtain that a similar result is true for all (probably, reducible) root systems.
\theop{Let $\Phi$ be an arbitrary root system. The following conditions on a primitive ideal $J\subset U(\nt)$ are equivalent\textup{:}
\begin{equation*}
\begin{split}
&\text{\textup{i)} $J$ is centrally generated};\\
&\text{\textup{ii)} the scalars $c_{\beta}$\textup, $\beta\in\Bu\setminus\Delta$\textup, are nonzero};\\
&\text{\textup{iii)} $J=J(f_{\xi})$ for a Kostant form $f_{\xi}\in\nt^*$}.\\
\end{split}
\end{equation*}
If these\label{theo:ideal_all} conditions are satisfied\textup{,} then the map $\xi$ can be reconstructed from $J$.\newpage}{Let $\Phi_i$, $1\leq i\leq k$, be the irreducible components of the root system $\Phi$, and let $\nt=\bigoplus_{i=1}^k\nt_k$ be the corresponding division of $\nt$ into a direct sum of its nilpotent ideals, so $U(\nt)=U(\nt_1)\otimes\ldots\otimes U(\nt_k)$ as associative algebras.

$\mathrm{(ii)}\Longrightarrow\mathrm{(iii)}$. Denote $\Delta_i=\Delta\cap\Phi_i$ and $\Bu_i=\Bu\cap\Phi_i$ for $1\leq i\leq k$, so $\Delta_i$ is a basis for $\Phi_i$ such that $\Phi_i^+=\Phi^+\cap\Phi_i$, and $\Bu_i$ is the Kostant cascade in $\Phi_i^+$. Put also $J_i=J\cap U(\nt_i)$, so $J_i$ is an ideal of $U(\nt_i)$ containing all $\Delta_{\beta}-c_{\beta}$, $\beta\in\Bu_i$. Since $c_{\beta}\neq0$ for $\beta\in\Bu_i\setminus\Delta_i$, Theorem~\ref{theo:ideal_exceptional}, \cite[Theorem 3.1]{IgnatyevPenkov1} and \cite[Theorem 2.4]{Ignatyev2} imply that $J_i$ is a primitive ideal of $U(\nt_i)$. Furthermore, for each $i$, there exists a map $\xi_i\colon \Bu_i\to\Cp$ such that $f_{\xi_i}$ is a Kostant form on $\nt_i$ and $J_i=J(f_{\xi_i})$. Define $\xi$ to be the map from $\Bu$ to $\Cp$ such that $\xi(\beta)=\xi_i(\beta)$ for $\beta\in\Bu_i$, so $f_{\xi}$ is a Kostant form on $\nt$. Now, let $J_c$ be, as above, the ideal of $U(\nt)$ generated by $\Delta_{\beta}-c_{\beta}$, $\beta\in\Bu$. Then $J_c$ is contained both in $J$ and in $J(f_{\xi})$. But the quotient algebra $U(\nt_i)/J_i$ is isomorphic to the Weyl algebra $\Au_{s_i}$ for certain $s_i\geq1$. Thus, $$U(\nt)/J\cong\Au_{s_1}\otimes\ldots\otimes\Au_{s_k}\cong\Au_s,$$ where $s=s_1+\ldots+s_k$, because $\Au_a\otimes\Au_b\cong\Au_{a+b}$. It follows that $J_c$ is primitive (and so maximal), hence $J=J_c=J(f_{\xi})$.

$\mathrm{(iii)}\Longrightarrow\mathrm{(i)}$. Denote by $\xi_i$ the restriction of $\xi$ to $\Bu_i$, $1\leq i\leq k$. Again by Theorem~\ref{theo:ideal_exceptional},\break \cite[Theorem~3.1]{IgnatyevPenkov1} and \cite[Theorem 2.4]{Ignatyev2}, each $J(f_{\xi_i})$ is a centrally generated primitive ideal of $U(\nt_i)$, and $c_{\beta}\neq0$ for all $\beta\in\Bu_i\setminus\Delta_i$ (and so for all $\beta\in\Bu\setminus\Delta$). And we see again that the centrally generated ideal $J_c$ is primitive and at the same time is contained in $J=J(f_{\xi})$, thus, $J=J_c$.

$\mathrm{(i)}\Longrightarrow\mathrm{(ii)}$. Assume that there exists $\beta\in\Bu\setminus\Delta$ such that $c_{\beta}=0$. Let $A$ be the set of all indices~$a$ between 1 and $k$ such that there exists $\beta\in\Bu_a\setminus\Delta_a$ for which $c_{\beta}=0$. It follows from the proofs of Theorem~\ref{theo:ideal_exceptional}, \cite[Theorem 3.1]{IgnatyevPenkov1} and \cite[Theorem 2.4]{Ignatyev2} that, given $a\in A$, there exist $D_a\subset\Phi_j^+$ and distinct maps $\xi_a^1$, $\xi_a^2$ from $D_a$ to $\Cp^{\times}$ such that $$\xi_{\beta}(f_{D_a,\xi_a^1})=\xi_{\beta}(f_{D_a,\xi_a^2})=c_{\beta}$$ for all $\beta\in\Bu_a$, and the coadjoint orbits of $f_{D_a,\xi_a^1}$ and $f_{D_a,\xi_a^2}$ (in $\nt_a^*$) are distinct.

One the other hand, if $i\in\{1,~\ldots,~k\}\setminus A$ then $c_{\beta}\neq0$ for all $\beta\in\Bu_i\setminus\Delta_i$. Hence, as above, $J\cap U(\nt_i)=J(f_{\xi_i})$ for a certain map $\xi_i\colon D\to\Cp$ (in other words, $f_{\xi_i}$ is a Kostant form on $\nt_i$ and $\xi_{\beta}(f_{\xi_i})=c_{\beta}$ for all $\beta\in\Bu_i$). Clearly, if $f\in\nt^*$ is a linear form on $\nt$ and $f_j$ is its restriction to $\nt_j$, $1\leq j\leq k$, then $$\Omega_f=\Omega_{f_1}\times\ldots\times\Omega_{f_k},$$ where $\Omega_f$ (respectively, $\Omega_{f_j}$) is the coadjoint orbit of the linear form $f$ in $\nt^*$ (respectively, of the linear form $f_j$ in~$\nt_j^*$). Now, put $$D=\bigcup_{a\in A}D_a\cup\bigcup_{i\notin A}\Bu_i$$ and define $\xi^j\colon D\to\Cp$, $j=1,2$, by the rule
\begin{equation*}
\xi^j(\beta)=\begin{cases}\xi_a^j(\beta),&\text{if }\beta\in\Bu_a\text{ for }a\in A,\\
\xi_i(\beta),&\text{if }\beta\in\Bu_i\text{ for }i\notin A.
\end{cases}
\end{equation*}
Put $f^j=f_{D,\xi^j}$ for $j=1,~2$, so $$\Omega_{f^j}=\prod_{a\in A}\Omega_{D_a,\xi_a^j}\times\prod_{i\notin A}\Omega_{\xi_i},$$ where $\Omega_{\xi_i}$ is the coadjoint orbit of the Kostant form $f_{\xi_i}$ in $\nt_i^*$. Since $\Omega_{D_a,\xi_a^1}\neq\Omega_{D_a,\xi_a^2}$ for at least one $a\in A$, one has $\Omega_{f^1}\neq\Omega_{f^2}$, so $J(f^1)\neq J(f^2)$. At the same time, both maximal ideals contain the maximal ideal $J$, a contradiction.

Finally, if conditions (i)--(iii) are satisfied then the map $\xi$ can be reconstructed from $J$, because the restriction of $\xi$ to $\Bu_i$ can be reconstructed from $J_i$ for all $i$. The proof is complete.}

\textsc{Supplementary material}. The online version of this article contains additional supplementary material (analogues of Tables \ref{table:D_D'_E_6}, \ref{table:D_D'_F_4}, \ref{table:c_D_E_6_F_4}, \ref{table:c_D_D'_G_2} for the root systems $E_7$, $E_8$ and the listing of the code generating these analogues). Please visit \texttt{https://arxiv.org/pdf/1907.04219.pdf}.

\medskip\textsc{Mikhail V. Ignatyev: Samara National Research University, Ak. Pavlova 1, 443011,\\\indent Samara, Russia}

\emph{E-mail address}: \texttt{mihail.ignatev@gmail.com}

\medskip\textsc{Aleksandr A. Shevchenko: Samara National Research University, Ak. Pavlova 1,\\\indent 443011, Samara, Russia}

\emph{E-mail address}: \texttt{shevchenko.alexander.1618@gmail.com }

\newpage
\begin{center}
\textbf{Appendix A}
\end{center}
\vspace{-0.7cm}\begin{center}
	\tabl{List of the subsets $D$ and $D'$ for $E_7$}
\end{center}
\vspace{-0.7cm}
\begin{longtable}{||l|l|l||l|l|l||l|l|l||}
	\hline
	&$D$&$D'$&&$D$&$D'$&&$D$&$D'$\\
	\hline\hline
1&\begin{tabular}[t]{l}
	$0100000$\end{tabular}&\begin{tabular}[t]{l}
	$0100000$\end{tabular}&2&\begin{tabular}[t]{l}
	$0101111$,\\$1234321$,\\$1011111$\end{tabular}&\begin{tabular}[t]{l}
	$0101111$,\\$1234321$,\\$1011111$\end{tabular}&3&\begin{tabular}[t]{l}
	$0101110$,\\$1112211$,\\$1234321$,\\$1011110$\end{tabular}&\begin{tabular}[t]{l}
	$0101110$,\\$1112211$,\\$1234321$,\\$1011110$\end{tabular}\\\hline4&\begin{tabular}[t]{l}
	$0101111$,\\$1112210$,\\$1234321$,\\$1011111$\end{tabular}&\begin{tabular}[t]{l}
	$0101111$,\\$1112210$,\\$1234321$,\\$1011111$\end{tabular}&5&\begin{tabular}[t]{l}
	$1223210$,\\$1122111$,\\$1123321$\end{tabular}&\begin{tabular}[t]{l}
	$1223210$,\\$1122111$,\\$1123321$\end{tabular}&6&\begin{tabular}[t]{l}
	$0001111$,\\$1223211$,\\$1122221$,\\$1123210$\end{tabular}&\begin{tabular}[t]{l}
	$0001111$,\\$1223211$,\\$1122221$,\\$1123210$\end{tabular}\\\hline7&\begin{tabular}[t]{l}
	$0001110$,\\$1223321$,\\$1122110$,\\$1123211$\end{tabular}&\begin{tabular}[t]{l}
	$0001110$,\\$1223321$,\\$1122110$,\\$1123211$\end{tabular}&8&\begin{tabular}[t]{l}
	$0001111$,\\$1223321$,\\$1122111$,\\$1123210$\end{tabular}&\begin{tabular}[t]{l}
	$0001111$,\\$1223321$,\\ \\$1123210$\end{tabular}&9&\begin{tabular}[t]{l}
	$0101000$,\\$1111110$,\\$1112211$,\\$1234321$\end{tabular}&\begin{tabular}[t]{l}
	$0101000$,\\$1111110$,\\$1112211$,\\$1234321$\end{tabular}\\\hline10&\begin{tabular}[t]{l}
	$0111111$,\\$1223210$,\\$1112211$,\\$1123221$\end{tabular}&\begin{tabular}[t]{l}
	$0111111$,\\ \\$1112211$,\\$1123221$\end{tabular}&11&\begin{tabular}[t]{l}
	$0111110$,\\$1223211$,\\$1112110$,\\$1123321$\end{tabular}&\begin{tabular}[t]{l}
	$0111110$,\\$1223211$,\\$1112110$\\ \end{tabular}&12&\begin{tabular}[t]{l}
	$0111111$,\\$1223210$,\\$1112111$,\\$1123321$\end{tabular}&\begin{tabular}[t]{l}
	$0111111$,\\ \\$1112111$\\ \end{tabular}\\
\hline
\end{longtable}

\begin{center}
	\tabl{List of the subsets $D$ and $D'$ for $E_8$}
\end{center}
\vspace{-0.7cm}{\small\begin{longtable}{||l|l|l||l|l|l||l|l|l||}
	\hline
	&$D$&$D'$&&$D$&$D'$&&$D$&$D'$\\
	\hline\hline
	1&\begin{tabular}[t]{l}
		$01000000$\end{tabular}&\begin{tabular}[t]{l}
		$01000000$\end{tabular}&2&\begin{tabular}[t]{l}
		$23354321$,\\$12343221$,\\$11233321$\end{tabular}&\begin{tabular}[t]{l}
		$23354321$,\\$12343221$,\\$11233321$\end{tabular}&3&\begin{tabular}[t]{l}
		$23354321$,\\$12344321$,\\$11232110$,\\$11232211$\end{tabular}&\begin{tabular}[t]{l}
		$23354321$,\\$12344321$,\\$11232110$,\\$11232211$\end{tabular}\\\hline4&\begin{tabular}[t]{l}
		$23354321$,\\$12344321$,\\$11232100$,\\$11232221$\end{tabular}&\begin{tabular}[t]{l}
		$23354321$,\\$12344321$,\\$11232100$,\\$11232221$\end{tabular}&5&\begin{tabular}[t]{l}
		$23454321$,\\$12243211$,\\$11232110$,\\$11233321$\end{tabular}&\begin{tabular}[t]{l}
		$23454321$,\\$12243211$,\\$11232110$,\\$11233321$\end{tabular}&6&\begin{tabular}[t]{l}
		$23454321$,\\$12243221$,\\$11232100$,\\$11233321$\end{tabular}&\begin{tabular}[t]{l}
		$23454321$,\\$12243221$,\\$11232100$,\\$11233321$\end{tabular}\\\hline7&\begin{tabular}[t]{l}
		$23454321$,\\$12244321$,\\$11232110$,\\$11232211$\end{tabular}&\begin{tabular}[t]{l}
		$23454321$,\\ \\$11232110$,\\$11232211$\end{tabular}&8&\begin{tabular}[t]{l}
		$23454321$,\\$12244321$,\\$11232100$,\\$11232221$\end{tabular}&\begin{tabular}[t]{l}
		$23454321$,\\ \\$11232100$,\\$11232221$\end{tabular}&9&\begin{tabular}[t]{l}
		$22354321$,\\$12232110$,\\$12233211$,\\$12343321$\end{tabular}&\begin{tabular}[t]{l}
		$22354321$,\\$12232110$,\\$12233211$,\\$12343321$\end{tabular}\\
\hline
\end{longtable}}

\begin{center}
	\tabl{Correspondence between $c$ and $D$ for $E_7$}
\end{center}

{\tabskip0pt\relax
\def\hline{\noalign{\hrule\allowbreak\vskip-.4pt\hrule}}
\def\Hline{\noalign{\hrule\kern2pt\hrule}}
\halign to \textwidth{\strut\raise1pt\copy\strutbox\lower1pt\copy\strutbox
\vrule\kern2pt\vrule#\unskip\tabskip0pt plus 1fil\relax
&\hfil\ignorespaces#\unskip&\vrule#\unskip
&\ignorespaces#\unskip\hfil&\vrule\kern2pt\vrule#\unskip
&\hfil\ignorespaces#\unskip&\vrule#\unskip
&\ignorespaces#\unskip\hfil&\vrule\kern2pt\vrule#\unskip
&\hfil\ignorespaces#\unskip&\vrule#\unskip
&\ignorespaces#\unskip\hfil&\vrule\kern2pt\vrule#\unskip\tabskip0pt\relax\cr\hline
&$D$\hfil\null&&\hfil Type of $c$&&$D$\hfil\null&&\hfil Type of $c$&&$D$\hfil\null&&\hfil Type of $c$&\cr\Hline
&1&&
$(0,0,0,0,0,0,0)$&&
2&&
$(0,0,0,c_{\beta_{4}},0,0,0)$&&
3&&
$(0,0,0,0,0,0,c_{\beta_{7}})$&\cr\hline&
4&&
$(0,0,0,c_{\beta_{4}},0,0,c_{\beta_{7}})$&&
5&&
$(0,0,0,0,0,c_{\beta_{6}},0)$&&
6&&
$(0,0,0,c_{\beta_{4}},0,c_{\beta_{6}},0)$&\cr\hline&
7&&
$(0,0,0,0,0,c_{\beta_{6}},c_{\beta_{7}})$&&
8&&
$(0,0,0,c_{\beta_{4}},0,c_{\beta_{6}},c_{\beta_{7}})$&&
9&&
$(0,0,0,0,c_{\beta_{5}},0,0)$&\cr\hline&
10&&
$(0,0,0,c_{\beta_{4}},c_{\beta_{5}},0,0)$&&
11&&
$(0,0,0,0,c_{\beta_{5}},0,c_{\beta_{7}})$&&
12&&
$(0,0,0,c_{\beta_{4}},c_{\beta_{5}},0,c_{\beta_{7}})$&\cr\hline&
13&&
$(0,0,0,0,c_{\beta_{5}},c_{\beta_{6}},0)$&&
14&&
$(0,0,0,c_{\beta_{4}},c_{\beta_{5}},c_{\beta_{6}},0)$&&
15&&
$(0,0,0,0,c_{\beta_{5}},c_{\beta_{6}},c_{\beta_{7}})$&\cr\hline&
16&&
$(0,0,0,c_{\beta_{4}},c_{\beta_{5}},c_{\beta_{6}},c_{\beta_{7}})$&&
17&&
$(0,0,c_{\beta_{3}},0,0,0,0)$&&
18&&
$(0,0,c_{\beta_{3}},c_{\beta_{4}},0,0,0)$&\cr\hline&
19&&
$(0,0,c_{\beta_{3}},0,0,0,c_{\beta_{7}})$&&
20&&
$(0,0,c_{\beta_{3}},c_{\beta_{4}},0,0,c_{\beta_{7}})$&&
21&&
$(0,0,c_{\beta_{3}},0,0,c_{\beta_{6}},0)$&\cr\hline&
22&&
$(0,0,c_{\beta_{3}},c_{\beta_{4}},0,c_{\beta_{6}},0)$&&
23&&
$(0,0,c_{\beta_{3}},0,0,c_{\beta_{6}},c_{\beta_{7}})$&&
24&&
$(0,0,c_{\beta_{3}},c_{\beta_{4}},0,c_{\beta_{6}},c_{\beta_{7}})$&\cr\hline&
25&&
$(0,0,c_{\beta_{3}},0,c_{\beta_{5}},0,0)$&&
26&&
$(0,0,c_{\beta_{3}},c_{\beta_{4}},c_{\beta_{5}},0,0)$&&
27&&
$(0,0,c_{\beta_{3}},0,c_{\beta_{5}},0,c_{\beta_{7}})$&\cr\hline&
28&&
$(0,0,c_{\beta_{3}},c_{\beta_{4}},c_{\beta_{5}},0,c_{\beta_{7}})$&&
29&&
$(0,0,c_{\beta_{3}},0,c_{\beta_{5}},c_{\beta_{6}},0)$&&
30&&
$(0,0,c_{\beta_{3}},c_{\beta_{4}},c_{\beta_{5}},c_{\beta_{6}},0)$&\cr\hline&
31&&
$(0,0,c_{\beta_{3}},0,c_{\beta_{5}},c_{\beta_{6}},c_{\beta_{7}})$&&
32&&
$(0,0,c_{\beta_{3}},c_{\beta_{4}},c_{\beta_{5}},c_{\beta_{6}},c_{\beta_{7}})$&&
33&&
$(0,c_{\beta_{2}},0,0,0,0,0)$&\cr\hline&
34&&
$(0,c_{\beta_{2}},0,c_{\beta_{4}},0,0,0)$&&
35&&
$(0,c_{\beta_{2}},0,0,0,0,c_{\beta_{7}})$&&
36&&
$(0,c_{\beta_{2}},0,c_{\beta_{4}},0,0,c_{\beta_{7}})$&\cr\hline&
37&&
$(0,c_{\beta_{2}},0,0,0,c_{\beta_{6}},0)$&&
38&&
$(0,c_{\beta_{2}},0,c_{\beta_{4}},0,c_{\beta_{6}},0)$&&
39&&
$(0,c_{\beta_{2}},0,0,0,c_{\beta_{6}},c_{\beta_{7}})$&\cr\hline&
40&&
$(0,c_{\beta_{2}},0,c_{\beta_{4}},0,c_{\beta_{6}},c_{\beta_{7}})$&&
41&&
$(0,c_{\beta_{2}},0,0,c_{\beta_{5}},0,0)$&&
42&&
$(0,c_{\beta_{2}},0,c_{\beta_{4}},c_{\beta_{5}},0,0)$&\cr\hline&
43&&
$(0,c_{\beta_{2}},0,0,c_{\beta_{5}},0,c_{\beta_{7}})$&&
44&&
$(0,c_{\beta_{2}},0,c_{\beta_{4}},c_{\beta_{5}},0,c_{\beta_{7}})$&&
45&&
$(0,c_{\beta_{2}},0,0,c_{\beta_{5}},c_{\beta_{6}},0)$&\cr\hline&
46&&
$(0,c_{\beta_{2}},0,c_{\beta_{4}},c_{\beta_{5}},c_{\beta_{6}},0)$&&
47&&
$(0,c_{\beta_{2}},0,0,c_{\beta_{5}},c_{\beta_{6}},c_{\beta_{7}})$&&
48&&
$(0,c_{\beta_{2}},0,c_{\beta_{4}},c_{\beta_{5}},c_{\beta_{6}},c_{\beta_{7}})$&\cr\hline&
49&&
$(0,c_{\beta_{2}},c_{\beta_{3}},0,0,0,0)$&&
50&&
$(0,c_{\beta_{2}},c_{\beta_{3}},c_{\beta_{4}},0,0,0)$&&
51&&
$(0,c_{\beta_{2}},c_{\beta_{3}},0,0,0,c_{\beta_{7}})$&\cr\hline&
52&&
$(0,c_{\beta_{2}},c_{\beta_{3}},c_{\beta_{4}},0,0,c_{\beta_{7}})$&&
53&&
$(0,c_{\beta_{2}},c_{\beta_{3}},0,0,c_{\beta_{6}},0)$&&
54&&
$(0,c_{\beta_{2}},c_{\beta_{3}},c_{\beta_{4}},0,c_{\beta_{6}},0)$&\cr\hline&
55&&
$(0,c_{\beta_{2}},c_{\beta_{3}},0,0,c_{\beta_{6}},c_{\beta_{7}})$&&
56&&
$(0,c_{\beta_{2}},c_{\beta_{3}},c_{\beta_{4}},0,c_{\beta_{6}},c_{\beta_{7}})$&&
57&&
$(0,c_{\beta_{2}},c_{\beta_{3}},0,c_{\beta_{5}},0,0)$&\cr\hline&
58&&
$(0,c_{\beta_{2}},c_{\beta_{3}},c_{\beta_{4}},c_{\beta_{5}},0,0)$&&
59&&
$(0,c_{\beta_{2}},c_{\beta_{3}},0,c_{\beta_{5}},0,c_{\beta_{7}})$&&
60&&
$(0,c_{\beta_{2}},c_{\beta_{3}},c_{\beta_{4}},c_{\beta_{5}},0,c_{\beta_{7}})$&\cr\hline&
61&&
$(0,c_{\beta_{2}},c_{\beta_{3}},0,c_{\beta_{5}},c_{\beta_{6}},0)$&&
62&&
$(0,c_{\beta_{2}},c_{\beta_{3}},c_{\beta_{4}},c_{\beta_{5}},c_{\beta_{6}},0)$&&
63&&
$(0,c_{\beta_{2}},c_{\beta_{3}},0,c_{\beta_{5}},c_{\beta_{6}},c_{\beta_{7}})$&\cr\hline&
64&&
$(0,c_{\beta_{2}},c_{\beta_{3}},c_{\beta_{4}},c_{\beta_{5}},c_{\beta_{6}},c_{\beta_{7}})$&&
65&&
$(c_{\beta_{1}},0,0,0,0,0,0)$&&
66&&
$(c_{\beta_{1}},0,0,c_{\beta_{4}},0,0,0)$&\cr\hline&
67&&
$(c_{\beta_{1}},0,0,0,0,0,c_{\beta_{7}})$&&
68&&
$(c_{\beta_{1}},0,0,c_{\beta_{4}},0,0,c_{\beta_{7}})$&&
69&&
$(c_{\beta_{1}},0,0,0,0,c_{\beta_{6}},0)$&\cr\hline&
70&&
$(c_{\beta_{1}},0,0,c_{\beta_{4}},0,c_{\beta_{6}},0)$&&
71&&
$(c_{\beta_{1}},0,0,0,0,c_{\beta_{6}},c_{\beta_{7}})$&&
72&&
$(c_{\beta_{1}},0,0,c_{\beta_{4}},0,c_{\beta_{6}},c_{\beta_{7}})$&\cr\hline&
73&&
$(c_{\beta_{1}},0,0,0,c_{\beta_{5}},0,0)$&&
74&&
$(c_{\beta_{1}},0,0,c_{\beta_{4}},c_{\beta_{5}},0,0)$&&
75&&
$(c_{\beta_{1}},0,0,0,c_{\beta_{5}},0,c_{\beta_{7}})$&\cr\hline&
76&&
$(c_{\beta_{1}},0,0,c_{\beta_{4}},c_{\beta_{5}},0,c_{\beta_{7}})$&&
77&&
$(c_{\beta_{1}},0,0,0,c_{\beta_{5}},c_{\beta_{6}},0)$&&
78&&
$(c_{\beta_{1}},0,0,c_{\beta_{4}},c_{\beta_{5}},c_{\beta_{6}},0)$&\cr\hline&
79&&
$(c_{\beta_{1}},0,0,0,c_{\beta_{5}},c_{\beta_{6}},c_{\beta_{7}})$&&
80&&
$(c_{\beta_{1}},0,0,c_{\beta_{4}},c_{\beta_{5}},c_{\beta_{6}},c_{\beta_{7}})$&&
81&&
$(c_{\beta_{1}},0,c_{\beta_{3}},0,0,0,0)$&\cr\hline&
82&&
$(c_{\beta_{1}},0,c_{\beta_{3}},c_{\beta_{4}},0,0,0)$&&
83&&
$(c_{\beta_{1}},0,c_{\beta_{3}},0,0,0,c_{\beta_{7}})$&&
84&&
$(c_{\beta_{1}},0,c_{\beta_{3}},c_{\beta_{4}},0,0,c_{\beta_{7}})$&\cr\hline&
85&&
$(c_{\beta_{1}},0,c_{\beta_{3}},0,0,c_{\beta_{6}},0)$&&
86&&
$(c_{\beta_{1}},0,c_{\beta_{3}},c_{\beta_{4}},0,c_{\beta_{6}},0)$&&
87&&
$(c_{\beta_{1}},0,c_{\beta_{3}},0,0,c_{\beta_{6}},c_{\beta_{7}})$&\cr\hline&
88&&
$(c_{\beta_{1}},0,c_{\beta_{3}},c_{\beta_{4}},0,c_{\beta_{6}},c_{\beta_{7}})$&&
89&&
$(c_{\beta_{1}},0,c_{\beta_{3}},0,c_{\beta_{5}},0,0)$&&
90&&
$(c_{\beta_{1}},0,c_{\beta_{3}},c_{\beta_{4}},c_{\beta_{5}},0,0)$&\cr\hline&
91&&
$(c_{\beta_{1}},0,c_{\beta_{3}},0,c_{\beta_{5}},0,c_{\beta_{7}})$&&
92&&
$(c_{\beta_{1}},0,c_{\beta_{3}},c_{\beta_{4}},c_{\beta_{5}},0,c_{\beta_{7}})$&&
93&&
$(c_{\beta_{1}},0,c_{\beta_{3}},0,c_{\beta_{5}},c_{\beta_{6}},0)$&\cr\hline&
94&&
$(c_{\beta_{1}},0,c_{\beta_{3}},c_{\beta_{4}},c_{\beta_{5}},c_{\beta_{6}},0)$&&
95&&
$(c_{\beta_{1}},0,c_{\beta_{3}},0,c_{\beta_{5}},c_{\beta_{6}},c_{\beta_{7}})$&&
96&&
$(c_{\beta_{1}},0,c_{\beta_{3}},c_{\beta_{4}},c_{\beta_{5}},c_{\beta_{6}},c_{\beta_{7}})$&\cr\hline&
97&&
$(c_{\beta_{1}},c_{\beta_{2}},0,0,0,0,0)$&&
98&&
$(c_{\beta_{1}},c_{\beta_{2}},0,c_{\beta_{4}},0,0,0)$&&
99&&
$(c_{\beta_{1}},c_{\beta_{2}},0,0,0,0,c_{\beta_{7}})$&\cr\hline&
100&&
$(c_{\beta_{1}},c_{\beta_{2}},0,c_{\beta_{4}},0,0,c_{\beta_{7}})$&&
101&&
$(c_{\beta_{1}},c_{\beta_{2}},0,0,0,c_{\beta_{6}},0)$&&
102&&
$(c_{\beta_{1}},c_{\beta_{2}},0,c_{\beta_{4}},0,c_{\beta_{6}},0)$&\cr\hline&
103&&
$(c_{\beta_{1}},c_{\beta_{2}},0,0,0,c_{\beta_{6}},c_{\beta_{7}})$&&
104&&
$(c_{\beta_{1}},c_{\beta_{2}},0,c_{\beta_{4}},0,c_{\beta_{6}},c_{\beta_{7}})$&&
105&&
$(c_{\beta_{1}},c_{\beta_{2}},0,0,c_{\beta_{5}},0,0)$&\cr\hline&
106&&
$(c_{\beta_{1}},c_{\beta_{2}},0,c_{\beta_{4}},c_{\beta_{5}},0,0)$&&
107&&
$(c_{\beta_{1}},c_{\beta_{2}},0,0,c_{\beta_{5}},0,c_{\beta_{7}})$&&
108&&
$(c_{\beta_{1}},c_{\beta_{2}},0,c_{\beta_{4}},c_{\beta_{5}},0,c_{\beta_{7}})$&\cr\hline&
109&&
$(c_{\beta_{1}},c_{\beta_{2}},0,0,c_{\beta_{5}},c_{\beta_{6}},0)$&&
110&&
$(c_{\beta_{1}},c_{\beta_{2}},0,c_{\beta_{4}},c_{\beta_{5}},c_{\beta_{6}},0)$&&
111&&
$(c_{\beta_{1}},c_{\beta_{2}},0,0,c_{\beta_{5}},c_{\beta_{6}},c_{\beta_{7}})$&\cr\hline&
112&&
$(c_{\beta_{1}},c_{\beta_{2}},0,c_{\beta_{4}},c_{\beta_{5}},c_{\beta_{6}},c_{\beta_{7}})$&&
&&
&&
&&
&\cr\hline}}

\begin{center}
	\tabl{Correspondence between $c$ and $D$ for $E_8$}
\end{center}

{\tabskip0pt\relax
\scriptspace=0pt\thinmuskip=1mu
\def\hline{\noalign{\hrule\allowbreak\vskip-.4pt\hrule}}
\def\Hline{\noalign{\hrule\kern2pt\hrule}}
\halign to \textwidth{\strut\raise1pt\copy\strutbox\lower1pt\copy\strutbox
\vrule\kern2pt\vrule#\unskip\tabskip0pt plus 1fil\relax
&\hfil\ignorespaces#\unskip&\vrule#\unskip
&\ignorespaces#\unskip\hfil&\vrule\kern2pt\vrule#\unskip
&\hfil\ignorespaces#\unskip&\vrule#\unskip
&\ignorespaces#\unskip\hfil&\vrule\kern2pt\vrule#\unskip
&\hfil\ignorespaces#\unskip&\vrule#\unskip
&\ignorespaces#\unskip\hfil&\vrule\kern2pt\vrule#\unskip\tabskip0pt\relax\cr\hline
&$D$\hfil\null&&\hfil Type of $c$&&$D$\hfil\null&&\hfil Type of $c$&&$D$\hfil\null&&\hfil Type of $c$&\cr\Hline
&1&&
$(0,0,0,0,0,0,0,0)$&&
2&&
$(0,0,0,0,c_{\beta_{5}},0,0,0)$&&
3&&
$(0,0,0,0,0,0,0,c_{\beta_{8}})$&\cr\hline&
4&&
$(0,0,0,0,c_{\beta_{5}},0,0,c_{\beta_{8}})$&&
5&&
$(0,0,0,0,0,0,c_{\beta_{7}},0)$&&
6&&
$(0,0,0,0,c_{\beta_{5}},0,c_{\beta_{7}},0)$&\cr\hline&
7&&
$(0,0,0,0,0,0,c_{\beta_{7}},c_{\beta_{8}})$&&
8&&
$(0,0,0,0,c_{\beta_{5}},0,c_{\beta_{7}},c_{\beta_{8}})$&&
9&&
$(0,0,0,0,0,c_{\beta_{6}},0,0)$&\cr\hline&
10&&
$(0,0,0,0,c_{\beta_{5}},c_{\beta_{6}},0,0)$&&
11&&
$(0,0,0,0,0,c_{\beta_{6}},0,c_{\beta_{8}})$&&
12&&
$(0,0,0,0,c_{\beta_{5}},c_{\beta_{6}},0,c_{\beta_{8}})$&\cr\hline&
13&&
$(0,0,0,0,0,c_{\beta_{6}},c_{\beta_{7}},0)$&&
14&&
$(0,0,0,0,c_{\beta_{5}},c_{\beta_{6}},c_{\beta_{7}},0)$&&
15&&
$(0,0,0,0,0,c_{\beta_{6}},c_{\beta_{7}},c_{\beta_{8}})$&\cr\hline&
16&&
$(0,0,0,0,c_{\beta_{5}},c_{\beta_{6}},c_{\beta_{7}},c_{\beta_{8}})$&&
17&&
$(0,0,0,c_{\beta_{4}},0,0,0,0)$&&
18&&
$(0,0,0,c_{\beta_{4}},c_{\beta_{5}},0,0,0)$&\cr\hline&
19&&
$(0,0,0,c_{\beta_{4}},0,0,0,c_{\beta_{8}})$&&
20&&
$(0,0,0,c_{\beta_{4}},c_{\beta_{5}},0,0,c_{\beta_{8}})$&&
21&&
$(0,0,0,c_{\beta_{4}},0,0,c_{\beta_{7}},0)$&\cr\hline&
22&&
$(0,0,0,c_{\beta_{4}},c_{\beta_{5}},0,c_{\beta_{7}},0)$&&
23&&
$(0,0,0,c_{\beta_{4}},0,0,c_{\beta_{7}},c_{\beta_{8}})$&&
24&&
$(0,0,0,c_{\beta_{4}},c_{\beta_{5}},0,c_{\beta_{7}},c_{\beta_{8}})$&\cr\hline&
25&&
$(0,0,0,c_{\beta_{4}},0,c_{\beta_{6}},0,0)$&&
26&&
$(0,0,0,c_{\beta_{4}},c_{\beta_{5}},c_{\beta_{6}},0,0)$&&
27&&
$(0,0,0,c_{\beta_{4}},0,c_{\beta_{6}},0,c_{\beta_{8}})$&\cr\hline&
28&&
$(0,0,0,c_{\beta_{4}},c_{\beta_{5}},c_{\beta_{6}},0,c_{\beta_{8}})$&&
29&&
$(0,0,0,c_{\beta_{4}},0,c_{\beta_{6}},c_{\beta_{7}},0)$&&
30&&
$(0,0,0,c_{\beta_{4}},c_{\beta_{5}},c_{\beta_{6}},c_{\beta_{7}},0)$&\cr\hline&
31&&
$(0,0,0,c_{\beta_{4}},0,c_{\beta_{6}},c_{\beta_{7}},c_{\beta_{8}})$&&
32&&
$(0,0,0,c_{\beta_{4}},c_{\beta_{5}},c_{\beta_{6}},c_{\beta_{7}},c_{\beta_{8}})$&&
33&&
$(0,0,c_{\beta_{3}},0,0,0,0,0)$&\cr\hline&
34&&
$(0,0,c_{\beta_{3}},0,c_{\beta_{5}},0,0,0)$&&
35&&
$(0,0,c_{\beta_{3}},0,0,0,0,c_{\beta_{8}})$&&
36&&
$(0,0,c_{\beta_{3}},0,c_{\beta_{5}},0,0,c_{\beta_{8}})$&\cr\hline&
37&&
$(0,0,c_{\beta_{3}},0,0,0,c_{\beta_{7}},0)$&&
38&&
$(0,0,c_{\beta_{3}},0,c_{\beta_{5}},0,c_{\beta_{7}},0)$&&
39&&
$(0,0,c_{\beta_{3}},0,0,0,c_{\beta_{7}},c_{\beta_{8}})$&\cr\hline&
40&&
$(0,0,c_{\beta_{3}},0,c_{\beta_{5}},0,c_{\beta_{7}},c_{\beta_{8}})$&&
41&&
$(0,0,c_{\beta_{3}},0,0,c_{\beta_{6}},0,0)$&&
42&&
$(0,0,c_{\beta_{3}},0,c_{\beta_{5}},c_{\beta_{6}},0,0)$&\cr\hline&
43&&
$(0,0,c_{\beta_{3}},0,0,c_{\beta_{6}},0,c_{\beta_{8}})$&&
44&&
$(0,0,c_{\beta_{3}},0,c_{\beta_{5}},c_{\beta_{6}},0,c_{\beta_{8}})$&&
45&&
$(0,0,c_{\beta_{3}},0,0,c_{\beta_{6}},c_{\beta_{7}},0)$&\cr\hline&
46&&
$(0,0,c_{\beta_{3}},0,c_{\beta_{5}},c_{\beta_{6}},c_{\beta_{7}},0)$&&
47&&
$(0,0,c_{\beta_{3}},0,0,c_{\beta_{6}},c_{\beta_{7}},c_{\beta_{8}})$&&
48&&
$(0,0,c_{\beta_{3}},0,c_{\beta_{5}},c_{\beta_{6}},c_{\beta_{7}},c_{\beta_{8}})$&\cr\hline&
49&&
$(0,0,c_{\beta_{3}},c_{\beta_{4}},0,0,0,0)$&&
50&&
$(0,0,c_{\beta_{3}},c_{\beta_{4}},c_{\beta_{5}},0,0,0)$&&
51&&
$(0,0,c_{\beta_{3}},c_{\beta_{4}},0,0,0,c_{\beta_{8}})$&\cr\hline&
52&&
$(0,0,c_{\beta_{3}},c_{\beta_{4}},c_{\beta_{5}},0,0,c_{\beta_{8}})$&&
53&&
$(0,0,c_{\beta_{3}},c_{\beta_{4}},0,0,c_{\beta_{7}},0)$&&
54&&
$(0,0,c_{\beta_{3}},c_{\beta_{4}},c_{\beta_{5}},0,c_{\beta_{7}},0)$&\cr\hline&
55&&
$(0,0,c_{\beta_{3}},c_{\beta_{4}},0,0,c_{\beta_{7}},c_{\beta_{8}})$&&
56&&
$(0,0,c_{\beta_{3}},c_{\beta_{4}},c_{\beta_{5}},0,c_{\beta_{7}},c_{\beta_{8}})$&&
57&&
$(0,0,c_{\beta_{3}},c_{\beta_{4}},0,c_{\beta_{6}},0,0)$&\cr\hline&
58&&
$(0,0,c_{\beta_{3}},c_{\beta_{4}},c_{\beta_{5}},c_{\beta_{6}},0,0)$&&
59&&
$(0,0,c_{\beta_{3}},c_{\beta_{4}},0,c_{\beta_{6}},0,c_{\beta_{8}})$&&
60&&
$(0,0,c_{\beta_{3}},c_{\beta_{4}},c_{\beta_{5}},c_{\beta_{6}},0,c_{\beta_{8}})$&\cr\hline&
61&&
$(0,0,c_{\beta_{3}},c_{\beta_{4}},0,c_{\beta_{6}},c_{\beta_{7}},0)$&&
62&&
$(0,0,c_{\beta_{3}},c_{\beta_{4}},c_{\beta_{5}},c_{\beta_{6}},c_{\beta_{7}},0)$&&
63&&
$(0,0,c_{\beta_{3}},c_{\beta_{4}},0,c_{\beta_{6}},c_{\beta_{7}},c_{\beta_{8}})$&\cr\hline&
64&&
$(0,0,c_{\beta_{3}},c_{\beta_{4}},c_{\beta_{5}},c_{\beta_{6}},c_{\beta_{7}},c_{\beta_{8}})$&&
65&&
$(0,c_{\beta_{2}},0,0,0,0,0,0)$&&
66&&
$(0,c_{\beta_{2}},0,0,c_{\beta_{5}},0,0,0)$&\cr\hline&
67&&
$(0,c_{\beta_{2}},0,0,0,0,0,c_{\beta_{8}})$&&
68&&
$(0,c_{\beta_{2}},0,0,c_{\beta_{5}},0,0,c_{\beta_{8}})$&&
69&&
$(0,c_{\beta_{2}},0,0,0,0,c_{\beta_{7}},0)$&\cr\hline&
70&&
$(0,c_{\beta_{2}},0,0,c_{\beta_{5}},0,c_{\beta_{7}},0)$&&
71&&
$(0,c_{\beta_{2}},0,0,0,0,c_{\beta_{7}},c_{\beta_{8}})$&&
72&&
$(0,c_{\beta_{2}},0,0,c_{\beta_{5}},0,c_{\beta_{7}},c_{\beta_{8}})$&\cr\hline&
73&&
$(0,c_{\beta_{2}},0,0,0,c_{\beta_{6}},0,0)$&&
74&&
$(0,c_{\beta_{2}},0,0,c_{\beta_{5}},c_{\beta_{6}},0,0)$&&
75&&
$(0,c_{\beta_{2}},0,0,0,c_{\beta_{6}},0,c_{\beta_{8}})$&\cr\hline&
76&&
$(0,c_{\beta_{2}},0,0,c_{\beta_{5}},c_{\beta_{6}},0,c_{\beta_{8}})$&&
77&&
$(0,c_{\beta_{2}},0,0,0,c_{\beta_{6}},c_{\beta_{7}},0)$&&
78&&
$(0,c_{\beta_{2}},0,0,c_{\beta_{5}},c_{\beta_{6}},c_{\beta_{7}},0)$&\cr\hline&
79&&
$(0,c_{\beta_{2}},0,0,0,c_{\beta_{6}},c_{\beta_{7}},c_{\beta_{8}})$&&
80&&
$(0,c_{\beta_{2}},0,0,c_{\beta_{5}},c_{\beta_{6}},c_{\beta_{7}},c_{\beta_{8}})$&&
81&&
$(0,c_{\beta_{2}},0,c_{\beta_{4}},0,0,0,0)$&\cr\hline&
82&&
$(0,c_{\beta_{2}},0,c_{\beta_{4}},c_{\beta_{5}},0,0,0)$&&
83&&
$(0,c_{\beta_{2}},0,c_{\beta_{4}},0,0,0,c_{\beta_{8}})$&&
84&&
$(0,c_{\beta_{2}},0,c_{\beta_{4}},c_{\beta_{5}},0,0,c_{\beta_{8}})$&\cr\hline&
85&&
$(0,c_{\beta_{2}},0,c_{\beta_{4}},0,0,c_{\beta_{7}},0)$&&
86&&
$(0,c_{\beta_{2}},0,c_{\beta_{4}},c_{\beta_{5}},0,c_{\beta_{7}},0)$&&
87&&
$(0,c_{\beta_{2}},0,c_{\beta_{4}},0,0,c_{\beta_{7}},c_{\beta_{8}})$&\cr\hline&
88&&
$(0,c_{\beta_{2}},0,c_{\beta_{4}},c_{\beta_{5}},0,c_{\beta_{7}},c_{\beta_{8}})$&&
89&&
$(0,c_{\beta_{2}},0,c_{\beta_{4}},0,c_{\beta_{6}},0,0)$&&
90&&
$(0,c_{\beta_{2}},0,c_{\beta_{4}},c_{\beta_{5}},c_{\beta_{6}},0,0)$&\cr\hline&
91&&
$(0,c_{\beta_{2}},0,c_{\beta_{4}},0,c_{\beta_{6}},0,c_{\beta_{8}})$&&
92&&
$(0,c_{\beta_{2}},0,c_{\beta_{4}},c_{\beta_{5}},c_{\beta_{6}},0,c_{\beta_{8}})$&&
93&&
$(0,c_{\beta_{2}},0,c_{\beta_{4}},0,c_{\beta_{6}},c_{\beta_{7}},0)$&\cr\hline&
94&&
$(0,c_{\beta_{2}},0,c_{\beta_{4}},c_{\beta_{5}},c_{\beta_{6}},c_{\beta_{7}},0)$&&
95&&
$(0,c_{\beta_{2}},0,c_{\beta_{4}},0,c_{\beta_{6}},c_{\beta_{7}},c_{\beta_{8}})$&&
96&&
$(0,c_{\beta_{2}},0,c_{\beta_{4}},c_{\beta_{5}},c_{\beta_{6}},c_{\beta_{7}},c_{\beta_{8}})$&\cr\hline&
97&&
$(0,c_{\beta_{2}},c_{\beta_{3}},0,0,0,0,0)$&&
98&&
$(0,c_{\beta_{2}},c_{\beta_{3}},0,c_{\beta_{5}},0,0,0)$&&
99&&
$(0,c_{\beta_{2}},c_{\beta_{3}},0,0,0,0,c_{\beta_{8}})$&\cr\hline&
100&&
$(0,c_{\beta_{2}},c_{\beta_{3}},0,c_{\beta_{5}},0,0,c_{\beta_{8}})$&&
101&&
$(0,c_{\beta_{2}},c_{\beta_{3}},0,0,0,c_{\beta_{7}},0)$&&
102&&
$(0,c_{\beta_{2}},c_{\beta_{3}},0,c_{\beta_{5}},0,c_{\beta_{7}},0)$&\cr\hline&
103&&
$(0,c_{\beta_{2}},c_{\beta_{3}},0,0,0,c_{\beta_{7}},c_{\beta_{8}})$&&
104&&
$(0,c_{\beta_{2}},c_{\beta_{3}},0,c_{\beta_{5}},0,c_{\beta_{7}},c_{\beta_{8}})$&&
105&&
$(0,c_{\beta_{2}},c_{\beta_{3}},0,0,c_{\beta_{6}},0,0)$&\cr\hline&
106&&
$(0,c_{\beta_{2}},c_{\beta_{3}},0,c_{\beta_{5}},c_{\beta_{6}},0,0)$&&
107&&
$(0,c_{\beta_{2}},c_{\beta_{3}},0,0,c_{\beta_{6}},0,c_{\beta_{8}})$&&
108&&
$(0,c_{\beta_{2}},c_{\beta_{3}},0,c_{\beta_{5}},c_{\beta_{6}},0,c_{\beta_{8}})$&\cr\hline&
109&&
$(0,c_{\beta_{2}},c_{\beta_{3}},0,0,c_{\beta_{6}},c_{\beta_{7}},0)$&&
110&&
$(0,c_{\beta_{2}},c_{\beta_{3}},0,c_{\beta_{5}},c_{\beta_{6}},c_{\beta_{7}},0)$&&
111&&
$(0,c_{\beta_{2}},c_{\beta_{3}},0,0,c_{\beta_{6}},c_{\beta_{7}},c_{\beta_{8}})$&\cr\hline&
112&&
$(0,c_{\beta_{2}},c_{\beta_{3}},0,c_{\beta_{5}},c_{\beta_{6}},c_{\beta_{7}},c_{\beta_{8}})$&&
113&&
$(0,c_{\beta_{2}},c_{\beta_{3}},c_{\beta_{4}},0,0,0,0)$&&
114&&
$(0,c_{\beta_{2}},c_{\beta_{3}},c_{\beta_{4}},c_{\beta_{5}},0,0,0)$&\cr\hline&
115&&
$(0,c_{\beta_{2}},c_{\beta_{3}},c_{\beta_{4}},0,0,0,c_{\beta_{8}})$&&
116&&
$(0,c_{\beta_{2}},c_{\beta_{3}},c_{\beta_{4}},c_{\beta_{5}},0,0,c_{\beta_{8}})$&&
117&&
$(0,c_{\beta_{2}},c_{\beta_{3}},c_{\beta_{4}},0,0,c_{\beta_{7}},0)$&\cr\hline&
118&&
$(0,c_{\beta_{2}},c_{\beta_{3}},c_{\beta_{4}},c_{\beta_{5}},0,c_{\beta_{7}},0)$&&
119&&
$(0,c_{\beta_{2}},c_{\beta_{3}},c_{\beta_{4}},0,0,c_{\beta_{7}},c_{\beta_{8}})$&&
120&&
$(0,c_{\beta_{2}},c_{\beta_{3}},c_{\beta_{4}},c_{\beta_{5}},0,c_{\beta_{7}},c_{\beta_{8}})$&\cr\hline&
121&&
$(0,c_{\beta_{2}},c_{\beta_{3}},c_{\beta_{4}},0,c_{\beta_{6}},0,0)$&&
122&&
$(0,c_{\beta_{2}},c_{\beta_{3}},c_{\beta_{4}},c_{\beta_{5}},c_{\beta_{6}},0,0)$&&
123&&
$(0,c_{\beta_{2}},c_{\beta_{3}},c_{\beta_{4}},0,c_{\beta_{6}},0,c_{\beta_{8}})$&\cr\hline&
124&&
$(0,c_{\beta_{2}},c_{\beta_{3}},c_{\beta_{4}},c_{\beta_{5}},c_{\beta_{6}},0,c_{\beta_{8}})$&&
125&&
$(0,c_{\beta_{2}},c_{\beta_{3}},c_{\beta_{4}},0,c_{\beta_{6}},c_{\beta_{7}},0)$&&
126&&
$(0,c_{\beta_{2}},c_{\beta_{3}},c_{\beta_{4}},c_{\beta_{5}},c_{\beta_{6}},c_{\beta_{7}},0)$&\cr\hline&
127&&
$(0,c_{\beta_{2}},c_{\beta_{3}},c_{\beta_{4}},0,c_{\beta_{6}},c_{\beta_{7}},c_{\beta_{8}})$&&
128&&
$(0,c_{\beta_{2}},c_{\beta_{3}},c_{\beta_{4}},c_{\beta_{5}},c_{\beta_{6}},c_{\beta_{7}},c_{\beta_{8}})$&&
129&&
$(c_{\beta_{1}},0,0,0,0,0,0,0)$&\cr\hline&
130&&
$(c_{\beta_{1}},0,0,0,c_{\beta_{5}},0,0,0)$&&
131&&
$(c_{\beta_{1}},0,0,0,0,0,0,c_{\beta_{8}})$&&
132&&
$(c_{\beta_{1}},0,0,0,c_{\beta_{5}},0,0,c_{\beta_{8}})$&\cr\hline&
133&&
$(c_{\beta_{1}},0,0,0,0,0,c_{\beta_{7}},0)$&&
134&&
$(c_{\beta_{1}},0,0,0,c_{\beta_{5}},0,c_{\beta_{7}},0)$&&
135&&
$(c_{\beta_{1}},0,0,0,0,0,c_{\beta_{7}},c_{\beta_{8}})$&\cr\hline&
136&&
$(c_{\beta_{1}},0,0,0,c_{\beta_{5}},0,c_{\beta_{7}},c_{\beta_{8}})$&&
137&&
$(c_{\beta_{1}},0,0,0,0,c_{\beta_{6}},0,0)$&&
138&&
$(c_{\beta_{1}},0,0,0,c_{\beta_{5}},c_{\beta_{6}},0,0)$&\cr\hline&
139&&
$(c_{\beta_{1}},0,0,0,0,c_{\beta_{6}},0,c_{\beta_{8}})$&&
140&&
$(c_{\beta_{1}},0,0,0,c_{\beta_{5}},c_{\beta_{6}},0,c_{\beta_{8}})$&&
141&&
$(c_{\beta_{1}},0,0,0,0,c_{\beta_{6}},c_{\beta_{7}},0)$&\cr\hline&
142&&
$(c_{\beta_{1}},0,0,0,c_{\beta_{5}},c_{\beta_{6}},c_{\beta_{7}},0)$&&
143&&
$(c_{\beta_{1}},0,0,0,0,c_{\beta_{6}},c_{\beta_{7}},c_{\beta_{8}})$&&
144&&
$(c_{\beta_{1}},0,0,0,c_{\beta_{5}},c_{\beta_{6}},c_{\beta_{7}},c_{\beta_{8}})$&\cr\hline&
145&&
$(c_{\beta_{1}},0,0,c_{\beta_{4}},0,0,0,0)$&&
146&&
$(c_{\beta_{1}},0,0,c_{\beta_{4}},c_{\beta_{5}},0,0,0)$&&
147&&
$(c_{\beta_{1}},0,0,c_{\beta_{4}},0,0,0,c_{\beta_{8}})$&\cr\hline&
148&&
$(c_{\beta_{1}},0,0,c_{\beta_{4}},c_{\beta_{5}},0,0,c_{\beta_{8}})$&&
149&&
$(c_{\beta_{1}},0,0,c_{\beta_{4}},0,0,c_{\beta_{7}},0)$&&
150&&
$(c_{\beta_{1}},0,0,c_{\beta_{4}},c_{\beta_{5}},0,c_{\beta_{7}},0)$&\cr\hline&
151&&
$(c_{\beta_{1}},0,0,c_{\beta_{4}},0,0,c_{\beta_{7}},c_{\beta_{8}})$&&
152&&
$(c_{\beta_{1}},0,0,c_{\beta_{4}},c_{\beta_{5}},0,c_{\beta_{7}},c_{\beta_{8}})$&&
153&&
$(c_{\beta_{1}},0,0,c_{\beta_{4}},0,c_{\beta_{6}},0,0)$&\cr\hline&
154&&
$(c_{\beta_{1}},0,0,c_{\beta_{4}},c_{\beta_{5}},c_{\beta_{6}},0,0)$&&
155&&
$(c_{\beta_{1}},0,0,c_{\beta_{4}},0,c_{\beta_{6}},0,c_{\beta_{8}})$&&
156&&
$(c_{\beta_{1}},0,0,c_{\beta_{4}},c_{\beta_{5}},c_{\beta_{6}},0,c_{\beta_{8}})$&\cr\hline&
157&&
$(c_{\beta_{1}},0,0,c_{\beta_{4}},0,c_{\beta_{6}},c_{\beta_{7}},0)$&&
158&&
$(c_{\beta_{1}},0,0,c_{\beta_{4}},c_{\beta_{5}},c_{\beta_{6}},c_{\beta_{7}},0)$&&
159&&
$(c_{\beta_{1}},0,0,c_{\beta_{4}},0,c_{\beta_{6}},c_{\beta_{7}},c_{\beta_{8}})$&\cr\hline&
160&&
$(c_{\beta_{1}},0,0,c_{\beta_{4}},c_{\beta_{5}},c_{\beta_{6}},c_{\beta_{7}},c_{\beta_{8}})$&&
161&&
$(c_{\beta_{1}},0,c_{\beta_{3}},0,0,0,0,0)$&&
162&&
$(c_{\beta_{1}},0,c_{\beta_{3}},0,c_{\beta_{5}},0,0,0)$&\cr\hline&
163&&
$(c_{\beta_{1}},0,c_{\beta_{3}},0,0,0,0,c_{\beta_{8}})$&&
164&&
$(c_{\beta_{1}},0,c_{\beta_{3}},0,c_{\beta_{5}},0,0,c_{\beta_{8}})$&&
165&&
$(c_{\beta_{1}},0,c_{\beta_{3}},0,0,0,c_{\beta_{7}},0)$&\cr\hline&
166&&
$(c_{\beta_{1}},0,c_{\beta_{3}},0,c_{\beta_{5}},0,c_{\beta_{7}},0)$&&
167&&
$(c_{\beta_{1}},0,c_{\beta_{3}},0,0,0,c_{\beta_{7}},c_{\beta_{8}})$&&
168&&
$(c_{\beta_{1}},0,c_{\beta_{3}},0,c_{\beta_{5}},0,c_{\beta_{7}},c_{\beta_{8}})$&\cr\hline&
169&&
$(c_{\beta_{1}},0,c_{\beta_{3}},0,0,c_{\beta_{6}},0,0)$&&
170&&
$(c_{\beta_{1}},0,c_{\beta_{3}},0,c_{\beta_{5}},c_{\beta_{6}},0,0)$&&
171&&
$(c_{\beta_{1}},0,c_{\beta_{3}},0,0,c_{\beta_{6}},0,c_{\beta_{8}})$&\cr\hline&
172&&
$(c_{\beta_{1}},0,c_{\beta_{3}},0,c_{\beta_{5}},c_{\beta_{6}},0,c_{\beta_{8}})$&&
173&&
$(c_{\beta_{1}},0,c_{\beta_{3}},0,0,c_{\beta_{6}},c_{\beta_{7}},0)$&&
174&&
$(c_{\beta_{1}},0,c_{\beta_{3}},0,c_{\beta_{5}},c_{\beta_{6}},c_{\beta_{7}},0)$&\cr\hline&
175&&
$(c_{\beta_{1}},0,c_{\beta_{3}},0,0,c_{\beta_{6}},c_{\beta_{7}},c_{\beta_{8}})$&&
176&&
$(c_{\beta_{1}},0,c_{\beta_{3}},0,c_{\beta_{5}},c_{\beta_{6}},c_{\beta_{7}},c_{\beta_{8}})$&&
177&&
$(c_{\beta_{1}},0,c_{\beta_{3}},c_{\beta_{4}},0,0,0,0)$&\cr\hline&
178&&
$(c_{\beta_{1}},0,c_{\beta_{3}},c_{\beta_{4}},c_{\beta_{5}},0,0,0)$&&
179&&
$(c_{\beta_{1}},0,c_{\beta_{3}},c_{\beta_{4}},0,0,0,c_{\beta_{8}})$&&
180&&
$(c_{\beta_{1}},0,c_{\beta_{3}},c_{\beta_{4}},c_{\beta_{5}},0,0,c_{\beta_{8}})$&\cr\hline&
181&&
$(c_{\beta_{1}},0,c_{\beta_{3}},c_{\beta_{4}},0,0,c_{\beta_{7}},0)$&&
182&&
$(c_{\beta_{1}},0,c_{\beta_{3}},c_{\beta_{4}},c_{\beta_{5}},0,c_{\beta_{7}},0)$&&
183&&
$(c_{\beta_{1}},0,c_{\beta_{3}},c_{\beta_{4}},0,0,c_{\beta_{7}},c_{\beta_{8}})$&\cr\hline&
184&&
$(c_{\beta_{1}},0,c_{\beta_{3}},c_{\beta_{4}},c_{\beta_{5}},0,c_{\beta_{7}},c_{\beta_{8}})$&&
185&&
$(c_{\beta_{1}},0,c_{\beta_{3}},c_{\beta_{4}},0,c_{\beta_{6}},0,0)$&&
186&&
$(c_{\beta_{1}},0,c_{\beta_{3}},c_{\beta_{4}},c_{\beta_{5}},c_{\beta_{6}},0,0)$&\cr\hline&
187&&
$(c_{\beta_{1}},0,c_{\beta_{3}},c_{\beta_{4}},0,c_{\beta_{6}},0,c_{\beta_{8}})$&&
188&&
$(c_{\beta_{1}},0,c_{\beta_{3}},c_{\beta_{4}},c_{\beta_{5}},c_{\beta_{6}},0,c_{\beta_{8}})$&&
189&&
$(c_{\beta_{1}},0,c_{\beta_{3}},c_{\beta_{4}},0,c_{\beta_{6}},c_{\beta_{7}},0)$&\cr\hline&
190&&
$(c_{\beta_{1}},0,c_{\beta_{3}},c_{\beta_{4}},c_{\beta_{5}},c_{\beta_{6}},c_{\beta_{7}},0)$&&
191&&
$(c_{\beta_{1}},0,c_{\beta_{3}},c_{\beta_{4}},0,c_{\beta_{6}},c_{\beta_{7}},c_{\beta_{8}})$&&
192&&
$(c_{\beta_{1}},0,c_{\beta_{3}},c_{\beta_{4}},c_{\beta_{5}},c_{\beta_{6}},c_{\beta_{7}},c_{\beta_{8}})$&\cr\hline&
193&&
$(c_{\beta_{1}},c_{\beta_{2}},0,0,0,0,0,0)$&&
194&&
$(c_{\beta_{1}},c_{\beta_{2}},0,0,c_{\beta_{5}},0,0,0)$&&
195&&
$(c_{\beta_{1}},c_{\beta_{2}},0,0,0,0,0,c_{\beta_{8}})$&\cr\hline&
196&&
$(c_{\beta_{1}},c_{\beta_{2}},0,0,c_{\beta_{5}},0,0,c_{\beta_{8}})$&&
197&&
$(c_{\beta_{1}},c_{\beta_{2}},0,0,0,0,c_{\beta_{7}},0)$&&
198&&
$(c_{\beta_{1}},c_{\beta_{2}},0,0,c_{\beta_{5}},0,c_{\beta_{7}},0)$&\cr\hline&
199&&
$(c_{\beta_{1}},c_{\beta_{2}},0,0,0,0,c_{\beta_{7}},c_{\beta_{8}})$&&
200&&
$(c_{\beta_{1}},c_{\beta_{2}},0,0,c_{\beta_{5}},0,c_{\beta_{7}},c_{\beta_{8}})$&&
201&&
$(c_{\beta_{1}},c_{\beta_{2}},0,0,0,c_{\beta_{6}},0,0)$&\cr\hline&
202&&
$(c_{\beta_{1}},c_{\beta_{2}},0,0,c_{\beta_{5}},c_{\beta_{6}},0,0)$&&
203&&
$(c_{\beta_{1}},c_{\beta_{2}},0,0,0,c_{\beta_{6}},0,c_{\beta_{8}})$&&
204&&
$(c_{\beta_{1}},c_{\beta_{2}},0,0,c_{\beta_{5}},c_{\beta_{6}},0,c_{\beta_{8}})$&\cr\hline&
205&&
$(c_{\beta_{1}},c_{\beta_{2}},0,0,0,c_{\beta_{6}},c_{\beta_{7}},0)$&&
206&&
$(c_{\beta_{1}},c_{\beta_{2}},0,0,c_{\beta_{5}},c_{\beta_{6}},c_{\beta_{7}},0)$&&
207&&
$(c_{\beta_{1}},c_{\beta_{2}},0,0,0,c_{\beta_{6}},c_{\beta_{7}},c_{\beta_{8}})$&\cr\hline&
208&&
$(c_{\beta_{1}},c_{\beta_{2}},0,0,c_{\beta_{5}},c_{\beta_{6}},c_{\beta_{7}},c_{\beta_{8}})$&&
209&&
$(c_{\beta_{1}},c_{\beta_{2}},0,c_{\beta_{4}},0,0,0,0)$&&
210&&
$(c_{\beta_{1}},c_{\beta_{2}},0,c_{\beta_{4}},c_{\beta_{5}},0,0,0)$&\cr\hline&
211&&
$(c_{\beta_{1}},c_{\beta_{2}},0,c_{\beta_{4}},0,0,0,c_{\beta_{8}})$&&
212&&
$(c_{\beta_{1}},c_{\beta_{2}},0,c_{\beta_{4}},c_{\beta_{5}},0,0,c_{\beta_{8}})$&&
213&&
$(c_{\beta_{1}},c_{\beta_{2}},0,c_{\beta_{4}},0,0,c_{\beta_{7}},0)$&\cr\hline&
214&&
$(c_{\beta_{1}},c_{\beta_{2}},0,c_{\beta_{4}},c_{\beta_{5}},0,c_{\beta_{7}},0)$&&
215&&
$(c_{\beta_{1}},c_{\beta_{2}},0,c_{\beta_{4}},0,0,c_{\beta_{7}},c_{\beta_{8}})$&&
216&&
$(c_{\beta_{1}},c_{\beta_{2}},0,c_{\beta_{4}},c_{\beta_{5}},0,c_{\beta_{7}},c_{\beta_{8}})$&\cr\hline&
217&&
$(c_{\beta_{1}},c_{\beta_{2}},0,c_{\beta_{4}},0,c_{\beta_{6}},0,0)$&&
218&&
$(c_{\beta_{1}},c_{\beta_{2}},0,c_{\beta_{4}},c_{\beta_{5}},c_{\beta_{6}},0,0)$&&
219&&
$(c_{\beta_{1}},c_{\beta_{2}},0,c_{\beta_{4}},0,c_{\beta_{6}},0,c_{\beta_{8}})$&\cr\hline&
220&&
$(c_{\beta_{1}},c_{\beta_{2}},0,c_{\beta_{4}},c_{\beta_{5}},c_{\beta_{6}},0,c_{\beta_{8}})$&&
221&&
$(c_{\beta_{1}},c_{\beta_{2}},0,c_{\beta_{4}},0,c_{\beta_{6}},c_{\beta_{7}},0)$&&
222&&
$(c_{\beta_{1}},c_{\beta_{2}},0,c_{\beta_{4}},c_{\beta_{5}},c_{\beta_{6}},c_{\beta_{7}},0)$&\cr\hline&
223&&
$(c_{\beta_{1}},c_{\beta_{2}},0,c_{\beta_{4}},0,c_{\beta_{6}},c_{\beta_{7}},c_{\beta_{8}})$&&
224&&
$(c_{\beta_{1}},c_{\beta_{2}},0,c_{\beta_{4}},c_{\beta_{5}},c_{\beta_{6}},c_{\beta_{7}},c_{\beta_{8}})$&&
225&&
$(c_{\beta_{1}},c_{\beta_{2}},c_{\beta_{3}},0,0,0,0,0)$&\cr\hline&
226&&
$(c_{\beta_{1}},c_{\beta_{2}},c_{\beta_{3}},0,c_{\beta_{5}},0,0,0)$&&
227&&
$(c_{\beta_{1}},c_{\beta_{2}},c_{\beta_{3}},0,0,0,0,c_{\beta_{8}})$&&
228&&
$(c_{\beta_{1}},c_{\beta_{2}},c_{\beta_{3}},0,c_{\beta_{5}},0,0,c_{\beta_{8}})$&\cr\hline&
229&&
$(c_{\beta_{1}},c_{\beta_{2}},c_{\beta_{3}},0,0,0,c_{\beta_{7}},0)$&&
230&&
$(c_{\beta_{1}},c_{\beta_{2}},c_{\beta_{3}},0,c_{\beta_{5}},0,c_{\beta_{7}},0)$&&
231&&
$(c_{\beta_{1}},c_{\beta_{2}},c_{\beta_{3}},0,0,0,c_{\beta_{7}},c_{\beta_{8}})$&\cr\hline&
232&&
$(c_{\beta_{1}},c_{\beta_{2}},c_{\beta_{3}},0,c_{\beta_{5}},0,c_{\beta_{7}},c_{\beta_{8}})$&&
233&&
$(c_{\beta_{1}},c_{\beta_{2}},c_{\beta_{3}},0,0,c_{\beta_{6}},0,0)$&&
234&&
$(c_{\beta_{1}},c_{\beta_{2}},c_{\beta_{3}},0,c_{\beta_{5}},c_{\beta_{6}},0,0)$&\cr\hline&
235&&
$(c_{\beta_{1}},c_{\beta_{2}},c_{\beta_{3}},0,0,c_{\beta_{6}},0,c_{\beta_{8}})$&&
236&&
$(c_{\beta_{1}},c_{\beta_{2}},c_{\beta_{3}},0,c_{\beta_{5}},c_{\beta_{6}},0,c_{\beta_{8}})$&&
237&&
$(c_{\beta_{1}},c_{\beta_{2}},c_{\beta_{3}},0,0,c_{\beta_{6}},c_{\beta_{7}},0)$&\cr\hline&
238&&
$(c_{\beta_{1}},c_{\beta_{2}},c_{\beta_{3}},0,c_{\beta_{5}},c_{\beta_{6}},c_{\beta_{7}},0)$&&
239&&
$(c_{\beta_{1}},c_{\beta_{2}},c_{\beta_{3}},0,0,c_{\beta_{6}},c_{\beta_{7}},c_{\beta_{8}})$&&
240&&
$(c_{\beta_{1}},c_{\beta_{2}},c_{\beta_{3}},0,c_{\beta_{5}},c_{\beta_{6}},c_{\beta_{7}},c_{\beta_{8}})$&\cr\hline}}

\newpage
\begin{center}
\textbf{Appendix B}
\end{center}

Below we present the listing of the code generating tables from Appendix A.

Code for $E_7$:

\begin{verbatim}
import time
import copy
print(time.ctime())
positive_roots=['0100000', '0101000', '0101100', '0101110', '0101111',
                '0111000', '0111100', '0111110', '0111111', '0112100',
                '0112110', '0112111', '0112210', '0112211', '0112221',
                '0010000', '0011000', '0011100', '0011110', '0011111',
                '0001000', '0001100', '0001110', '0001111', '0000100',
                '0000110', '0000111', '0000010', '0000011', '0000001',
                '2234321', '1223210', '1223211', '1223221', '1111000',
                '1223321', '1111100', '1111110', '1111111', '1224321',
                '1112100', '1112110', '1112111', '1112210', '1112211',
                '1112221', '1000000', '1234321', '1122100', '1122110',
                '1122111', '1122210', '1122211', '1122221', '1010000',
                '1123210', '1123211', '1123221', '1011000', '1123321',
                '1011100', '1011110', '1011111']
dim=7
scalar_product_matrix=[[2, 0, -1, 0, 0, 0, 0], [0, 2, 0, -1, 0, 0, 0],
                       [-1, 0, 2, -1, 0, 0, 0], [0, -1, -1, 2, -1, 0, 0],
                       [0, 0, 0, -1, 2, -1, 0], [0, 0, 0, 0, -1, 2, -1],
                       [0, 0, 0, 0, 0, -1, 2]]
def dot_product(alpha,beta):
    sum=0
    for i in range(0,len(alpha)):
        for j in range(0,len(beta)):
            sum=sum+int(alpha[i])*int(beta[j])*scalar_product_matrix[i][j]
    return(sum)
def dot_product_matrix():
    list_dot_products=[]
    n=len(positive_roots)
    for i in range(0,n):
        list_dot_products.append([])
        for j in range(0,n):
            list_dot_products[i].append(dot_product(positive_roots[i],
                                                    positive_roots[j]))
    return(list_dot_products)
dpm=dot_product_matrix()
print('matrix calculated')
print(time.ctime())
def test_orth(list_of_roots):
    n=len(list_of_roots)
    b=1
    for i in range(0,n-1):
        for j in range(i+1,n):
            if (dpm[list_of_roots[i]][list_of_roots[j]])!=0.0:
                b=0
    return(b)
def root_to_list(alpha):
    result_list=[]
    for i in range(0,len(alpha)):
        result_list.append(int(alpha[i]))
    return(result_list)
def positive_roots_list():
    result_list=[]
    for i in range(0,len(positive_roots)):
        result_list.append(root_to_list(positive_roots[i]))
    return(result_list)
prl=positive_roots_list()
def roots_to_list(list_of_roots):
    result_list=[]
    for i in range(0,len(list_of_roots)):
        result_list.append(prl[list_of_roots[i]])
    return(result_list)
def all_orth_subset():
    result_list=[]
    result_list_numbers=[[]]
    n=len(positive_roots)
    for i in range(0,n):
        result_list.append([root_to_list(positive_roots[i])])
        result_list_numbers[0].append([i])
    for i in range(1,dim):
        result_list_numbers.append([])
        for j in range(0,len(result_list_numbers[i-1])):
            test_list=copy.copy(result_list_numbers[i-1][j])
            test_list.append(result_list_numbers[i-1][j][i-1])
            while test_list[i]+1<n:
                test_list[i]=test_list[i]+1
                if test_orth(test_list)==1:
                    result_list_numbers[i].append(copy.copy(test_list))
                    result_list.append(roots_to_list(test_list))
    return(result_list)
aos=all_orth_subset()
print(len(aos))
print(time.ctime())
def sum_list(list1,list2):
    list3=[]
    l=len(list1)
    for i in range(0,l):
        list3.append(list1[i]+list2[i])
    return(list3)
def list_sum(list):
    s=0
    n=len(list)
    for i in range(0,n):
        s=s+list[i]
    return(s)
def list_mult_const(list1,const):
    list2=[]
    l=len(list1)
    for i in range(0,l):
        list2.append(const*list1[i])
    return(list2)
def list_weight_sum(list1,coef_list):
    result_list=list_mult_const(list1[0],coef_list[0])
    for i in range(1,len(list1)):
        result_list=sum_list(result_list,list_mult_const(list1[i],coef_list[i]))
    return(result_list)
cascade=[[2, 2, 3, 4, 3, 2, 1], [0, 1, 1, 2, 2, 2, 1], [0, 1, 1, 2, 1, 0, 0],
         [0, 1, 0, 0, 0, 0, 0], [0, 0, 1, 0, 0, 0, 0], [0, 0, 0, 0, 1, 0, 0],
         [0, 0, 0, 0, 0, 0, 1]]
weights=[]
list_cascade_weights=[[1, 0, 0, 0, 0, 0, 0],[1, 1, 0, 0, 0, 0, 0],
                      [2, 1, 1, 0, 0, 0, 0],[2, 1, 1, 1, 0, 0, 0],
                      [3, 1, 1, 0, 1, 0, 0],[3, 2, 1, 0, 0, 1, 0],
                      [1, 1, 0, 0, 0, 0, 1]]
for i in range(0,len(list_cascade_weights)):
    weights.append(list_weight_sum(cascade,list_cascade_weights[i]))
degree=[]
for i in range(0,len(weights)):
    degree.append(list_sum(list_cascade_weights[i]))
print(degree)
def test_orth_subset(orth_subset,weights,degree,test_cons):
    n=len(orth_subset)
    result_list=[]
    for j in range(0,len(weights)):
        result_list.append(0)
    for j in range(0,len(weights)):
        list1=[]
        for i in range(0,n):
            list1.append(0)
        b=0
        while (b==0):
            list1[n-1]=list1[n-1]+1
            for i in range(n-1,0,-1):
                if list1[i]>degree[j]:
                    list1[i]=0
                    list1[i-1]=list1[i-1]+1
            if list_sum(list1)==degree[j]:
                list2=list_mult_const(orth_subset[0],list1[0])
                for i in range(1,n):
                    list2=sum_list(list2,list_mult_const(orth_subset[i],list1[i]))
                if list2==weights[j]:
                    result_list[j]=result_list[j]+1
            if list1[0]==degree[j]:
                b=1
    if result_list==test_cons:
        res=1
    else:
        res=0
    return(res)
def calculate_test_cons_for_orth_subset(orth_subset,weights,degree):
    n=len(orth_subset)
    result_list=[]
    for j in range(0,len(weights)):
        result_list.append(0)
    for j in range(0,len(weights)):
        list1=[]
        for i in range(0,n):
            list1.append(0)
        b=0
        while (b==0):
            list1[n-1]=list1[n-1]+1
            for i in range(n-1,0,-1):
                if list1[i]>degree[j]:
                    list1[i]=0
                    list1[i-1]=list1[i-1]+1
            if list_sum(list1)==degree[j]:
                list2=list_mult_const(orth_subset[0],list1[0])
                for i in range(1,n):
                    list2=sum_list(list2,list_mult_const(orth_subset[i],list1[i]))
                if list2==weights[j]:
                    result_list[j]=result_list[j]+1
            if list1[0]==degree[j]:
                b=1
    return(result_list)
def border_tc(aos):
    result_list=[[0,0]]
    length=0
    for i in range(0,len(aos)):
        length1=len(aos[i])-1
        if length1>length:
            result_list[length][1]=i
            result_list.append([i,0])
            length=length1
    result_list[length][1]=len(aos)
    return(result_list)
print(border_tc(aos))
btc=border_tc(aos)
def calculate_all_test_cons_for_all_orth_subset(aos,weights,degree):
    atc=[]
    percent=0
    for i in range(0,len(aos)):
        percent1=int(100*(i+1)/len(aos))
        if percent1>percent:
            percent=percent1
            print(str(percent)+'%')
            print(time.ctime())
        atc.append(calculate_test_cons_for_orth_subset(aos[i],weights,degree))
    return atc
atc=calculate_all_test_cons_for_all_orth_subset(aos,weights,degree)
print(weights)
def test_orth_subset_1(orth_subset,weights,degree,test_cons):
    n=len(orth_subset)
    result_list=[]
    for j in range(0,len(weights)):
        result_list.append(0)
    for j in range(0,len(weights)):
        list1=[]
        for i in range(0,n):
            list1.append(0)
        b=0
        while (b==0):
            list1[n-1]=list1[n-1]+1
            for i in range(n-1,0,-1):
                if list1[i]>degree[j]:
                    list1[i]=0
                    list1[i-1]=list1[i-1]+1
            if list_sum(list1)==degree[j]:
                list2=list_mult_const(orth_subset[0],list1[0])
                for i in range(1,n):
                    list2=sum_list(list2,list_mult_const(orth_subset[i],list1[i]))
                if list2==weights[j]:
                    result_list[j]=result_list[j]+1
            if list1[0]==degree[j]:
                b=1
    if result_list==test_cons:
        res=1
    else:
        res=0
    return(res)
def test_all_orth_subset(aos,weights,degree,test_cons,test_coef):
    s=list_sum(test_cons)+test_coef
    result_list=[]
    for i in range(0,len(aos)):
        if len(aos[i])==s:
            if test_orth_subset(aos[i],weights,degree,test_cons)==1:
                if [2, 2, 3, 4, 3, 2, 1] not in aos[i]:
                    test_list=copy.copy(aos[i])
                    test_list.append([2, 2, 3, 4, 3, 2, 1])
                    tc=copy.copy(test_cons)
                    tc[0]=1
                    if test_orth_subset(test_list,weights,degree,tc)==1:
                        result_list.append(aos[i])
                else:
                    result_list.append(aos[i])
    return(result_list)
def test_all_orth_subset_1(aos,weights,degree,test_cons,test_coef):
    s=list_sum(test_cons)+test_coef
    result_list=[]
    for i in range(0,len(aos)):
        if len(aos[i])==s:
            if test_orth_subset_1(aos[i],weights,degree,test_cons)==1:
                if [2, 2, 3, 4, 3, 2, 1] not in aos[i]:
                    test_list=copy.copy(aos[i])
                    test_list.append([2, 2, 3, 4, 3, 2, 1])
                    tc=copy.copy(test_cons)
                    tc[0]=1
                    if test_orth_subset_1(test_list,weights,degree,tc)==1:
                        result_list.append(aos[i])
                else:
                    result_list.append(aos[i])

    return(result_list)
def test_all_orth_subset_all_test_cons(aos,weights,degree):
    test_cons=[]
    counter=1
    for i in range(0,len(weights)):
        test_cons.append(0)
    test_cons=[1, 1, 1, 1, 0, 0, 0]
    b=0
    print('CASE '+str(counter))
    print(test_cons)
    test_coef=1
    case_list=test_all_orth_subset(aos,weights,degree,test_cons,test_coef)
    while case_list==[] and list_sum(test_cons)+test_coef<=dim:
        test_coef=test_coef+1
        case_list=test_all_orth_subset(aos,weights,degree,
                                       test_cons,test_coef)
    print(case_list)
    while (b==0):
        test_cons[len(weights)-1]=test_cons[len(weights)-1]+1
        for i in range(len(weights)-1,0,-1):
            if test_cons[i]>1:
                test_cons[i]=0
                test_cons[i-1]=test_cons[i-1]+1
        if test_cons[0]>1:
            b=1
        if b!=1:
            counter=counter+1
            print(time.ctime())
            print('CASE '+str(counter))
            print(test_cons)
            test_coef=1
            case_list=test_all_orth_subset(aos,weights,degree,
                                           test_cons,test_coef)
            while case_list==[] and list_sum(test_cons)+test_coef<=dim:
                test_coef=test_coef+1
                case_list=test_all_orth_subset(aos,weights,degree,
                                               test_cons,test_coef)
            print(case_list)
def test_all_orth_subset_all_test_cons1(aos,atc):
    test_cons=[]
    counter=1
    for i in range(0,len(weights)):
        test_cons.append(0)
    b=0
    print('CASE '+str(counter))
    print(test_cons)
    case_list=[]
    test_coef=1
    for i in range(btc[list_sum(test_cons)+test_coef-1][0],
                   btc[list_sum(test_cons)+test_coef-1][1]):
        if atc[i]==test_cons:
            case_list.append(aos[i])
    print(case_list)
    while (b==0):
        test_cons[len(weights)-1]=test_cons[len(weights)-1]+1
        for i in range(len(weights)-1,0,-1):
            if test_cons[i]>1:
                test_cons[i]=0
                test_cons[i-1]=test_cons[i-1]+1
        if test_cons[0]>1:
            b=1
        if b!=1:
            counter=counter+1
            print(time.ctime())
            print('CASE '+str(counter))
            print(test_cons)
            case_list=[]
            test_coef=1
            if list_sum(test_cons)+test_coef<=dim:
                for i in range(btc[list_sum(test_cons)+test_coef-1][0],
                               btc[list_sum(test_cons)+test_coef-1][1]):
                    if atc[i]==test_cons:
                        case_list.append(aos[i])
            else:
                case_list=[]
            test_coef=test_coef+1
            while case_list==[] and list_sum(test_cons)+test_coef<=dim:
                for i in range(btc[list_sum(test_cons)+test_coef-1][0],
                               btc[list_sum(test_cons)+test_coef-1][1]):
                    if atc[i]==test_cons:
                        case_list.append(aos[i])
                test_coef=test_coef+1
            if case_list==[]:
                test_coef=0
                if list_sum(test_cons)+test_coef<=dim:
                    for i in range(btc[list_sum(test_cons)+test_coef-1][0],
                                   btc[list_sum(test_cons)+test_coef-1][1]):
                        if atc[i]==test_cons:
                            case_list.append(aos[i])
                else:
                    case_list=[]
                print(case_list)
                case_list=[]
            case_list1=[]
            for i in range(0,len(case_list)):
                if [2, 2, 3, 4, 3, 2, 1] not in case_list[i]:
                    test_list=copy.copy(case_list[i])
                    test_list.append([2, 2, 3, 4, 3, 2, 1])
                    tc=copy.copy(test_cons)
                    tc[0]=1
                    if test_orth_subset_1(test_list,weights,degree,tc)==1:
                        case_list1.append(case_list[i])
                else:
                    case_list1.append(case_list[i])
            print(case_list1)
print(time.ctime())
test_all_orth_subset_all_test_cons1(aos,atc)
print(time.ctime())

\end{verbatim}

Code for $E_8$:

\begin{verbatim}
import time
import copy
print(time.ctime())
positive_roots=['01000000', '01010000', '01011000', '01011100', '01011110',
                '01011111', '23354321', '01110000', '01111000', '01111100',
                '01111110', '01111111', '23454321', '01121000', '01121100',
                '01121110', '01121111', '23464321', '01122100', '01122110',
                '01122111', '23465321', '01122210', '01122211', '23465421',
                '01122221', '23465431', '23465432', '00100000', '00110000',
                '00111000', '00111100', '00111110', '00111111', '22454321',
                '00010000', '00011000', '00011100', '00011110', '00011111',
                '22354321', '00001000', '00001100', '00001110', '00001111',
                '22344321', '00000100', '00000110', '00000111', '22343321',
                '00000010', '00000011', '22343221', '00000001', '22343211',
                '22343210', '13354321', '12232100', '12232110', '12232111',
                '12232210', '12232211', '12232221', '11110000', '12233210',
                '12233211', '12233221', '11111000', '12233321', '11111100',
                '11111110', '11111111', '12243210', '12243211', '12243221',
                '11121000', '12243321', '11121100', '11121110', '11121111',
                '12244321', '11122100', '11122110', '11122111', '11122210',
                '11122211', '11122221', '10000000', '12343210', '12343211',
                '12343221', '11221000', '12343321', '11221100', '11221110',
                '11221111', '12344321', '11222100', '11222110', '11222111',
                '11222210', '11222211', '11222221', '10100000', '12354321',
                '11232100', '11232110', '11232111', '11232210', '11232211',
                '11232221', '10110000', '11233210', '11233211', '11233221',
                '10111000', '11233321', '10111100', '10111110', '10111111']
for i in range(0,len(positive_roots)):
    if positive_roots[i]=='23465432':
        print(i)
dim=8
scalar_product_matrix=[[2, 0, -1, 0, 0, 0, 0, 0],[0, 2, 0, -1, 0, 0, 0, 0],
                       [-1, 0, 2, -1, 0, 0, 0, 0],[0, -1, -1, 2, -1, 0, 0, 0],
                       [0, 0, 0, -1, 2, -1, 0, 0],[0, 0, 0, 0, -1, 2, -1, 0],
                       [0, 0, 0, 0, 0, -1, 2, -1],[0, 0, 0, 0, 0, 0, -1, 2]]
def dot_product(alpha,beta):
    sum=0
    for i in range(0,len(alpha)):
        for j in range(0,len(beta)):
            sum=sum+int(alpha[i])*int(beta[j])*scalar_product_matrix[i][j]
    return(sum)
def dot_product_matrix():
    list_dot_products=[]
    n=len(positive_roots)
    for i in range(0,n):
        list_dot_products.append([])
        for j in range(0,n):
            list_dot_products[i].append(dot_product(positive_roots[i],
                                                    positive_roots[j]))
    return(list_dot_products)
dpm=dot_product_matrix()
print('matrix calculated')
print(time.ctime())
def test_orth(list_of_roots):
    n=len(list_of_roots)
    b=1
    for i in range(0,n-1):
        for j in range(i+1,n):
            if (dpm[list_of_roots[i]][list_of_roots[j]])!=0.0:
                b=0
    return(b)
def root_to_list(alpha):
    result_list=[]
    for i in range(0,len(alpha)):
        result_list.append(int(alpha[i]))
    return(result_list)
def positive_roots_list(positive_roots):
    result_list=[]
    for i in range(0,len(positive_roots)):
        result_list.append(root_to_list(positive_roots[i]))
    return(result_list)
prl=positive_roots_list(positive_roots)
def roots_to_list(list_of_roots):
    result_list=[]
    for i in range(0,len(list_of_roots)):
        result_list.append(prl[list_of_roots[i]])

    return(result_list)
def all_orth_subset():
    result_list=[]
    result_list_numbers=[[]]
    n=len(positive_roots)
    for i in range(0,n):
        result_list.append([root_to_list(positive_roots[i])])
        result_list_numbers[0].append([i])
    for i in range(1,dim+1):
        print(i)
        print(time.ctime())
        result_list_numbers.append([])
        for j in range(0,len(result_list_numbers[i-1])):
            test_list=copy.copy(result_list_numbers[i-1][j])
            test_list.append(result_list_numbers[i-1][j][i-1])
            while test_list[i]+1<n:
                test_list[i]=test_list[i]+1
                if test_orth(test_list)==1:
                    result_list_numbers[i].append(copy.copy(test_list))
                    result_list.append(roots_to_list(test_list))
    return(result_list)
aos=all_orth_subset()
print(len(aos))
print(time.ctime())
def sum_list(list1,list2):
    list3=[]
    l=len(list1)
    for i in range(0,l):
        list3.append(list1[i]+list2[i])
    return(list3)
def list_sum(list):
    s=0
    n=len(list)
    for i in range(0,n):
        s=s+list[i]
    return(s)
def list_mult_const(list1,const):
    list2=[]
    l=len(list1)
    for i in range(0,l):
        list2.append(const*list1[i])
    return(list2)
def list_weight_sum(list1,coef_list):
    result_list=list_mult_const(list1[0],coef_list[0])
    for i in range(1,len(list1)):
        result_list=sum_list(result_list,list_mult_const(list1[i],coef_list[i]))
    return(result_list)
cascade=[[2, 3, 4, 6, 5, 4, 3, 2], [2, 2, 3, 4, 3, 2, 1, 0],
         [0, 1, 1, 2, 2, 2, 1, 0], [0, 1, 1, 2, 1, 0, 0, 0],
         [0, 1, 0, 0, 0, 0, 0, 0], [0, 0, 1, 0, 0, 0, 0, 0],
         [0, 0, 0, 0, 1, 0, 0, 0], [0, 0, 0, 0, 0, 0, 1, 0]]
weights=[]
list_cascade_weights=[[1, 0, 0, 0, 0, 0, 0, 0], [1, 1, 0, 0, 0, 0, 0, 0],
                      [2, 1, 1, 0, 0, 0, 0, 0], [3, 2, 1, 1, 0, 0, 0, 0],
                      [3, 2, 1, 1, 1, 0, 0, 0], [4, 3, 1, 1, 0, 1, 0, 0],
                      [5, 3, 2, 1, 0, 0, 1, 0], [3, 1, 1, 0, 0, 0, 0, 1]]
for i in range(0,len(list_cascade_weights)):
    weights.append(list_weight_sum(cascade,list_cascade_weights[i]))
degree=[]
for i in range(0,len(weights)):
    degree.append(list_sum(list_cascade_weights[i]))
print(degree)
def calculate_test_cons_for_orth_subset(orth_subset,weights,degree):
    n=len(orth_subset)
    result_list=[]
    for j in range(0,len(weights)):
        result_list.append(0)
    for j in range(0,len(weights)):
        list1=[]
        for i in range(0,n):
            list1.append(0)
        b=0
        while (b==0):
            list1[n-1]=list1[n-1]+1
            for i in range(n-1,0,-1):
                if list1[i]>degree[j]:
                    list1[i]=0
                    list1[i-1]=list1[i-1]+1
            if list_sum(list1)==degree[j]:
                list2=list_mult_const(orth_subset[0],list1[0])
                for i in range(1,n):
                    list2=sum_list(list2,list_mult_const(orth_subset[i],list1[i]))
                if list2==weights[j]:
                    result_list[j]=result_list[j]+1
            if list1[0]==degree[j]:
                b=1
    return(result_list)
def gen_c_wave_n_k(k,n):
    result_list=[]
    list1=[]
    for i in range(0,k):
        list1.append(0)
    result_list.append(copy.copy(list1))
    b=0
    while (b==0):
        list1[k-1]=list1[k-1]+1
        for i in range(k-1,0,-1):
            if list1[i]>n-1:
                list1[i-1]=list1[i-1]+1
                for l in range(i,k):
                    list1[l]=list1[i-1]
        if list1[0]<=n-1:
            result_list.append(copy.copy(list1))
        if list1[0]>n-1:
            b=1
    return result_list
keys=[]
for i in range(0,dim+1):
    for j in range(0,len(degree)):
        keys.append((degree[j],i+1))
print(keys)
dict_of_c_wave=dict()
for i in range(0,len(keys)):
    dict_of_c_wave[keys[i]]=copy.copy(gen_c_wave_n_k(keys[i][0],keys[i][1]))
print('dict_calculated')
print(time.ctime())
def calculate_test_cons_for_orth_subset(orth_subset,weights,degree):
    n=len(orth_subset)
    result_list=[]
    for j in range(0,len(weights)):
        result_list.append(0)
    for j in range(0,len(weights)):
        list1=copy.copy(dict_of_c_wave[(degree[j],n)])
        for i in range(0,len(list1)):
            result_list1=[]
            for m in range(0,dim):
                sum2=0
                for l in range(0,degree[j]):
                    sum2=sum2+orth_subset[list1[i][l]][m]
                result_list1.append(sum2)
            if result_list1==weights[j]:
                #print(list1[i])
                result_list[j]=result_list[j]+1
    return result_list
def border_tc(aos):
    result_list=[[0,0]]
    length=0
    for i in range(0,len(aos)):
        length1=len(aos[i])-1
        if length1>length:
            result_list[length][1]=i
            result_list.append([i,0])
            length=length1
    result_list[length][1]=len(aos)
    return(result_list)
print(border_tc(aos))
btc=border_tc(aos)
print('start time')
print(time.ctime())
def calculate_all_test_cons_for_all_orth_subset(aos,weights,degree):
    atc=[]
    percent=0
    f = open('text.txt', 'w')
    for i in range(0,len(atc)):
        f.write(str(atc[i])+'\n')
    sp=len(atc)
    for i in range(sp,352875):
        percent1=int(100*(i+1)/len(aos))
        if percent1>percent:
            percent=percent1
            print(str(percent)+'%')
            print(time.ctime())
        result_list=copy.copy(calculate_test_cons_for_orth_subset(aos[i],
                                                                  weights,degree))
        atc.append(result_list)
        f.write(str(result_list)+'\n')
    f.close()
    return atc
atc=calculate_all_test_cons_for_all_orth_subset(aos,weights,degree)
print('finish time')
print(time.ctime())
def test_all_orth_subset_all_test_cons1(aos,atc):
    test_cons=[]
    counter=1
    for i in range(0,len(weights)):
        test_cons.append(0)
    b=0
    print('CASE '+str(counter))
    print(test_cons)
    case_list=[]
    test_coef=1
    for i in range(btc[list_sum(test_cons)+test_coef-1][0],
                   btc[list_sum(test_cons)+test_coef-1][1]):
        if atc[i]==test_cons:
            case_list.append(aos[i])
    print(case_list)
    while (b==0):
        test_cons[len(weights)-1]=test_cons[len(weights)-1]+1
        for i in range(len(weights)-1,0,-1):
            if test_cons[i]>1:
                test_cons[i]=0
                test_cons[i-1]=test_cons[i-1]+1
        if test_cons[0]>1:
            b=1
        if b!=1:
            counter=counter+1
            print(time.ctime())
            print('CASE '+str(counter))
            print(test_cons)
            case_list=[]
            test_coef=1
            if list_sum(test_cons)+test_coef<=dim:
                for i in range(btc[list_sum(test_cons)+test_coef-1][0],
                               btc[list_sum(test_cons)+test_coef-1][1]):
                    if atc[i]==test_cons:
                        case_list.append(aos[i])
            else:
                case_list=[]
            test_coef=test_coef+1
            while case_list==[] and list_sum(test_cons)+test_coef<=dim:
                for i in range(btc[list_sum(test_cons)+test_coef-1][0],
                               btc[list_sum(test_cons)+test_coef-1][1]):
                    if atc[i]==test_cons:
                        case_list.append(aos[i])
                test_coef=test_coef+1
            case_list1=[]
            for i in range(0,len(case_list)):
                if [2, 3, 4, 6, 5, 4, 3, 2] not in case_list[i]:
                    test_list=copy.copy(case_list[i])
                    test_list.append([2, 3, 4, 6, 5, 4, 3, 2])
                    tc=copy.copy(test_cons)
                    tc[0]=1
                    if calculate_test_cons_for_orth_subset(test_list,
                                                           weights,degree)==tc:
                        case_list1.append(case_list[i])
                else:
                    case_list1.append(case_list[i])
            print(case_list1)
test_all_orth_subset_all_test_cons1(aos,atc)

\end{verbatim}

\end{document}